\documentclass[11pt,reqno]{article}
\usepackage{amsmath,amssymb,amscd}
\usepackage{epsfig}
\usepackage{parskip} 

\newcommand\qed{{\unskip\nobreak\hfil\penalty50\hskip2em\vadjust{}
    \nobreak\hfil$\Box$\parfillskip=0pt\finalhyphendemerits=0\par}}

\newtheorem{theorem}{Theorem}[section]
\newtheorem{lemma}[theorem]{Lemma}

\newtheorem{proposition}[theorem]{Proposition}
\newtheorem{definition}[theorem]{Definition}
\newtheorem{cor}[theorem]{Corollary}
\newtheorem{prop}[theorem]{Proposition}
\newtheorem{claim}[theorem]{Claim}

\let\qqed=\qed

\let\qed=\QED

\numberwithin{equation}{section}

\renewcommand{\d}{\delta}
\newcommand{\D}{\Delta}
\newcommand{\e}{\varepsilon}
\newcommand{\g}{\gamma}

\renewcommand{\i}{\infty}

\newcommand{\macp}{\beta}
\newcommand{\micp}{\alpha}

\newcommand{\eps}{\epsilon}

\newcommand{\N}{\mathbb{N}}
\newcommand{\Z}{\mathbb{Z}}
\newcommand{\R}{\mathbb{R}}
\newcommand{\E}{\mathbb{E}}
\newcommand{\PP}{\mathbb{P}}
\newcommand{\tord}{\mathbb{T}^d}

\newcommand{\dd}{{\rm d}}
\newcommand{\surv}{{\mathsf{S}}}
\newcommand{\coup}{{\mathcal{C}}}

\newcommand{\collision}{{\mathfrak{C}}}

\newcommand{\dist}{\vert \vert}
\newcommand{\ald}{\tfrac{\micp(n,m)}{d(n) + d(m)}}

\newcommand{\almeq}{\simeq_n}
\newcommand{\freem}{\mathfrak{A}_F}
\newcommand{\collop}{\mathfrak{A}_C}
\newcommand{\markovg}{\mathfrak{M}}

\newcommand{\empmeas}{\mu}

\newcommand{\tgel}{t_{{\rm gel}}}

\newcommand{\ball}{\mathsf{B}}

\def\xtp{x_1}
\def\mtp{m_1}

\def\mass{n}

\def\FF{\mathcal{F}}

\def\ue{u^\eps}

\begin{document}

\title{Coagulation and diffusion: a probabilistic perspective on the Smoluchowski PDE}
\author{Alan Hammond
\thanks{Departments of Mathematics and Statistics, U.C. Berkeley. 
This survey developed from a graduate course given at the University of Geneva in the autumn of 2012, at which time, the author worked in the Department of Statistics at the University of Oxford. 
The  course was supported by the Swiss Doctoral Program in Mathematics. The author was at that time supported principally by  EPSRC grant EP/I004378/1. He is now supported by NSF grant DMS-1512908.
}}

\maketitle

\begin{abstract}
The Smoluchowski coagulation-diffusion PDE is a system of partial differential equations modelling the evolution in time of mass-bearing Brownian particles which are subject to short-range pairwise coagulation. 
This survey presents a fairly detailed exposition of the kinetic limit derivation of the Smoluchowski PDE from a microscopic model of many coagulating Brownian particles that was undertaken in~\cite{HR3d}. It presents heuristic explanations of the form of the main theorem before discussing the proof, and presents key estimates in that proof using a novel probabilistic technique. The survey's principal aim is an exposition of this kinetic limit derivation, but it also contains an overview of several topics which either motivate or are motivated by this derivation. 
\end{abstract}

\newpage





\tableofcontents

\section{Introduction}
\subsection{Microscopic particles and macroscopic descriptions}

An important aim in statistical mechanics is to explain how the huge amount of information available in a microscopic description of a physical object, such as the positions and momenta of all the molecules comprising the air in a room, may be accurately summarised by first specifying a small number of physical parameters which are functions of macroscopic location, such as the density, temperature and pressure of this body of air at different points in the room, and then determining how these parameters evolve in space and time. 

\subsubsection{The elastic billiards model and the heat equation}

\begin{figure}
    \begin{center}
      \includegraphics[width=0.6\textwidth]{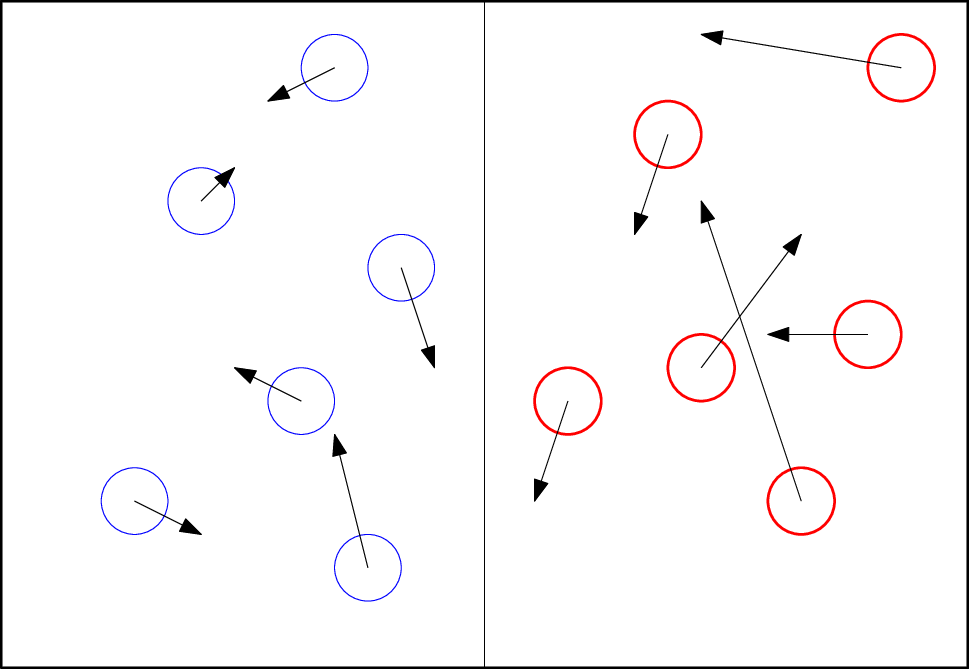}
    \end{center}
    \caption{This schematic figure depicts the left chamber $[-1,0] \times [-1,1]^2$ and the right chamber $[0,1] \times [-1,1]^2$ and the particles they contain at a negative moment at time at which the separating partition remains in place. The arrows indicate present velocities. The greater average magnitude of these velocities in the right chamber reflects the higher temperature of the body of air enclosed in that chamber.}\label{f.twochambers}
  \end{figure}

The microscopic system may begin out of equilbrium: for example, a still body of warm air in one room may be separated by a partition from another still body of cooler air in another, and then the partition instantaneously removed, so that air molecules from one side and the other intermingle over time, and an equilibrium is eventually approached in which the body of air in the whole room is again close to still, at a temperature which is some average of those of the two isolated systems at the original time. 
In such a case as this, it is a natural task to seek to summarise the evolution of a few suitable macroscopic physical quantities as the solution of partial differential equations. In the example, our object of study might be the temperature of the gas, and our aim to show that it is the heat equation, $\tfrac{\partial}{\partial t} T(x,t) =  \Delta T (x,t)$,  which models the macroscopic evolution $T(x,t)$ (with $x$ varying over the whole room, $[-1,1]^3$, say) of the temperature from the moment of the removal of the partition at time $t =0$ until a late time, $t \to \infty$, at which a new equilibrium is approached. In an idealized and very classical choice of microscopic description of the gas, we might model the ensemble of air molecules as a system of tiny spheres of equal radius and mass, each moving according to some velocity, and each pair of which undergoes a perfectly elastic collision on contact, in the same manner that a pair of billiards would. On each of the walls that comprise the boundary $\partial [-1,1]^3$ of the room, each sphere bounces elastically. The partition is modelled by the immobile sheet $\{ 0 \} \times [-1,1]^2$  on which spheres on either side also bounce elastically before time zero; the partition is removed instantaneously at that time. The initial instant of time may be taken to be zero, or some negative time. At that moment, we may scatter the spheres in an independent Poissonian manner throughout the room $[-1,1]^3$ (the reader may notice that in fact some extra rule is needed here to ensure the spheres' disjointness); and on one side and the other of the partition, choose their velocities independently, those to the right of the partition according to a non-degenerate law of zero mean, and those on the left according to another such law of lower variance than the first; in this way, we model two bodies of still air, a warm one in the right chamber $[0,1] \times [-1,1]^2$, and a cooler one in the left $[-1,0] \times [-1,1]^2$: see Figure~\ref{f.twochambers}.

In the microscopic model, there are huge numbers of tiny 
spheres in the system. Indeed, we may seek to understand the macroscopic evolution of temperature by in fact considering a whole sequence of microscopic models indexed by total particle number $N$, in a limit of high $N$. In the $N\textsuperscript{th}$ model, spheres are initially scattered as we described, with a Poissonian intensity $N$ thoughout $[-1,1]^3$. To carry out this task of understanding the large-scale evolution, we would wish to specify a microscopic definition of the notion of temperature, and then explain how it is in the high $N$ limit that the microscopic temperature data may be meaningfully reduced to a macroscopic description, and that this latter description indeed evolves according to the heat equation.  Microscopically, temperature is interpreted \cite[Section 1.1]{villanireview} as the average kinetic energy of particles, where here the velocity of particles is measured relative to the average velocity of nearby particles. Since our particle systems are large microscopically, considering as we do a 
high $N$ limit, we may specify in our $N\textsuperscript{th}$ microscopic model a definition of temperature at any given location $x \in [-1,1]^3$ as follows: first we may compute the mean velocity $v_{N,\delta}(x)$ of the set of spheres whose centres lie within some small distance $\delta$ of a given location $x$ in the room, and then we are able to define the microscopic temperature $T_{N,\delta}(x)$ to be the average of the square of the particle velocity minus $v_{N,\delta(x)}$, where the average is taken over the same set of spheres. Of course, the value $T_{N,\delta}(x)$ will change in time. As $N$ approaches infinity with $\delta$ being fixed but small, huge numbers of particles are involved in the empirical counts used for averaging. Our aim is to consider the space-time evolution of the microscopically specified temperatures after the high $N$ limit is taken, at which point, the weak law of large numbers might suggest that these empirical counts behave non-randomly to first order, so that our description becomes 
deterministic: the limit $T_{N,\delta}(x)$ will be some non-random function $T_\delta(x)$. In fact, since $\delta$ is fixed, we should not yet expect our system to approximate the heat equation, since there is an effect of macroscopic smearing in our calculation of microscopic temperature. Rather, one might expect the heat equation description to emerge if we take a $\delta \searrow 0$ limit of $T_\delta(x)$, after the first high $N$ limit has been taken. Moreover, to hope to obtain this description, we will also need to scale time appropriately in the $N\textsuperscript{th}$ microscopic model, as we take the first, high $N$, limit. In the scaled time coordinates, the microscopic models should make their approach to the new thermal equilibrium at the same rate, as $N \to \infty$. What rate this is in fact depends on another important consideration concerning the microscopic models which our brief description left unspecified: the radius $r_N$ of each sphere in the $N\textsuperscript{th}$ model must certainly be chosen to satisfy $r_N^3 \leq 
c N^{-1}$ for some constant $c > 0$, if only to permit all of the spheres to inhabit the room disjointly; our choice of decay rate for $r_N$ as a function of $N$, subject to this constraint, will determine the factor by which we scale time in the $N\textsuperscript{th}$ model in order to seek a heat equation description in the large.         

To implement the programme proposed in the preceding paragraph is an open problem, and in all likelihood, an extremely difficult one. There is no randomness in the model except in the initial selection of particle locations and velocities: from that time on, the deterministic laws of Newtonian mechanics govern the evolution of the microscopic models. Moreover, some choices for density admissible in the above description -- such as when $N r_N^3$ converges in high $N$, to a suitably small constant -- lead to rather dense systems of particles. The derivation may be less inordinately hard were more dilute choices of limit considered, where $r_N$ converges to zero more, and perhaps much more, quickly than does  $N^{-1/3}$.

It is important to note, however, that, if a choice of $r_N$ as a function of $N$ is made which is too rapidly decaying, we may leave the realm in which the heat equation is the appropriate macroscopic description. For example, if $r_N = o(N^{-1/2})$, it is a simple matter to check that a typical sphere after time zero will traverse the entire room on many occasions before meeting any other particle. The system will reach equilibrium after the removal of the partition simply by the free motion of the particles. The heat equation is only a suitable description when a typical particle experiences the thermal agitation caused by its collision with many other particles in short periods of macroscopic time. 

\subsubsection{The elastic billiards model and Boltzmann's equation}\label{s.elastic}

Moreover, the elastic billiards model crosses at least one interesting regime as it is diluted from the dense $r_n = \Theta \big(N^{-1/3}\big)$ phase towards the trivial free motion phase $r_N = o\big(N^{-1/2}\big)$. Consider the choice $r_N = N^{-1/2}$. A moment's thought shows that, in this regime, a typical sphere will travel (at unit-order velocity) for a duration before its first collision with another particle which on average neither tends to zero nor to infinity as $N \to \infty$. This is the regime of {\em constant mean free path}. The heat equation will not offer a suitable description for the evolution of temperature in this regime, because the mechanism providing for thermal agitation of particles -- manifest only when a typical particle has suffered many collisions -- occurs on a time scale which is marginally too slow. However, the programme of deriving a macroscopic description by means of a PDE does make sense, and in this case, offers a powerful model of gas dynamics. Suppose that, instead of using the microscopic data to form a description of temperature, we use 
it to describe the density of particles having a given velocity $v \in \R^d$ nearby a given location $x \in [-1,1]^3$. Particles may be scattered in a Poissonian fashion as before at the initial time, but with inhomogeneities in the intensity of this scattering permitted in both the space and velocity variables. With the macroscopic smearing parameter now being used to approximate velocity $v \in \R^d$ as well as location $x \in [-1,1]^3$, we may record a microscopic description $f^\delta_N(x,v)$ for the $\delta$-smeared density of spheres at space-velocity location $(x,v)$. Taking a high $N$ and then low $\delta$ limit as above, our macroscopic evolution is modelled by the fundamental system of equations in gas dynamics, Boltzmann's equation, valid for $t \geq 0$, $x \in \R^d$ and $v \in \R^d$:
\begin{equation}\label{e.boltzmann}
  \tfrac{\partial}{\partial t} f(x,v) = - v \cdot \nabla_x f (x,v) + Q(f,f) \, .
\end{equation}

Here, $- v \cdot \nabla_x$ is the free motion operator associated to particles of velocity $v$, while $Q(f,f)$ is a binary collision operator that reflects the microscopic elastic collision and whose form we will specify when  we return to Boltzmann's equation in a brief discussion in Section~\ref{s.neighbouring}. 
For now, note that the time evolution of the macroscopic densities is governed both by the free motion and by the collision operator. This is what is to be expected in the regime of constant mean free path, where the typical particle experiences unit-order durations free of collision and other such periods where several collisions occur.   

Boltzmann carried out a derivation of~(\ref{e.boltzmann}) as a model of gas dynamics in 1872, based on several assumptions, including one of molecule chaos that he called the Stosszahlansatz and which we will later discuss. (See \cite{boltzmann} for an English translation of his 1872 article.) The validity of his derivation was a matter of controversy, not least due to Loschmidt's paradox concerning precollisional particle independence (see Subsection~\ref{s.loschmidt}), and it was a fundamental advance made in 1975 by Lanford~\cite{Lanford} when the programme of rigorously deriving Boltzmann's equation from the elastic billiards model in the regime of mean free path was successfully implemented, for a short initial duration of time. By the latter condition, we mean that the validity of the description was established for some non-zero finite period, whose value depends on the form of the initial density profile of particles in space-velocity.  

Lanford derived Boltzmann's equation by establishing that the correlation functions concerning several particles in the model satisfy a hierarchy of equations called the BBGKY hierarchy, where the index of an equation in the hierarchy is the number of particles whose correlation is being considered, and by showing that when the correlation functions adhere to the BBGKY hierarchy, the density profile follows Boltzmann's equation. Illner and Pulvirenti implemented this approach in~\cite{IllnerPulvirenti} in order to derive Boltzmann's equation in a similar sense, but now globally in time, although with a comparable smallness condition, now on sparseness of the initial particle distribution; the cited derivation concerns a two dimensional gas, but this restriction on dimension was later lifted by the same authors. 

\subsubsection{Our main goal: coagulating Brownian particles and the Smoluchowski PDE}

This survey is intended to offer a detailed overview of a programme for deriving the macroscopic description of a gas of particles in the same vein as the descriptions above propose. However, our microscopic particles will diffuse, each following a Brownian trajectory, and as such their evolution is random, not deterministic; the mechanism of interaction will be pairwise as above, but a coagulation in which only one particle survives rather than a collision in which both do. On the other hand, in an effort to provide some generality in the microscopic description and richness in the macroscopic one, each of the particles will bear a mass, which the pairwise coagulation will conserve; and, moreover, we will permit the diffusivity of the Brownian trajectory of each particle to depend on the particle's mass. 

The partial differential equation which the programme seeks to obtain in this case -- the analogue of the heat equation or Boltzmann's equation in our opening examples -- is, like Boltzmann's equation, in fact a system of PDE, in our case coupled in the mass parameter, known as the Smoluchowski coagulation-diffusion PDE. 
The choice made for diluteness in the high particle number limit will be that of the regime of constant mean free path. The programme of deriving the PDE in the case of constant mean free path is sometimes called a {\em kinetic} limit derivation.  

In the special case of mass-independent diffusion rates, the kinetic limit derivation was carried out in 1980 by Lang and Nguyen~\cite{Lang}, 
who followed the method of showing that the correlation functions between several particles are described by the BBGKY hierarchy which Lanford had employed. 

Introduced to the problem of generalizing Lang and Nguyen's derivation of the Smoluchowski PDE by James Norris, the author collaborated on it with Fraydoun Rezakhanlou. The principal aim of these notes is to give an informal but fairly detailed exposition of the kinetic limit derivation of the Smoluchowski PDE that was undertaken for dimension $d \geq 3$ in \cite{HR3d}. The treatment also first presents heuristic arguments with the aim that the reader may understand why the main theorem should be true before beginning a presentation of the proof of the theorem, and it also uses some novel probabilistic techniques to obtain key estimates used in the proof. The survey also touches on some related topics.  


\subsubsection{Acknowledgments} The author is very grateful to James Norris for introducing him to the topic of diffusive coagulating systems and for valuable discussions.  He thanks Fraydoun Rezakhanlou for comments and guidance regarding the article's structure and approach; he further thanks Omer Angel, Nathana\"el Berestycki, Pierre Germain and Alain-Sol Sznitman for useful discussions, Dan Erdmann-Pham and Soumendu Mukherjee for comments on a draft version of the article, and the participants of the graduate class in Geneva for their interest and enthusiasm.

\subsection{The Smoluchowski coagulation-diffusion PDE}\label{s.smolcoag}

We begin by recording the form of these equations and offering a brief explanation of the phenomenon that they may be expected to describe. 

Let the dimension $d \geq 2$ be given. A collection of functions $f_n:\R^d \times [0,\infty) \to [0,\infty)$, $n \in \N$, is a strong solution of the discrete Smoluchowski coagulation-diffusion PDE with initial data 
$h_n:\R^d \to [0,\infty)$, $n \in \N$, if, for each $n \in \N$ and $x \in \R^d$, $f_n(x,0) = h_n(x)$; and, for each $n \in \N$ and 
$(x,t) \in \R^d \times [0,\infty)$,
\begin{equation}\label{syspde}
\frac{\partial{f_n}}{\partial t}(x,t)  = d(n)\Delta f_n (x,t)  +
Q^n_1(f) (x,t) -
Q^n_2(f) (x,t) \, ,
\end{equation} 
where the Laplacian acts on the spatial variable $x \in \R^d$. 
The final two terms are interaction terms, a gain term given by
\begin{equation}\label{gainterm}
Q^n_1(f) (x,t)  = \tfrac{1}{2}  \sum_{m=1}^{n-1} \macp (m,n-m) f_m(x,t)
f_{n-m}(x,t) \, ,
\end{equation}
and a loss term by
\begin{equation}\label{lossterm}
    Q^n_2(f) =   f_n (x,t) \sum_{m=1}^{\infty} \macp (m,n) f_m (x,t) \, .
\end{equation}
(When $t=0$, the partial time derivative on the left-hand side in~(\ref{syspde}) is interpreted as a right derivative.)

Note that the equations have two sets of parameters: the diffusion rates $d:\N \to (0,\infty)$ and the coagulation propensities $\beta:\N^2 \to [0,\infty)$. The equations have a continuous counterpart, where the mass variable is now a positive real, and the above sums are replaced in an evident way by integrals, which we will not consider in this survey except in passing.

To interpret the solution, consider a large number of minute particles in space $\R^d$, each carrying an integer mass. In a similar manner to our opening discussion, the quantity $f_n(x,t)$ is interpreted as the density of particles of mass $n \in \N$ in the immediate vicinity of location $x \in \R^d$ at time $t \geq 0$. The form of the right-hand side~(\ref{syspde}) reflects the two dynamics for the particles: diffusive transport and binary coagulation. Particles of mass $m$ diffuse at rate $2d(m)$, so that such a particle's displacement is given by $B(2d(m)t)$, $t \geq 0$, where $B:[0,\infty) \to \R^d$ is a standard Brownian motion. (The factor of two appears because the infinitesimal generator of standard Brownian motion is a {\em one-half multiple} of the Laplacian; when we call $d(n)$ the diffusion rate, this is thus strictly speaking a misnomer.) When a pair of particles are {\it microscopically} close, they may collide, disappearing from the model, to be replaced by a newcomer, whose mass is the sum 
of the two exiting particles'. The coagulation 
gain term~(\ref{gainterm}) expresses the possible means by which a new particle of mass $n$ may appear in the immediate vicinity of location $x$ at time $t$: by the coagulation of some pair of particles of masses $(1,n-1)$, 
or $(2,n-2)$ ... 
or $(n-1,1)$. The product form  $f_m(x,t) f_{n-m}(x,t)$ in the interaction term reflects an assumption that the particles in the immediate vicinity of $x$ are well mixed, and the coefficient $\beta(m,n-m)$ models the tendency of particles at close range of pair-type $(m,n-m)$ to coagulate in the immediate future. In the loss term~(\ref{lossterm}), we see the means by which the density $f_n(x,t)$ may fall due to coagulation: a particle of mass $n$ may drop out of the count for this density due to coagulation with another particle, and that other particle may have any mass~$m \in \N$.  

Our aim in this survey is to explain how the system~(\ref{syspde}) may be derived in a kinetic limit from a collection of microscopic random models of diffusing mass-bearing particles that are liable to coagulate in pairs at close range. We now describe in precise terms the elements for this programme; for the case at hand, we are thus presenting an instance of the type of programme  which we hazily sketched in our two opening examples. First, we specify the sequence of microscopic models, including their initial particle distributions, as well as their dynamics: the free motion of individual particles, and the mechanism of pairwise coagulation at close range.
In the main body of the article, we discuss only the derivation made in dimension $d \geq 3$, which was undertaken in \cite{HR3d}. Thus $d \geq 3$ may be assumed, except on one occasion when we make a short comment about the case when~$d=2$.

\subsection{The microscopic models}\label{sec.micromodels}

The sequence of microscopic random models will be indexed by the total number $N$ of particles intially present, at time zero. 
The $N$-indexed model will be specified by a probability measure~$\PP_N$. It is a measure not only on initial particle locations and masses but also on particle dynamics throughout $[0,\infty)$.  

\medskip

\noindent{\bf Initial particle distribution under $\PP_N$.} The quantity $f_n(x,0) = h_n(x)$ may be interpreted as the density of particles of mass $n$ in a tiny neighbourhood of $x \in \R^d$ at time zero. Thus, $\int_{\R^d} h_n(x) \, \dd x$ is interpreted as being proportional to the total number of particles of mass $n$ and the constant $Z \in (0,\infty)$, which we define by $Z = \sum_{n \in \N} \int_{\R^d} h_n(x) \, \dd x$, as being proportional to the total number of initial particles.

We will index the time-zero particle set under $\PP_N$ by $[N] : = \{1,\cdots,N\}$; the initial mass and location of particle $i$ will be denoted by $\big(x_i(0),m_i(0)\big)$. Reflecting the above density interpretation, we choose $\big(x_i(0),m_i(0)\big)$ independently, so that $\big(x_i(0),m_i(0)\big)$ has density $Z^{-1} h_n(x)$ at $(x,n) \in \R^d \times \N$.\\

\noindent{\bf Notation for particle trajectories under $\PP_N$.} We wish to describe the subsequent evolution of each of the initial particles under $\PP_N$. The trajectory of the $i\textsuperscript{th}$ particle will be described by $\big( x_i,m_i \big): [0,\infty) \to \big( \R^d \times \N \big) \cup \{ c \}$, where here $c$ is an element called a {\em cemetery} state whose role, which we will shortly describe in precise terms, is to house particles that have disappeared from the model due to being on the wrong side of a pairwise collision.  

As such, at any given time $t \geq 0$, the particle configuration under $\PP_N$ is described by a map $[N] \to \big( \R^d \times \N \big) \cup \{ c \}$, where $i \in [N]$ maps to $\big( x_i(t), m_i(t) \big)$ (or to $c$). 

To define the Markov process $\PP_N$ precisely, we will specify its Markov generator, which acts on test functions $F: \Big( \big( \R^d \times \N \big) \cup \{ c \} \Big)^{[N]} \to \R$. The action will be comprised of two parts: free motion of individual particles, and pairwise collision. We discuss our choice of each of these in words before providing the definition of the Markov generator. \\


\noindent{\bf Free motion.} A particle of mass $n \in \N$ follows, independently of other particles, the trajectory $t \to B\big(2d(n)t\big)$ relative to its starting point, where $B$ is a standard Brownian motion on $\R^d$.\\

\noindent{\bf Pairwise collision.} 
Any two particles will be liable to collide when their locations differ by order~$\eps$. Here, $\eps = \eps_N$, the interaction range, is determined by $N$ in a manner that we explain shortly. We introduce a compactly supported smooth interaction kernel $V:\R^d \to [0,\infty)$ and a collection of microscopic interaction strengths $\alpha: \N^2 \to (0,\infty)$, and declare that, at time $t$, particles $i$ and $j$ collide at infinitesimal rate $\alpha(m_i,m_j) V_\eps(\cdot)$,
where we adopt the convention that $V_\eps(\cdot) = \eps^{-2} V \big( \cdot/\eps \big)$. 
The argument $\cdot/\eps$ for $V$ indeed entails that collision may occur only between particles whose locations differ by order $\eps$; the prefactor of $\eps^{-2}$ is introduced because, in dimension $d \geq 3$, once a pair of particles have approached to distance of order $\eps$, they are liable to remain at such a displacement for order $\eps^2$ of time, since their relative displacement evolves as a Brownian motion of rate $2\big(d(n) + d(m)\big)$; thus, the role of this prefactor is to ensure that the proportion of instances of particle pairs approaching into the interaction range that result in collision is of unit order, uniformly in $N$. The role of the factor $\alpha(m_i,m_j)$ is to control whether this proportion is close to one for a given particle mass pair (which would be ensured by choosing the value of $\alpha$ in question to be high) or closer to zero.\\

\noindent{\bf The precise mechanism of collision.} 
On collision of $(x_i,m_i)$ and $(x_j,m_j)$ at time $t$, each of the pair of particles disappears, to be replaced by a new particle of mass $m_i + m_j$ in the vicinity. As a matter of convenience for the ensuing proofs, the precise rule we pick for the appearance of the new particle is to choose its location to be $x_i$ or $x_j$, with probabilities $\tfrac{m_i}{m_i + m_j}$ and  $\tfrac{m_j}{m_i + m_j}$. This rule permits the interpretation that, when two particles collide, one survives the collision and the other perishes; the probability of survival is proportional to incoming particle mass; the particle surviving collects the mass of the perishing particle, and the perishing particle vanishes from space. 

In a formal device, the perishing particle's location and mass are each sent to the cemetery state~$c$, where they remain forever. As such, for each $i \in [N]$, the $i\textsuperscript{th}$ particle's trajectory is described by setting the vanishing time $v_i \in [0,\infty]$ equal to the first time at which particle~$i$ experiences a collision in which it perishes. The trajectory is then given by 
$\big( x_i,m_i \big) \to \R^d \times \N$ on $[0,v_i)$ and $\big( x_i,m_i \big) = c$ on $[v_i,\infty)$.

\noindent{\bf The Markov generator of the dynamics.} For any configuration  $q \in  \Big( \big( \R^d \times \N \big) \cup \{ c \} \Big)^{[N]}$, write $I_q$, the surviving particle set, for those $i \in [N]$ such that $(x_i,m_i)$ lies in $\R^d \times \N$ (rather than equalling $c$). 
Let $F: \Big( \big( \R^d \times \N \big) \cup \{ c \} \Big)^{[N]} \to \R$ be smooth (in each hyperplane given by specifying the $c$-valued coordinates of the argument of $F$).
Then the Markov generator $\markovg$ for $\PP_N$ is given as follows. 
For each $q \in  \Big( \big( \R^d \times \N \big) \cup \{ c \} \Big)^{[N]}$,
$\markovg F (q) = \freem F (q) + \collop F (q)$, with the free-motion operator being given by 
$$
\freem F (q) = \sum_{i \in I_q} d(m_i) \Delta_{x_i} F (q) \, ,
$$ 
where $\Delta_{x_i}$ is the $d$-dimensional Laplacian acting on $F$ viewed as a function of $x_i \in \R^d$; and, recalling that $V_\eps(\cdot) = \eps^{-2}V\big(\cdot/\eps\big)$, with the collision operator being given by
\begin{equation}\label{e.collisionoperator}
\collop F (q) = \tfrac{1}{2} \sum_{i,j \in I_q}  \alpha(m_i,m_j)  V_\eps \big( x_i - x_j \big)
 \Big[ \tfrac{m_i}{m_i +
    m_j} F \big( S^1_{i,j} q \big) + \tfrac{m_j}{m_i +
    m_j} F \big( S^2_{i,j} q \big) - F(q)  \Big] \, .
\end{equation}
Here, $S^1_{i,j}(q)$, the configuration adopted in the event that particle $i$ survives collision with particle~$j$, is given by
\begin{equation}
 S^1_{i,j}(q)(k) = 
 \left\{ \begin{aligned}
     q(k) &  \, \, \, \, \, \textrm{for $k \in [N] \setminus \{ i,j \}$,} \\
    \big( x_i, m_i + m_j \big)   & \, \, \, \, \, \textrm{for $k = i$,} \\
         c  & \, \, \, \, \, \textrm{for $k = j$,}  \end{aligned} \right.
\end{equation}
while $S^2_{i,j}(q)$ is given by the same formula with the roles of $i$ and $j$ being reversed.

(A point of notation deserves mention. When we write $\sum_{i,j \in I_q}$ in specifying $\collop F (q)$ in~(\ref{e.collisionoperator}), we are using a slightly imprecise notation which refers to a sum over {\em distinct} pairs of indices $(i,j)$ lying in $I_q^2$.
Since the pairs are ordered, each is counted twice. The factor of one-half outside the sum is introduced to cancel this double counting.  Thus the $(i,j)$-indexed pair of particles is coagulating at the desired rate $\alpha(m_i,m_j) V_\e(x_i-x_j)$ that we specified in the pairwise collision description. We will use such abusive notation for double or triple sums later, but will comment at potential moments of confusion.)

\subsection{The regime of constant mean free path and the choice of interaction range}\label{sec.meanfreepath}

It remains to specify how the interaction range $\eps$ is determined by total initial particle number $N$. This choice is made to be in the regime of constant mean free path:
for dimension $d \geq 3$, $\eps = \eps_N$ will be chosen to satisfy
\begin{equation}\label{e.intrange}
   N =  Z \eps^{2 - d} \, .
\end{equation}

(This formula breaks down when $d=2$, and this is the basic reason why the two-dimensional case differs. Recall that we are focussing on the case $d \geq 3$.)
To explain why the regime for the length of the free path given by~(\ref{e.intrange}) is suitable, note that, since diffusion and coagulation terms are each present in the Smoluchowski PDE~(\ref{syspde}), we expect that the evolution of a typical particle will be determined both by its free motion and its collision with other particles. It will neither diffuse without collision nor collide repeatedly before diffusing a macroscopic distance. 

The consideration that this regime be adopted forces the choice of scaling of $\eps$ as a function of $N$: picking a uniformly random particle index $i \in [N]$ at the outset, $\eps$ should be chosen so that the mean time to first collision of particle $i$ converges as $N \to \infty$ to some strictly positive and finite constant.

A heuristic argument explains why (\ref{e.intrange}) produces this outcome. We anticipate that, at any given time $t \geq 0$, a positive (although $t$-dependent) proportion of particles are surviving (rather than in the cemetery state). Assume that the surviving particles at time $t$ are distributed so that the location and mass of  each is chosen independently; the law of the location-mass statistic $(x,n)$ of any given particle is equal to $f_n(x,t)$ (normalized to make the integral of this density equal to one). In other words, we are assuming in a very strong sense that the density profile of particles under~$\PP_N$ mimics the solution of~(\ref{syspde}). 

\begin{figure}
    \begin{center}
      \includegraphics[width=0.6\textwidth]{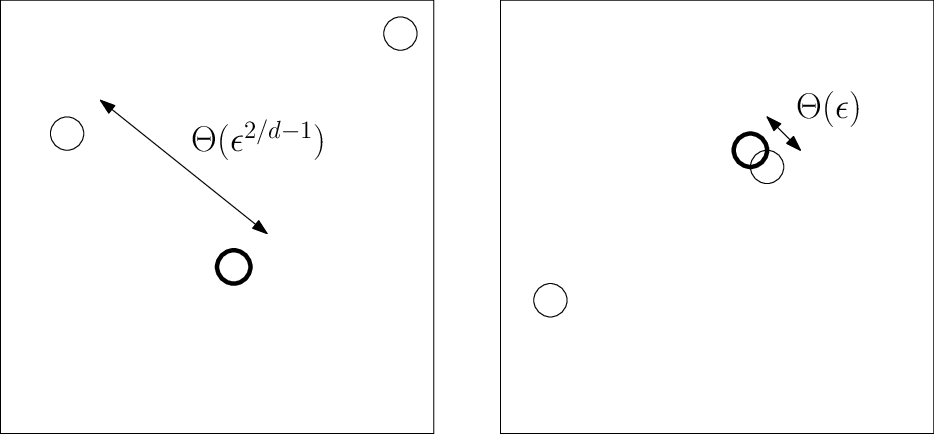}
    \end{center}
    \caption{The tracer particle is indicated by a bold circle in each sketch. {\em Left:} at a generic unit-order time, the order of distance $r$ of the nearest particle to the tracer particle may be expected to equal $r = \eps^{2/d - 1}$, since a ball of radius $r$ will contain $N r^d$ particles, and this choice of $r$ dictates a unit-order size for this quantity. {\em Right:} at a unit-order of special moments during a unit interval of time, this nearest distance drops to be less than~$\eps$, heralding a very short interaction window of duration $\Theta(\e^2)$ in which coagulation has a unit-order probability.}\label{f.twotimes}  
    \end{figure}

Pick a particle uniformly at random at the initial time and call the selected particle the {\em tracer} particle. We would like to estimate the mean number of collisions suffered by the tracer particle during $[0,1]$ in terms of $N$ and $\eps$. As we briefly discussed in the paragraph under the heading ``pairwise collision'' in the preceding section, this quantity is expected to have the same order as the number of other particles which enter the $\eps$-neighbourhood of the given particle during $[0,1]$. At any given time, our assumption on the distribution of other particles means that the probability that there is some other particle at distance less than~$\eps$ from the tracer particle is of order $N\eps^d$. Thus, the mean total amount of time during $[0,1]$  that some other particle is at distance less than~$\eps$ from the tracer particle is also of order $N \eps^d$. 
Whenever another particle approaches the tracer particle to distance $\eps$, it remains at the order of that distance for time of order $\eps^2$ (since $d \geq 3$). Thus, the mean number of different particles which during $[0,1]$ approach to within distance $\eps$ the tracer particle is of order  $N \eps^d \cdot \eps^{-2} = N \eps^{d-2}$. See Figure~\ref{f.twotimes}.

We thus see that imposing the relation~(\ref{e.intrange}) may indeed be expected to ensure that the mean number of collisions suffered by the tracer particle in unit time is bounded away from zero and infinity uniformly in $N$. 

\subsection{The recipe for the macroscopic coagulation propensities}\label{s.recipe}

The macroscopic coagulation propensities $\beta:\N^2 \to (0,\infty)$ that appear in the limiting system~(\ref{syspde}) depend in a non-trivial fashion on the microscopic parameters $V:\R^d \to [0,\infty)$, $\alpha:\N^2 \to (0,\infty)$ and $d:\N \to (0,\infty)$. 
Here is the recipe for obtaining $\beta(n,m)$ from these ingredients.
As we will later explain, there exists a unique solution $u = u_{n,m}:
\mathbb{R}^d \to (0,\infty)$ 
 of the equation
\begin{equation}\label{pdew}
   - \,  \Delta u_{n,m}(x)  \ = \  \frac{\micp(n,m)}{d(n) + d(m)} V(x) \Big[
1 - u_{n,m}
    (x) \Big]
\end{equation}
that satisfies 
$u_{n,m}(x) \to 0$ as $x \to \infty$. In fact, $0 \leq u_{n,m}(x) \leq 1$ for all $x \in \R^d$, and $u_{n,m}(x) =O\big(\dist x \dist^{2-d}\big)$ as $x \to \infty$.  (Here, as we will later, we write $\dist \cdot \dist$ for the Euclidean norm on $\R^d$.)

The quantities $\macp: \mathbb{N}^2  \to
(0,\infty)$ in (\ref{syspde}) are then specified by the formula
\begin{equation}\label{recp}
  \macp (n,m)= \micp(n,m) \int_{\mathbb{R}^d}
    V(x)\big(1 - u_{n,m} (x) \big) \, \dd x \, .
\end{equation}
We mention that the minus sign appearing on the left-hand side of (\ref{pdew}) was not used in the original treatment in~\cite{HR3d}. A positive choice for $u_{n,m}$ permits an attractive probabilistic interpretation of this quantity.
We wish to continue assembling the elements needed to state the main theorem concerning the kinetic derivation; when this is done, however, we will return to discuss the probabilitistic interpretation of~$u_{n,m}$: see Sections~\ref{secsimplecomp} and~\ref{s.killedbm}. It is tempting however before we continue to give a brief spoiler explaining the form of~(\ref{recp}): a fuller heuristic explanation will appear in Section~\ref{s.torus}.
In an instant of time beginning at a given moment $t$, the macroscopic rate 
of coagulation of pairs of particles of masses $n$ and $m$ near a given macroscopic location $z \in \R^d$ 
is equal to $\beta(n,m) f_n(z,t) f_m(z,t)$. To find a form for $\beta(n,m)$, note that this macroscopic rate 
should be computed as the integral over unit-order $x \in \R^d$ of an integrand given by the product of the coagulation rate associated to each pair of particles of such masses at negligible macroscopic distance from  $z$ that enjoy a relative displacement $\e x$,  and the density of presence of pairs of particles of these masses near $z$ and at such relative displacement.
In this `rate times density' description of the integrand, the rate should equal $\alpha(n,m) V(x)$ (up to a power of $\e$ that we neglect to mention here and in discussing the density term). Naively one might use a product structure ansatz to describe the density term, in which the role of $x$ would be irrelevant, and on this basis one would conclude that the density equals $f_n(z,t) f_m(z,t)$.
However, $x$ is of unit order, and in this case, particle pair presence is depleted due to the consideration that particles within the $\e$-radius interaction range may not be present because of the possiblity that they already coagulated in the last few instances of time (of duration of order~$\e^2$) leading up to time~$t$. In its accurate form, the density, which is  $\big(1 - u_{n,m}(x) \big)  f_n(z,t) f_m(z,t)$, contains an additional factor of $1 - u_{n,m}(x)$. As we will see in Sections~\ref{secsimplecomp} and~\ref{s.killedbm}, $u_{n,m}$ has an interpretation as a collision probability for a pair of Brownian particles. As such, the term  $1 - u_{n,m}(x)$ is interpreted as a survival probability for such particles: it is included to reflect the event that our nearby, $\e x$-displaced, particles survived their close encounter in the moments of time leading to the time~$t$ at which we consider the prospects of their imminent coagulation. 

\subsection{The weak formulation of the Smoluchowski PDE}

Pursuing the route to stating our main theorem, we now recast the Smoluchowski PDE~(\ref{syspde}) in weak form, since it is to this form of the equations that we will prove convergence. To do so, let $\mathfrak{J}$ be the space of sequences $J = \big\{ J_n : n \in \N \big\}$ of smooth compactly supported functions $J_n: \R^d \times [0,\infty) \to [0,\infty)$. 
Then we say that $f = \big\{ f_n: n \in \N \big\}$,
with $f_n: \R^d \times [0,\infty) \to [0,\infty)$ measurable for each $n \in \N$, is a weak solution of~(\ref{syspde}) if, for each $J \in \mathfrak{J}$, it satisfies the formula obtained from~(\ref{syspde}) by multiplication by $J_n$, integration in space-time, and integration by parts. Namely, such an $f$ solves~(\ref{syspde}) weakly if, for each $J \in \mathfrak{J}$ and $T \in (0,\infty)$,
\begin{eqnarray}
   & & \int_{\R^d} \Big( J_n(x,T) f_n(x,T) - J_n(x,0)  f_n(x,0) \Big) \, \dd x \label{syspdeweak} \\ 
   & = & 
    \int_{\R^d \times [0,T)}  \frac{\partial J_n(x,t)}{\partial t}  f_n(x,t) \, \dd x \dd t  \nonumber \\
   &  & \quad + \quad  
    \int_{\R^d \times [0,T)} \Big( d(n) f_n(x,t) \Delta J_n(x,t) 
      +  \big( Q_1^n(f)(x,t) -  Q_2^n(f)(x,t) \big)  J_n(x,t) \Big) \, \dd x \dd t \, . \nonumber
\end{eqnarray}

\subsection{Empirical densities}\label{s.empirical}

In our opening discussion of the programme for deriving a macroscopic limiting PDE, we suggested the use of $\delta$-macroscopically smeared particle counts as candidates to approximate the limiting evolution. Such counts play an important role in our derivation, and we will introduce them under the name {\em microscopic candidate densities} when we give an overview of the derivation of our main theorem, Theorem~\ref{thmo}, in Section~\ref{s.routetothm}.

However, to state this theorem, we will not use them. Rather, we will use a close cousin, empirical density measures defined under the microscopic models $\PP_N$. We now define these.

Under the law $\mathbb{P}_N$, let $\empmeas_N$ denote the $\PP_N$-random variable, valued in measures on space-mass-time $\mathbb{R}^d \times \N \times [0,\infty)$ such that, for each $t \geq 0$, its time-$t$ marginal $\empmeas_n(\cdot,t)$ is given by
$$
\empmeas_N(\cdot,t)
=\e^{d-2}\sum_{i\in I_{q(t)}}\d_{\big(x_i(t),m_i(t)\big)} \, .
$$

For given $n \in \N$, let $\empmeas_{N,n}$ denote the $\PP_N$-random variable, valued in measures on space-time $\mathbb{R}^d \times [0,\infty)$ such that, for each $t \geq 0$, its time-$t$ marginal $\empmeas_{N,n}(\cdot,t)$ is given by
$$
\empmeas_{N,n}(\cdot,t)
=\e^{d-2}\sum_{i\in I_{q(t)}}\d_{x_i(t)} \mathbf{1}_{x_i(t) = n} \, .
$$

Let $\mathcal{M}$ denote the space of measures $\mu$  on  $\mathbb{R}^d \times \N \times [0,\infty)$
such that $0 \leq \mu \left(\R^d \times \N \times [0,T]\right) \le TZ$, and 
note that $\empmeas_N$ is $\mathbb{P}_N$-a.s. valued in $\mathcal{M}$. (Recall that the constant $Z \in (0,\infty)$ was specified in Section~\ref{sec.micromodels}.)
We equip $\mathcal{M}$ with the topology of vague convergence,
under which a sequence $\{ \chi_n : n \in \N \}$
of measures converges to a limit $\chi$ precisely when $\int f \dd \chi_n$
converges to $\int f \dd \chi$ for all continuous $f: \mathbb{R}^d \times \N \times [0,\infty) \to \R$ of compact support. 
We make use of this topology because it makes $\mathcal{M}$ metrizable and compact.

\subsection{Hypotheses on microscopic parameters}

Our microscopic parameters are $d:\N \to (0,\infty)$, $\alpha:\N^2 \to [0,\infty)$ and $V:\R^d \to [0,\infty)$. 
In the original paper~\cite{HR3d} and in the detailed overview of proof that we give in this survey, some hypotheses on these parameters must be imposed to enable the derivation to be made. We make some comments about the hypotheses made in~\cite{HR3d} and then specify and discuss those we make here. The two sets of assumptions will be called the {\em original} and the {\em survey} assumptions throughout.

\subsubsection{Original assumptions}
The hypotheses governing the derivation in~\cite{HR3d} are now stated or at least roughly described. We will not follow the original derivation at a fine enough level of detail that the reasons for the form of these assumptions will be apparent; the survey assumptions deputise for the original ones in this regard. We do however summarise the original assumptions because they are significantly weaker than the survey ones. 

\noindent{\bf On the diffusion rate and the microscopic interaction strengths.} 
Suppose that there exists a function $\g : \mathbb{N}^2 \to 
(0,\infty)$ such that $\micp(n,m) \leq \g(n,m)$ for all $(n,m) \in \N^2$,
with $\g$ satisfying
\begin{equation}\label{hypo}
n_2 \cdot \g \big( n_1, n_2 + n_3 \big) \cdot \max  \left\{ 1 \, , \, \Big[ \frac{d (
n_2 + n_3 )}{d(n_2)}
\Big]^{2d-1} \right\}  \leq \big( n_2 + n_3 \big) \cdot \g(n_1,n_2) \, .
\end{equation}
(This condition is the same as that stated in equation~(1.9) of~\cite{HR3d}, but it has been simplified from its form in~\cite{HR3d}.)

\noindent{\bf On the initial condition.} A technical-to-state but fairly weak assumption is needed, of the membership in local $L^\infty$ space of some sums over $n$ of certain averages of $h_n$: see \cite[Section 1]{HR3d}. The assumption is certainly satisfied if $h_n$ is non-zero for only finitely many $n$, and each $h_n$ is compactly supported with bounded supremum.

\medskip

It is physically reasonable to think that the Brownian motion that is the free trajectory of the constituent particles in the models $\PP_N$ arises due to thermal agitation caused by many collisions with the constituents of an ambient environment of much smaller air molecules. Viewed in these terms, it is very natural to suppose that the diffusion rate $d(\cdot)$ will decrease as a function of the mass. Accepting this, the assumption~(\ref{hypo}) is rather weak. If the diffusion rate is indeed decreasing, then~(\ref{hypo}) is satisfied provided that there exists a function $C:\N \to \infty$ for which  $\micp(n,m)\le C(n)m$ for all $(n,m) \in \N$. 
Also, if the microscopic interaction strength $\micp$ is identically constant, then we may choose $\g$ equal to that constant in~(\ref{hypo}); if we then consider pure-power diffusion-rate choices $d(n) = n^\phi$, we find~(\ref{hypo}) to be satisfied whenever $\phi \in \R$ satisfies $\phi \leq (2d-1)^{-1}$. When $d=3$, then,  we are permitted choices of $d(n)$ that grow as quickly as $n^{1/5}$.


\subsubsection{Survey assumptions}


These assumptions are the following set of conditions. 

\medskip

\noindent{\bf On $V$.} The function $V:\R^d \to [0,\infty)$ is smooth and  compactly supported. 

\vspace{1mm}

\noindent{\bf On $h_n: \R^d\to [0,\infty)$.} 
Setting $\ell_n = \dist h_n \dist_{L^\infty(\R^d)}$ for $n \geq 1$,
each $\ell_n$ is finite, and the functions $h_n$ are all supported in a common given compact region of $\R^d$.

\vspace{1mm}

\noindent{\bf On $d:\N \to (0,\infty)$.} The function $d:\N \to (0,\infty)$ is non-increasing and $\sup_{m \in \N} m^{-1} d(m)^{- d/2} < \infty$.

\vspace{1mm}

 \noindent{\bf On $\alpha : \N^2  \to (0,\infty)$.} The supremum $\sup_{(n,m) \in \N^2} \alpha(n,m)$ is finite.

\vspace{1mm}

 \noindent{\bf A further condition, on $(h_n,d)$.} The sum $\sum_{n=1}^\infty \ell_n n d(n)^{d/2}$ is finite.

 Among these, the assumptions on the diffusion rates are genuinely restrictive: we must suppose that $d(m)$ grows less slowly than $m^{-2/d}$, which is not a particularly fast decay in any dimension $d \geq 3$. 
  No such imposition was made in the original assumptions.  It must also be admitted that the uniform bound demanded on $\alpha(\cdot,\cdot)$, is another  significant restriction. The final assumption limits the possibility for a heavy tail of high mass particles at the initial time, particularly since $d(n)$ must be supposed to decrease none too rapidly. Despite these limitations, the survey assumptions will permit us to offer a method of proof of key estimates needed for the main result which is largely self-contained, as well as being novel and very probabilistic in nature; since it serves our expository purpose, we have decided to accept the more limited domain of validity demanded by these assumptions.

\subsection{Statement of main theorem}

Here is our main result.

\begin{theorem}\label{thmo}
Let $d \geq 3$ and suppose that either of the above set of assumptions is in force. Let~$\mathcal P_N$ denote the law on $\mathcal{M}$ given by the law of the random measure $\empmeas_N$ under $\mathbb P_N$;  recall that $\epsilon$ is related to~$N$ by means of the
formula $N
\epsilon^{d-2} = Z$, with the constant $Z \in (0,\infty)$ being
given by the expression $Z = \sum_{n \in
\mathbb{N}}{\int_{\mathbb{R}^d}h_n}$.

Recall that the space of measures $\mathcal{M}$ has been given the topology  of vague convergence. The sequence~$\{\mathcal P_N : N \in \N \}$ is 
tight in~$\mathcal{M}$. Moreover, any law $\mathcal P$ on $\mathcal{M}$ that is   a weak limit point of the sequence $\{\mathcal P_N\}$ is concentrated on the space of measures taking the form $\sum_{n=0}^{\infty} f_n(x,t) \, \dd x \times \delta_n \times \dd t$ where 
$\{f_n: n \in \mathbb N\}$ is $\mathcal{P}$-almost surely a weak solution of (\ref{syspde}) that satisfies 
the initial condition $f_n(\cdot,0) = h_n(\cdot)$; recall that  the collection of
constants $\macp: \mathbb{N}^2 \to [0,\infty)$ is given by (\ref{recp}).
\end{theorem}  

The first assertion made by the theorem is trivial: since $\mathcal{M}$
is compact when equipped with the vague topology, any sequence of probability measures on~$\mathcal{M}$ is tight.

The reader may wonder what the meaning of the theorem is if it is not known that~(\ref{syspde}) has a weak (global in time) solution for the relevant parameter choices of $d(\cdot)$ and $\beta(\cdot,\cdot)$. In fact, the method of proof furnishes the existence of at least one weak solution. In any case, Lauren\c{c}ot and Mischler~\cite{LaurencotMischler} have established the existence of a global in time weak solution of~(\ref{syspde}) whenever $\lim_m m^{-1} \beta(n,m) = 0$, and $d(n) > 0$, for each $n \in \N$, conditions which are significantly weaker than those demanded by the theorem.

Theorem~\ref{thmo} describes the evolution of the density profiles of particles of various masses in the limit of large particle number by means of the Smoluchowski PDE, and in this way it realizes the derivation programme that we began this article by outlining, for the diffusive coagulating system in question.  The derivation has the merit of being global in time. However, note that, in general, there are limitations in the description offered of the large-scale behaviour of the system. If the weak solution of this system of PDE is not known to be unique, we merely demonstrate convergence in a subsequential sense to the space of solutions. For example, admitting the possibility that the system~(\ref{syspde}) has two distinct weak solutions $\big\{ f'_n: n \in \N \big\}$ and $\big\{ \hat{f}_n: n \in \N \big\}$ with initial condition $f_n(\cdot,0) = h_n(\cdot)$, each of the following behaviours is consistent with Theorem~\ref{thmo}:
\begin{itemize}
 \item the empirical densities under the microscopic models $\mathbb{P}_N$ may converge weakly to the solution $\big\{ f'_n: n \in \N \big\}$ as $N \to \infty$ along the subsequence of even integers, and to $\big\{ \hat{f}_n: n \in \N \big\}$ as $N \to \infty$ along the subsequence of odd integers;  
 \item it may be that evolution of these densities is accurately approximated by flipping a fair coin, with the densities converging weakly to  $\big\{ f'_n: n \in \N \big\}$ as $N \to \infty$ should the outcome be heads, and to   $\big\{ \hat{f}_n: n \in \N \big\}$ as $N \to \infty$ should the outcome be tails.
\end{itemize}

These peculiar scenarios are excluded if uniqueness of solutions to~(\ref{syspde}) is known. Some conditions for uniqueness are furnished by \cite[Proposition 2.6]{wrzd}; after deriving the kinetic limit of the PDE in 
\cite{HR3d}, Fraydoun Rezakhanlou and the author in  \cite{momentbounds} provided uniqueness under rather weaker hypotheses. Indeed, as~\cite[Remark 1.2]{momentbounds} discusses, the next proposition is a consequence of Theorems 1.1, 1.2, 1.3 and 1.4 of~\cite{momentbounds}. 
\begin{prop}\label{propuniqueness}
Let the dimension satisfy $d \geq 1$. 
For $a,b > 0$ such that $a + b < 1$, and for positive constants $c_1$ and $c_2$, assume that 
$\macp(n,m)\le c_1(n^a+m^a)$ and $d(n)\ge c_2 n^{-b}$ for
all $n,m\in \N$. Also assume that $d:\N \to (0,\infty)$ is non-increasing. There exists $e > 0$ such that $\sum_nn^e \|h_n\|_{L^\infty(\R^d)}<\i$
 and $\|\sum_nn^e h_n\|_{L^1(\R^d)}<\i$ imply that~(\ref{syspde}) has a unique weak solution. 
\end{prop}

Note that the survey assumptions in fact imply the hypotheses of Proposition~\ref{propuniqueness}. This means that, in working with these assumptions, we automatically obtain the simpler statement of convergence available when uniqueness of the PDE system is known (and which we are about to state).

It is a simple corollary of Theorem~\ref{thmo} and Proposition~\ref{propuniqueness} that convergence to~(\ref{syspde}) in fact holds in the following stronger sense.

\begin{cor}\label{coro}
Let $d \geq 3$ and suppose that the original assumptions, and the assumptions of Proposition~\ref{propuniqueness}, are in force.
Let $J: \mathbb{R}^d 
\times [0,\infty) \to \mathbb{R}$ be a 
 bounded and continuous test function. Then, for each
$n \in \mathbb{N}$ and $T \in (0,\infty)$,
\begin{equation}\label{resu}
\limsup_{N \to \infty } \mathbb{E}_N  \bigg\vert \int_{[0,T)} \int_{\mathbb{R}^d}
J(x,t) \big( \empmeas_{N,n}(\dd x,t) - f_n (x,t) \dd x  \big) \, \dd t \, \bigg\vert = \ 0,
\end{equation}
where again $N
\epsilon^{d-2} = Z$, with $Z = \sum_{n \in
\mathbb{N}}{\int_{\mathbb{R}^d}h_n}$. In (\ref{resu}),
    $\{ f_n: \mathbb{R}^d \times [0,\infty) \to [0,\infty), n \in \mathbb{N} \}$
denotes the unique weak solution to the system of
partial differential equations (\ref{syspde}), with $\macp: \mathbb{N}^2 \to [0,\infty)$ again given by (\ref{recp}).
\end{cor}

\subsection{A simple computation about the collision of two particles}\label{secsimplecomp}

The basic mechanism of interaction in our model concerns a pair of particles. Here, we explain a brief computation concerning such a pair, which offers a probabilistic interpretation of the function $u_{n,m}:\R^d \to [0,1]$ used in the recipe~(\ref{recp}) for the macroscopic coagulation propensity $\beta$.

Suppose at a certain time, a particle of mass $n$ is located at $0$ and another, of mass $m$, is located at $x \eps$, where $x \in \R^d$. The pair are thus prone to interact shortly, in the next order $\eps^2$ of time. Note also that, assuming uniform and independent placement of other particles in a compact region (in order to make an inference which we may find plausible for the actual model $\PP_N$ at any given time), the typical distance from a particle to the set of other particles is of order $N^{-1/d} = \eps^{1-2/d}$, which is far greater than the $\eps$ distance between the two particles in question. This means that in discussing the possible upcoming collision of this particle pair, we may harmlessly remove all other particles from the model.

Left with a two particle problem, we set $u^\eps_{n,m}(x)$ equal to the probability of subsequent collision of the pair. We may now use Brownian scaling, zooming in by a factor of $\eps^{-1}$ and slowing down time by a factor of $\eps^{-2}$, to obtain a particle of mass $n$ at the origin, one of mass $m$ at $x$, with the trajectories $X_1,X_2:[0,\infty) \to \R^d$ being Brownian motions of speeds $d(n)$ and $d(m)$, and collision occurring at rate $\alpha(n,m) V(X_1 - X_2)$. That is, $u^\eps_{n,m}(x)$ is independent of $\eps > 0$, and we may take~$\eps = 1$.   

As our notation suggests, $u_{n,m}^1$ is nothing other than $u_{n,m}$ from (\ref{pdew}): 

\begin{lemma}\label{l.g}
If $d \geq 3$, then $u_{n,m}^1$ is the unique solution $u_{n,m}:\R^d \to [0,1]$ of (\ref{pdew}).
\end{lemma}

For occasional later use, we further define $u_{n,m}^{[t]} :\R^d \to [0,1]$ for each $t > 0$ to be the probability that the two particles specified in the definition of $u_{n,m}^1(x)$ collide during $[0,t)$. Thus, $u_{n,m}^{[\infty]} = u_{n,m}^1$. In the expository discussion in Section~\ref{s.torus} (though not for the proof of Theorem~\ref{thmo}), we will need the next result. 
\begin{lemma}\label{l.gfin}
Suppose that $d \geq 3$. Then  $\dist u_{n,m}^{[t]} - u_{n,m}^1  \dist_\infty \to 0$ as $t \to \infty$.
\end{lemma}

We present the proofs of these two lemmas by using a more general notation which we now present.

\subsection{Killed Brownian motion and the Feynman-Kac formula}\label{s.killedbm}

In our two particle problem after scaling, the displacement of the particles performs a Brownian motion at rate $2 \big( d(n) + d(m) \big)$ until a random collision time. Slowing time by a factor of $d(n) + d(m)$, this process is rate two Brownian motion {\em killed at rate $\tfrac{\alpha(n,m)}{d(n) + d(m)} V$} in a sense we now explain. 

Let $W:\R^d \to [0,\infty)$ denote a smooth and compactly supported function. Let $x \in \R^d$. Rate two Brownian motion in $\R^d$ begun at $x$ and killed at rate $W$ is the stochastic process $X$ that we now specify. The process $X$ maps $[0,\infty)$ into $\R^d \cup \{ c \}$ where, as before, $c$ is a formal cemetery state. To define $X$, let $B:[0,\infty) \to \R^d$ denote rate two Brownian motion with $B(0) = x$. (Thus, $B(t/2) - B(0)$ is standard Brownian motion.) 
Define  its interaction until time $t$, $I_t$, to be equal to $\int_0^t W \big( B(s) \big) \dd s$. Let  $E$ denote an independent exponential random variable of rate one, and set  the killing time $K_x \in [0,\infty]$ equal to $\inf \big\{ t \geq 0: I_t \geq E \big\}$, with the convention that $\inf \emptyset = \infty$. Then 
\begin{equation*}
 X(s) =
\begin{cases}
  B(s)  \textrm{  for $s < K_x$} \, ,\\
   c  \textrm{  for $s \geq K_x$}\,  .
\end{cases}
\end{equation*}

We say that killing occurs if $K_x < \infty$ and 
let $u_W:\R^d \to [0,\infty)$ be such that $u_W(x)$ is the  probability that killing occurs.

\begin{lemma}\label{l.s}
For $d \geq 3$, $u = u_W$ is a solution of the modified Poisson equation
\begin{equation}\label{e.u}
 - \Delta u (x) =  W(x) (1-u)(x) \, .
\end{equation}
satisfying $u \to 0$ as $x \to \infty$. 
\end{lemma}
\noindent{\bf Remark.} The solution is unique subject to $u \to 0$ as $x \to \infty$. In a formal sense, this is verified by observing that the difference $v$ of two solutions solves  $\Delta v = W v$ on $\R^d$ and then noting that 
 \begin{equation}\label{e.formal}
  \int_{\R^d} \vert \vert \nabla v \vert \vert^2 \, \dd x  = - \int_{\R^d} v \Delta v \, \dd x \, ,
 \end{equation}
 whose right-hand side is  $-\int_{\R^d} W v^2 \dd x$ and is thus at most zero. Hence,  $\int_{\R^d} \vert \vert \nabla v \vert \vert^2 \dd x = 0$ and so $\nabla v$ is identically zero on $\R^d$. We thus see that $v$ is a constant function, and, since $v \to 0$ as $\dist x \dist \to \infty$, $v$ is identically equal to $0$. This would prove uniqueness, except that~(\ref{e.formal}) is a formal identity; if we integrate instead over the Euclidean ball $\ball_R$ and take $R \to \infty$, then the boundary term in Green's theorem vanishes in the limit provided a decay condition such as $\dist x \dist^{d-1} v(x) \nabla v (x) \to 0$ uniformly as $x \to \infty$ obtains. It is easy enough to confirm that this is the case: indeed, from the form of the fundamental solution of Laplace's equation in~\cite[Subsection 2.2.1.a]{EvansPDE}, we find that 
 $$
 v(x) = - \, c_0 \int_{\R^d} v(y) \, \vert \vert x - y \vert\vert^{2-d} W(y) \, \dd y \, ,
 $$ 
 where $c_0 = c_0(d)$ equals $d^{-1}(d-2)^{-1} \omega_d^{-1}$, with $\omega_d$ being the volume of the Euclidean unit ball in~$\R^d$. We see then that, since $W$ has compact support, $v$ has a decay at infinity at least as fast as~$\vert \vert x \vert\vert^{2-d}$. Differentiating this formula for $v$, we similarly learn that $\nabla v$ decays as quickly as $\vert \vert x \vert\vert^{1-d}$. Thus,  $\dist x \dist^{d-1} v(x) \nabla v (x)$ has decay as fast as $\vert \vert x \vert\vert^{2-d}$. The reader may also consult \cite[Section 6]{HR3d} for a proof of existence and uniqueness of the solution of~(\ref{e.u}) (subject to $u \to 0$ as $x \to \infty$) that uses Fredholm theory and compactness arguments.
 
 \medskip

\noindent{\bf Proof of Lemma~\ref{l.s}.} Let $v:\R^d \times [0,\infty) \to  [0,1]$ given by $v(x,t) = \E \, e^{- \int_0^t W\big(x + B(s) \big) \dd s}$, where the mean is taken over trajectories of rate two Brownian motion $B:[0,\infty) \to \R^d$ begun at zero. The Feynman-Kac 
formula~\cite[Section III.19]{RogersWilliams} shows that $v$ satisfies the partial differential equation
\begin{equation}\label{e.feynmankac}
 \tfrac{\partial}{\partial t} v(x,t) = \Delta v(x,t) - W(x) v(x,t) 
\end{equation}
for $x \in \R^d$ and $t > 0$. 
Note that for any $s > 0$, 
$$
\big\vert v(x,t+s) - v(x,t) \big\vert \leq s \, \dist W \dist_\infty \sup_{r \in [t,t+s]} \PP \Big( x + B(r) \in {\rm supp}(W) \Big) \, ;
$$
as $t \to \infty$, this probability tends to zero uniformly in $x$, so that we find that 
 $\tfrac{\partial}{\partial t} v(x,t) \to 0$ as $t \to \infty$, uniformly in $x \in \R^d$.

Note then that  $1 - u_W(x)$, which is the probability that  Brownian motion $X$ begun at $x$ and killed at rate $W$ is never killed, is equal to~$v(x,\infty)$. That $u_W$ solves  $-\Delta u(x) = W(x)\big( 1 - u(x)\big)$ in a distributional sense follows by taking a high $t$ limit of~(\ref{e.feynmankac}), since $v_t$ converges to $0$ locally in $L^1$. Since $W$ is smooth, $u$ being in local $L^2$ implies that $\Delta u$ is also in this space; thus, $u$ is locally in~$H^2$. Iterating, we find that in fact $u \in C^\infty$, and so $u$ solves~(\ref{e.u}) in strong sense. \qed 

The reader may consult Section~5.2 of the graduate PDE text~\cite{EvansPDE}
for a discussion of Sobolev spaces including~$H^2$. 

   

\medskip


\noindent{\bf Proof of Lemma~\ref{l.g}.} Note that, by the spatial-temporal scaling satisfied by Brownian motion, $u_{n,m}^1$ equals 
$u_W$ where $W = \tfrac{\alpha(n,m)V}{d(n) + d(m)}$. Hence Lemma~\ref{l.s} and the remark that follows it yield the result.  \qed

\medskip 

\noindent{\bf Proof of Lemma~\ref{l.gfin}.} We have that $u_{n,m}^{[t]}(x)$ equals $1 - v(x,t)$, and $u_{n,m}^1(x)$ equals $1- v(x,\infty)$. Note that $\big\vert v(x,\infty) - v(x,t) \big\vert$ is at most the probability that Brownian motion begun at $x$ visits the support of $W$ after time $t$. With $\mu$ denoting $d$-dimensional Lebesgue measure, this probability is at most a constant multiple of $\mu \big( {\rm supp} W \big) \cdot t^{-d/2}$, independently of $x \in \R^d$. \qed

In summary of this section and the preceding one, we have exhibited $u_{n,m}:\R^d \to [0,\infty)$ from~(\ref{pdew}) as a collision probability for a pair of Brownian particles. As we mentioned after stating the macroscopic coagulation propensity~(\ref{recp}), we are thus able to interpret the factor $1 - u_{n,m}$ 
in the integrand in~(\ref{recp}) in a way that offers 
 a heuristical explanation of the form of~(\ref{recp}). As the ensuing guide makes clear, we will explain these heuristics more carefully in Section~\ref{s.torus}.

\subsection{A guide to the rest of the survey}

We have now set up the microscopic models $\PP_N$ and laid out the programme for deriving their macroscopic evolution, including our main result, Theorem~\ref{thmo}. Our principal goal is to explain at a rather high, though not complete, level of detail, the proof of this theorem. We pause from pursuing this goal to explore two other directions first, however. First, in Section~\ref{s.neighbouring}, we offer a glimpse of several topics which are tangentially related to this principal goal; this discussion is intended only to whet the reader's appetite for perhaps some of these topics and problems, and, for this reason as well as owing to limitations in the author's knowledge, it is brief and very inexhaustive. Second, in Section~\ref{s.torus}, we offer a leisurely heuristic overview of our kinetic limit derivation in a very simplified special case, of annihilating constant diffusivity Brownian particles on a torus with a translation invariant initial condition. The argument is not rigorous at each step here, 
and in its method it does not provide a template for the derivation of Theorem~\ref{thmo}; rather, its main goal is to provide an intuitive explanation for the form~(\ref{recp}) of the recipe for the macroscopic coagulation rates $\beta:\N^2 \to (0,\infty)$; (the explanation elaborates that offered after the statement of this equation).  
We thus hope that, at the end of Section~\ref{s.torus}, the reader will have a fuller understanding of why the statement of Theorem~\ref{thmo} is true, if not yet of how it may be proved.

We then return to the survey's principal goal.  Section~\ref{s.routetothm} explains how Theorem~\ref{thmo} will be proved,
and the reader whose main interest is to see this proof explained may turn directly to this section. 
Therein, we introduce $\delta$-smeared approximations to the particle densities defined in the microscopic models $\PP_N$, called microscopic candidate densities. We state a fundamental estimate, the Stosszahlansatz, which expresses total coagulation propensity in $\PP_N$ approximately in terms of integrated products of microscopic candidate densities.

The next Section~\ref{s.outlinestoss} describes the method of proof of the Stosszahlansatz. While so doing, it gives an alternative explanation to that of Section~\ref{s.torus} for the form (\ref{recp}) of $\beta:\N^2 \to (0,\infty)$.

The actual proof of the Stosszahlansatz is given in Section~\ref{s.esterr}. The proof relies on several estimates concerning various integrated sums of test functions over pairs and triples of particle indices. These bounds in turn are reduced to two key estimates (or more accurately two sets of such estimates). The first of the two are particle concentration bounds, which state that if we have $L^\infty$-control at the initial time for the joint behaviour of $k$-tuples of particles in the models $\PP_N$ (with $k \in \N$ fixed, such as $k$ equal to two or three), then this control propagates to all later times; it is here that the more restrictive aspect of the survey assumptions, on the decay rate of the diffusion rates $d: \N \to (0,\infty)$, is invoked. The second key estimate takes the form of bounds on killing probabilities $u_W$ (that we introduced in Section~\ref{s.killedbm}), which are uniform over $W:\R^d \to [0,\infty)$ with given compact support. The proofs for this second key estimate also appear in this section. 

Section~\ref{s.partcon} provides a proof of the first key estimate, the particle concentration bounds. The proof occupies several pages, but we hope that it is probabilistically interesting and intuitive.

Finally, in Section~\ref{s.final}, we provide a summary of those points in our derivation where some steps were skipped, and mention where these omissions are treated in the original derivation in~\cite{HR3d}. We also take this opportunity to explain how the proof in~\cite{HR3d} is obscure at a certain moment, and highlight how in the present paper we have endeavoured to structure the arguments to shed light on this obscurity.

\section{A short foray into some neighbouring topics}\label{s.neighbouring}

\subsection{The Smoluchowski coagulation ordinary differential equation}

The spatially homogeneous analogue of the PDE we study, the Smoluchowski coagulation equation, was, to the author's best knowledge, originally formulated in Smoulchowski's seminal work~\cite[equation (67)]{smol}, and has been an object of attention for theoretical probabilists for a long time. The equation is now an ODE; it has a discrete and continuous mass version, as does the PDE. The solution $f: \N \times [0,\infty) \to [0,\infty)$ of the discrete form of this ODE may be written
\begin{equation}\label{smolode}
 \tfrac{\dd}{\dd t} f(n,t) = \tfrac{1}{2} \sum_{m=1}^{n-1} K(m,n-m) f(m,t) f(n - m,t) \, - \, f(n,t) \sum_{m=1}^\infty K(m,n) f(m,t) \, .  
\end{equation}
In the continuous counterpart, the space of masses is now $[0,\infty)$ rather than $\N$, and the two sums are replaced by integrals with the natural ranges of integration. Existence and uniqueness results were obtained by McLeod \cite{mcleoddiscrete}, White \cite{white} and Ball and Carr \cite{BallCarr} in the discrete case, the latter also addressing the mass conservation of solutions. 
For such results in the continuous case, see McLeod \cite{mcleodcts}, and Norris \cite{norrisone} who also shows an example of non-uniquneness.

Aldous' 1999 survey \cite{aldoussurvey} discusses how special choices of the interaction kernel $K$ have interesting probabilistic interpretations, in terms of point processes and related constructions on the complete graph (with $K(x,y) = xy$, discrete), the uniform measure on large trees ($K(x,y) = x+y$, discrete), and on the continuum random tree ($K(x,y) = x+y$, continuous). 

It is natural to pose the question of whether the ODE may be derived from a sequence of random models with diverging initial particle number. In such a model, the analogue of the microscopic law~$\PP_N$ in the programme we discuss would consist of a probability measure under which $N$ particles carry integer or non-negative real valued masses, and any pair coalesces at infinitesimal rate $K(\cdot,\cdot)$, where the arguments are the masses of the concerned pair; at the moment of their coalescence, the two particles leave the system, to be replaced by a newcomer whose mass is the sum of the departing pair's. 
In \cite{norrisone,norristwo}, Norris derives the Smoluchowski ODE from this model,  considers more general mechanisms for coalescence and carries out corresponding derivations for them. Rezakhanlou in~\cite{Rezak13} presents sufficient conditions for gelation (a concept that we will shortly discuss in the spatial setting) in such models.

In the case where $K$ is identically constant, this random system of coalescing particles is called Kingman's coalescent \cite{kingmanone,kingmantwo}. Although we have defined it only with an initial condition with a finite number of particles, the system {\em comes down from infinity}, in the sense that there is a well defined stochastic process $S$ that maps $(0,\infty)$ to $\N$ with the property that,
for any $t > 0$, the process $[0,\infty) \to \N: s \to S(t+s)$ has the distribution of the process of the number of surviving particles in Kingman's coalescent given that the initial number is $S(t)$.

See \cite{Berestycki} for a recent survey of coalescence theory, including a treatment of Kingman's coalescents, a more general class of coalescents, $\Lambda$-coalescents, in which several particles may combine simultaneously, spatial models and their applications to population genetics. 

\subsection{Coalescing random walkers on $\Z^d$}

Suppose that at time zero, a finite collection of walkers are located, each at some site in $\Z^d$; there may be several walkers at any given site. Each walker performs a continuous time simple random walk, staying at her present site for a duration which is exponentially distributed, of mean one, independently of other decisions, and then jumping with equal probability to one of her $2d$ neighbours. For any instance of a pair of walkers occupying the same site at any given moment of time, one walker in the pair is annihilated at exponential rate one. The description is informal and we do not provide a precise formulation here.

\medskip

\noindent{\bf Many walkers initially: coming down from infinity}

Naturally enough, these models have much in common with coagulating diffusive systems. Defining this model on a singleton set, where all walkers must occupy the same site, note that the model reduces to Kingman's coalescent, which, begun with infinitely many particles, has at any positive time only finitely many. Does this phenomenon also take place in $\Z^d$? A natural starting condition for the model on $\Z^d$ is to begin with $N$ walkers, all located at the origin, and consider high $N$ behaviour. In \cite{ABL}, it is shown that, in contrast to the non-spatial case, infinitely many walkers survive asymptotically: if there are $N$ initially, there are of the order of $(\log^*N)^d$ at any given positive time, where the function $\log^*N$ in essence denotes the number of iterations of the logarithm which reduces the value of $N$ to below zero. At constant time, the surviving particles roughly fill up a ball of radius $\log^*N$ with a tight number of particles present at each site in this ball. Since the model 
asymptotically manufactures an arbitrarily large number of surviving particles, however slowly it does so in the high $N$ limit, it is a natural extension of the 
derivation of Theorem~\ref{thmo} to enquire as to whether a counterpart result holds for this model. We may expect to squeeze space by a factor of $\log^*N$ and slow time by the square of this quantity to reach the regime described by the PDE. The recipe~(\ref{recp}) for the macroscopic coagulation rates will presumably be altered so that a discrete Laplacian appears instead.   

\medskip

\noindent{\bf The model on a discrete torus: Smoluchowski and Kingman}

 Similarly to the model on $\Z^d$, we may consider the model of annihilating random walkers on the discrete torus 
 $\tord_N$, formed by  quotienting $\Z^d$ by $N \Z^d$. Imagine for the sake of simplicity that, at time zero, the walkers sit one apiece at each of the $N^d$ sites of $\tord_N$. In the first few units of time, a positive fraction of walkers are annihilated, and the system becomes sparser. Writing $s_N:[0,\infty) \to [0,1]$ for the proportion of surviving particles in $\tord_N$, we may define $s:[0,\infty) \to [0,1]$, $s = \lim_N s_N$, (an almost sure limit which clearly exists and is non-random), and ask questions about the rate at which $s(t) \searrow 0$ as unscaled time $t$ tends to infinity. 
 
 On the other hand, by scaling space and time, we may find a regime analogous to that identified in the principal aim of this survey. Scaling $\tord_N$ to the unit torus by squeezing space by a factor of $N$,  Brownian scaling dictates that time be sped up by a factor of $N^2$. In the new coordinates, there is a rapid obliteration of walkers that reduces their number from $N^d$ to an order of $N^{d-2}$ (by any given positive time). Note that the relation (\ref{e.intrange}), that $n \e^{d-2}$ is of unit order, with $n$ the survivng particle number, is satisfied in the sense that $n$ is order $N^{d-2}$ and $\eps$, the interaction range, is $N^{-1}$. In scaled coordinates then, the system crashes down from infinity, and naturally slows down into the regime of constant mean free path. Our translation invariant choice of initial condition should manifest itself in this regime by convergence of the surviving particle number to a solution of the Smoluchowski coagulation ODE~(\ref{smolode}); since 
particles are indistinguishable, the kernel $K$ is identically equal to one.    
 
 At time scales beyond the $N^2$ rescaling, surviving particles cross the torus many times between collisions. If we choose a speeding up of time by a factor of $N^d$, then we enter a regime where only finitely many particles survive. Indeed, Cox~\cite{Cox} proved that, when dimension $d \geq 3$, the rescaled process of surviving particle number, $[0,\infty) \to \N$, $t \to s_N(2 N^d t/G)$, converges in law to the process of surviving particle number in Kingman's coaelscent. The constant factor $G$ is the mean total amount of time that a continuous-time simple random walk in $\Z^d$ begun at the origin spends at the origin in all positive time. This regime is one where a finite population of walkers each mixes spatially at an infinite rate, thus becoming indistinguishable; the presence of the factor~$G$ is explained by noting that the time of any pair of such particles to meet should be gauged in a clock which advances at a speed which is double (since a displacement between two walks is considered) that at 
which the unit-rate walker on $\Z^d$ at late time encounters previously unvisited vertices. Cox also studied the problem in dimension two, and noted the implications of the solution for voter model consensus times.

\subsection{The elastic billiards model and Boltzmann's equation}

Our discussion draws heavily on Chapter 1 of Villani's review~\cite{villanireview} of collisional kinetic theory. 

\subsubsection{The form of the equations}

In Subsection~\ref{s.elastic}, we mentioned Lanford's 1975 derivation for short times of Boltzmann's equation~(\ref{e.boltzmann})
in a kinetic limit from a system of elastically colliding billiards. In the form of the equation suitable for such a billiard model in dimension $d=3$, the collision operator in~(\ref{e.boltzmann}) is given by 
\begin{equation}\label{e.boltzmanncoll}
 Q(f,f)(t,x,v) = C \int_{\R^d} \dd v_* \int_{S^{d-1}} \, \dd \sigma \,  \dist v - v_* \dist \Big( f(t,x,v') f(t,x,v_*') - f(t,x,v) f(t,x,v_*) \Big) \, ,
\end{equation}
where $C \in (0,\infty)$ is a constant. To explain the notation, consider a collision that a sphere of velocity $v \in \R^d$ may undergo. The particle with which it collides has some velocity $v_* \in \R^d$. Impact may occur over a hemisphere in a surface of the velocity $v$ sphere. Denoting the outgoing velocities of the two spheres by $v'$ and $v_*'$, conservation of momentum and kinetic  energy imply that
\begin{equation*}
\begin{cases}
    \, v' + v_*' = v + v_* \\
    \,  \dist v' \dist^2 + \dist v_*' \dist^2 = \dist v \dist^2 + \dist v_* \dist^2 \,  .
\end{cases}
\end{equation*}
The form of the outgoing velocities is determined by the angle of impact. The possibilities may be parameterized as follows:
\begin{equation*}
\begin{cases}
    \, v'  =  \tfrac{v + v_*}{2} +  \tfrac{\dist v + v_* \dist}{2} \sigma \\
    \,  v_*'  =  \tfrac{v + v_*}{2} -  \tfrac{\dist v + v_* \dist}{2} \sigma  \, , 
\end{cases}
\end{equation*}
as $\sigma$ varies over $S^{d-1}$. 

We mention in passing one notable feature of (\ref{e.boltzmanncoll}): the term $\dist v - v_* \dist$, which is the {\em Boltzmann collision kernel}, and which in a more general setting may depend non-trivially on $\sigma$, has no such dependence in the present case of elastic collisions and dimension $d = 3$. 

The collision operator $Q(f,f)$ may be written as a difference of non-negative gain and loss terms, $Q^+(f,f) - Q^-(f,f)$, by splitting~(\ref{e.boltzmanncoll}) across the minus sign in the big bracket. We obtain a phenomenological description of~(\ref{e.boltzmann}) akin to that offered for the Smoluchowski PDE appearing after the equations in Section~\ref{s.smolcoag}.
Particles of velocity $v$ near location $x$ are subject to collision, and contribute to the loss term $Q^-(f,f)$; collisions of particles with other velocity pairs may occur which produce new velocity $v$ particles near $x$, as the gain term $Q^+(f,f)$ records.

\subsubsection{Loschmidt's paradox and the Stosszahlansatz}\label{s.loschmidt}

A system of elastically colliding billiards in a box is at equilibrium a reversible continuous-time Markov chain. However,
Boltzmann's equation begun from generic initial data do not share this reversibility. Indeed, 
Boltzmann's $H$ theorem shows that the Boltzmann $H$ functional, 
$$
 H(f) = \int_{(x,v) \in \R^d \times \R^d} f \log f \, ,
$$
satisfies $\tfrac{\dd}{\dd t} H \big( f(t,\cdot,\cdot) \big) \leq 0$. The quantity $H$ may be viewed as a measure of information; information dissipates monotonically as time evolves, in accordance with the second law of thermodynamics. At some late time, this rate of dissipation may slow as the system approaches equilibrium. However, for generic initial data for Boltzmann's equation, $H$ may decrease in a strictly monotonic fashion as time advances. This irreversible property of the macroscopic evolution seems to be in tension with the reversible nature of the basic collision event that two billiards may undergo revealed by reversing in time a viewing of the collision. Concerns such as these caused Boltzmann's claim that the equation offered an accurate macroscopic description of classical many body systems such as elastic billiards to be treated with much scepticism. Loschmidt found a paradox which brought these concerns into a sharper relief. Accepting that a large system of elastic billiards is 
accurately modelled by Boltzmann's equation for all typical choices of initial data, begin with some such data and run the deterministic dynamical rules for the billiard system for some fixed time $t$. The density profiles will approximately follow the solution of Boltzmann's equation, and the $H$ functional will drop from its initial value. Stop the evolution at time $t$ and then reverse the velocity of each particle, leaving each particle's location unchanged. Then run the system for a further $t$ units of time. Clearly the resulting evolution will be a rerun of the dynamics we just witnessed in the sense of reversed time. At the end of this second dynamics, the collection of billiards has its original set of locations, with reversed velocities. Note however that during this second dynamics, the $H$ functional was rising, not falling. However, this is impossible for a system which is accurately approximating a solution of Boltzmann's equation. 

Loschmidt's paradox indicates that not all microscopic data consistent with a given macroscopic density profile may result in an evolution for which Boltzmann's equation is an accurate model. The velocity-reversed time-$t$ particle data is a counterexample to the hypothesis that Boltzmann's equation may be so derived from all such microscopic data. However, there is no contradiction to the hypothesis that all but a tiny fraction of particle configurations approximating a given density profile begin a dynamics whose evolution is accurately described by Boltzmann's equation. 

The paradox also has implications for methods of proof that may be proposed for deriving Boltzmann's equation from microscopic models. In Boltzmann's own derivation, he invoked an assumption of independence on the part of colliding particles, which he called the Stosszahlansatz, or the {\em collision number hypothesis}. This asserts roughly that, in the neighbourhood of a location $x$ at any time $t$, the distribution of the numbers of collisions of particles of two given velocities, $v$ and $v_*$, in a many body system of elastic billiards, is accurately specified by knowing the densities (macroscopically denoted by $f(x,v,t)$ and $f(x,v_*,t)$) of velocity $v$ and $v_*$ particles near $x$ at time $t$; the particles' histories until this time does not significantly bias the local structure of these families of particles away from that of a Poissonian system of such particles at these two densities; so that, for example, the rate of collision of a randomly picked particle of velocity $v$ near $x$ at time $t$ 
with a particle of velocity $v_*$, and the distribution of the impact parameter on collision, are accurately modelled by the Poisson systems at these densities. Loschmidt's paradox indicates a subtlety about the Stosszahlansatz. It may be valid for precollisional particle velocities, but in cannot be for postcollisional ones; for the latter, the history of the concerned particles has a great deal that biases their distribution from a Poissonian cloud model for the two velocity types. In other words, the mechanism of elastic collision may propagate chaos, taking independent randomness present at an initial time and preserving it at given later times, but the chaos propagated is one-sided, not double-sided, referring to statistical inferences about the particles' future, and not their past.

In our more humble setting of coagulating Brownian particles, a key role is played by a result, Proposition~\ref{propsz}, concerning collision propensity for the microscopic models, which we interpret as the Stosszahlansatz, as we will see in Section~\ref{s.routetothm}. However, the random and reversible nature of the free motion of the individual particles means that there is no analogue of Loschmidt's paradox and no need to formulate the Stosszahlansatz as a statement concerning merely one-sided, rather than double-sided, chaos.

\subsection{Gelation and mass conservation}

\subsubsection{Mass conservation for the Smoluchowski PDE}\label{s.masscons}

The collision event in the microscopic models $\PP_N$ conserves mass. How does mass conservation manifest itself macroscopically, in a solution of the Smoluchoski PDE? For a solution $\big\{ f_n: n \in \R^d \big\}$ of~(\ref{syspde}), we may intepret $M_n(t) : = n \int_{\R^d} f_n(x,t) \dd t$ as the total mass among particles of mass $n \in \N$ at time $t \in [0,\infty)$, and thus
$M(t) : = \sum_{n \in \N} M_n(t)$ to be the cumulative mass of particles at this time. For any $T > 0$, a solution of~(\ref{syspde}) is said to conserve mass during $[0,T]$ if $M(s) = M(0)$ for all $s \in [0,T]$. The passage from the microscopic to the macroscopic 
might lead one to expect solutions to be mass conserving on all such intervals. In fact, only the inference that $M:[0,\infty) \to [0,\infty)$ is non-increasing is readily available. We may define then the gelation time $\tgel \in [0,\infty]$, 
$\tgel = \inf \big\{ t \geq 0: M(t) < M(0) \big\}$, with $\inf \emptyset = \infty$. It is shown in \cite{momentbounds} that the unique weak solution of the PDE which 
Proposition~\ref{propuniqueness} provides is mass conserving in the sense that $\tgel = \infty$. Certainly under the survey assumptions, the resulting solution of the PDE satisfies the hypotheses of  Proposition~\ref{propuniqueness}, and so is mass conserving. Indeed, this is true in every circumstance under which a kinetic limit derivation of the Smoluchowski PDE has been carried out. 

\subsubsection{The meaning of gelation for the microscopic models}\label{s.gelation}

Nonetheless, it is natural to ask what behaviour we would expect to see under the laws $\PP_N$ in a system which converges to a solution of the Smoluchowski PDE with gelation. After the gelation time $\tgel$, a positive fraction of the initial particle mass in a high-$N$ indexed model $\PP_N$ will be present in particles above any given $K \in \N$; this fraction is independent of the value of $K$, although, for high values of $K$, we will have to increase $N$ in order to witness this effect. A gel, composed of super-heavy particles, is forming microscopically beyond the gelation time. 

Does this phenomenon actually take place in a model $\PP_N$ for some choice of parameters $V$, $d(\cdot)$ and $\alpha(\cdot,\cdot)$, or in some variant of this model? 

To prepare to answer this, we first discuss a natural extension to our definition of microscopic model. Under $\PP_N$, all particles have an equal interaction range $\eps = \eps_N$, irrespective of their mass. It is natural to introduce a mass-dependent interaction range, of the form $r_n \eps$ for particles of mass $n$; presumably $r_n$ would be increasing, and the choice $r_n = n^{1/d}$ would correspond to solid ball particles composed of a common material which instantaneously merge to form a larger such ball on collision. Other choices $r_n = n^\chi$, for $\chi \in [1/d,1]$, may be possible, corresponding to fractal geometries for the internal particle structure (as we will discuss in Section~\ref{s.onetwo}). Without the assumption of additional and non-local attractive inter-particle forces, it is hard to see, however, how a choice of the form $\chi > 1$ would be physically meaningful. The choice $\chi = 1$ is already a little beyond the border of the plausible realm:  in a farfetched effort to justify 
this choice, we may model each particle as a long and very thin bar, and imagine that each bar rotates rapidly and chaotically about its centre of mass, while diffusing on a slower time-scale; when two bars touch, they instantaneously and rigidly join. Because of their rapid rotation, this will tend to happen when they are closely aligned, so that the new particle also resembles a long thin bar.

Whatever the physically reasonable range of choices for radial parameters $\big\{ r_n: n \in \N \big\}$ may be, the natural adaptation of particle dynamics when they are introduced is a change in the pairwise collision rule discussed in Section~\ref{sec.micromodels}. Where before particles $x_i$ and $x_j$ of mass $n$ and $m$ coagulate at rate $\alpha(n,m) \eps^{-2} V\big( (x_i - x_j)/\eps \big)$, we now stipulate that this rate is  $\alpha(n,m) \eps^{-2} (r_n + r_m)^{-2} V\big( \tfrac{x_i - x_j}{(r_n + r_m)\eps} \big)$; the presence of the term $(r_n + r_m)^{-2}$ allows the microscopic coagulation propensity $\alpha(n,m)$ to retain its interpretation of determining the proportion of particle pair overlaps that lead to coagulation (uniformly as the masses of the pair vary).

The kinetic limit derivation of Theorem~\ref{thmo} is undertaken after these changes are made by Rezakhanlou~\cite{RezakhanlouMPRF}. When $d \geq 3$ and the relation $r_n = n^\chi$ is imposed, the derivation is made when $\chi \in \big[0,1/(d-2) \big)$. The macroscopic coagulation rates 
$\beta: \N^2 \to (0,\infty)$ are then found to satisfy
$$
 \beta(n,m) \leq C \big( d(n) + d(m) \big) \big( r_n + r_m \big)^{d-2} \cdot {\rm Cap}\big({\rm supp}(V) \big) \, ,
$$
where the latter term denotes the Newtonian capacity of the support of $V$. 

When dimension $d$ equals three, and we suppose, very reasonably, that $\sup_{n \in \N} d(n) < \infty$, we see that $\beta(n,m) = O(n+m)$ whenever $r_n = O(n)$. In such a regime for $\beta$, it is reasonable to believe that the Smoluchowski PDE is mass conserving for all time; indeed, Proposition~\ref{propuniqueness} comes close to showing this if $d(\cdot)$ decreases gradually enough.

We may tentatively conclude then that the perturbation of our model which includes radial dependence of particles without making more profound changes to inter-particle interaction is not a suitable physical context to study the phenomenon of gelation.

\subsubsection{A weaker notion of gelation: an analogue of weak turbulence}

A weaker notion of solution blowup than finite gelation time is the condition that 
\begin{equation}\label{e.weakturb}
\int_0^{\infty} \int_{\mathbb{R}^d}  m^r f_m(x,t)  \, \dd x \dd m \to \infty \, \, \textrm{as $t \to \infty$,}
\end{equation}
for some $r > 1$. In \cite[Appendix]{YRH}, an analogy is drawn between the non-linear Schr\"odinger equation and the Smoluchowski PDE under which, in the case of the cubic defocussing NLS, the criterion above corresponds to weak turbulence. The condition~(\ref{e.weakturb}) corresponds to ongoing coagulation under which a positive fraction of mass reaches arbitrarily high mass nodes at sufficiently late time. It is argued non-rigorously  in~\cite{YRH} on the basis of scaling considerations for the PDE that, modelling $\beta(n,m) = n^\eta + m^\eta$ and $d(n) = n^{-\phi}$, the behaviour~(\ref{e.weakturb}) is not expected to occur provided that $\eta + \phi < 1$.  


\subsection{The kinetic limit derivation when $d=2$ and with other variants}\label{s.planar}

In \cite{HR2d}, the kinetic limit derivation of the PDE from the models $\PP_N$ was undertaken in dimension $d=2$. We mention here the key changes needed in the models $\PP_N$, and the changes in the recipe for determining $\beta:\N^2 \to (0,\infty)$ from the microscopic parameters. We will also return to the discussion of case $d=2$ in Section~\ref{s.final}, in order to discuss how the proofs in this case differ from when $d \geq 3$.

\medskip

\noindent{\bf Interaction range.} The relation (\ref{e.intrange}) becomes  $N \big| \log \eps \big|^{-1} = Z$, for a given constant $Z \in (0,\infty)$. The interaction range is now exponentially small in $N$, far smaller than it was in the case $d \geq 3$. It is the same regime of constant mean free path that dictates the scale, but now particles are readily available to each other due to Brownian recurrence; small interaction range acts as a countervailing effect.

\medskip

\noindent{\bf Pairwise collision rule.} The infinitesimal rate of coagulation between two particles of mass $n$ and $m$ located at $x_i$ and $x_j$ is now taken to be $\alpha(n,m) \eps^{-2} \big\vert \log \eps \big\vert^{-1} V \big( \tfrac{x_i - x_j}{\eps} \big)$, for a collection of microscopic interaction strengths $\alpha:\N^2 \to (0,\infty)$. The change from the case $d \geq 3$ is the appearance of the factor  $\big\vert \log \eps \big\vert^{-1}$. Its role is to preserve the interpretation of $\alpha(\cdot,\cdot)$ as specifying the proportion of particle overlaps leading to coagulation: were it absent, Brownian recurrence would offer overlapping particles endless opportunities to coagulate, and the proportion of coagulation would be one, for any positive value of $\alpha$. 

\medskip

\noindent{\bf The recipe for the macroscopic coagulation rates.} With a choice of compactly supported interaction kernel $V:\R^2 \to [0,\infty)$ for which $\int_{\R^2} V = 1$, the formula~(\ref{recp}) becomes 
\begin{equation}\label{recptwod}
 \beta(n,m)= \frac{2\pi \cdot \big( d(n) + d(m) \big) \cdot \alpha(n,m)}{2\pi \cdot \big( d(n)
+ d(m) \big) \, + \, \alpha(n,m)} \, .
\end{equation}
Thus, the nature of $V$ is manifest in the macroscopic evolution only through the value of its $L^1$ norm. The reason for this is that, having accepted the presence of a new factor of $\big\vert \log \eps \big\vert^{-1}$ into the formula for $\alpha(\cdot,\cdot)$, any overlapping pair of particles is now not likely to coagulate in any particular excursion into each other's interaction range of duration of order $\eps^2$. Rather, many such opportunities to visit occur for the particles before they move to a large distance from one another, and one among these many visits may cause collision. During all the visits, the details of the form of $V$ no longer really matter, except in a weak law of large numbers' sense, where the average rate of interaction is determined by the $L^1$ norm of $V$.

\subsection{Diffusion and coagulation: two effects from one random dynamics}\label{s.onetwo}

In our microscopic models, the Brownian motion that is the free motion of individual particles is simply a definition, as is the binary coalescence mechanism. Might it be possible to find a microscopic model in which these two phenomena emerge from one microscopic description? Here are two possible answers.

\subsubsection{Physical Brownian motion} 

Physically, Brownian motion arises by the thermal agitation of a particle caused by many random collisions with its neighbours, in a similar vein to the way that the heat equation is expected to arise as a macroscopic evolution equation which we discussed in this survey's opening paragraphs. A physically natural but mathematically presumably intractable microscopic model might suspend comparatively large spheres in an ambient environment of much smaller particles, with a dynamic of elastic collision, and a mechanism of sticking of the large particles on mutual contact. In this sense, the Smoluchowski PDE is sometimes called a model of a colloid. Regarding the important question of the physical derivation of Brownian motion, we mention the recent advance~\cite{BGSR}, concerning the long time behaviour of a tracer billiard in a system of elastic billiards, in a dilute, constant mean free path limit.

\subsubsection{Random walker clusters} A less classical but probabilistically interesting model is the following. Consider a Markov chain whose state space consists of a finite collection of {\em occupied} sites in $\Z^d$. Think that each site is occupied by one walker. Each walker decides to make a transition at the ring times of independent Poisson rate one clocks. Any given walker's transition takes place instantaneously. If a walker is selected to make a transition, he may not move -- his transition is in place -- if his removal from the lattice disconnects a connected component of occupied sites in the nearest neighbour structure. Otherwise, the walker considers making a uniformly random move to one of the nearest-neighbour or diagonally adjacent sites. He does so if the move is to an unoccupied site and the move does not disconnect any two occupied sites that were connected before the move. Otherwise, he stays in place. 

Suppose for a moment that initially the occupied sites are nearest-neighbour connected. The rules are rigged so that this remains so at later times. The ergodicity of the system indicates that the centre of mass of the connected component diffuses in the long term. It is natural to pose the question as to how the diffusivity depends on the mass. Anyway, we obtain a collection of mass-dependent diffusion rates $d:\N \to (0,\infty)$, where now mass means the number of occupied sites in the cluster.

Suppose instead that we begin with a collection of comparatively well separated pairs of nearest neighbours. Each pair begins a random journey which in the large is Brownian. When two clusters meet, they combine, and never break apart subsequently. 

All in all, then, it would seem that with a suitable initial condition and a parabolic scaling of space-time, the model may converge to a solution of the Smoluchowski PDE for some choice of its parameters. Note that the appropriate form of the equation may include the mass-dependent interacting range which we discussed in Subsection~\ref{s.gelation}. One may speculate that $r_n$ should be chosen to scale as $n \to \infty$ according to the scaling satisfied by the typical diameter of an isolated cluster with $n$ occupied sites at equilibrium. Presumably, the cluster has a fractal structure that contributes an exponent of the form $r_n = n^{\chi(1 + o(1))}$.

This microscopic model could be altered so that variants of the Smoluchowski PDE emerge where mass conservation is replaced by conservation of several quantities. Suppose instead that sites are instead occupied either by red or blue particles (but not both), and that the rules are as before, except that blue particles are selected at a rate which is double that for red ones. The diffusion rate of a cluster is now specified by the pair of natural numbers given by the number of constituent red, and blue, particles. This pair replaces the mass as the natural conserved quantity for cluster collision. Convergence to a variant of the Smoluchowski PDE may be expected, and indeed a framework, which may be expected to include the limiting PDE for this example, that specifies variants of the Smoluchowski PDE in which the notion of mass conservation is generalized to conservation of possibly more complex particle characteristics is introduced and analysed in~\cite{Norris14}.  

\section{Homogeneously distributed particles in the torus}\label{s.torus}
We now study the convergence of the microscopic models $\PP_N$ to the limiting system~(\ref{syspde}) in a very particular special case. The choice is made so that, while most of the technical subtleties of definition and proof in the convergence are eliminated, the recipe~(\ref{recp}) for the macroscopic coagulation propensities will be maintained. The main aim of our study of the special case is to explain why the relation~(\ref{recp}) holds; intimately tied to this is a certain microscopic repulsion experienced by the particles at positive macroscopic times, which we also take the opportunity to discuss. Despite these various simplifications, our discussion here is heuristic, with the derivation of several intuitively plausible steps only sketched or omitted entirely; the model in this section is a special case of the annihilating system studied by Sznitman~\cite{sznitman}, and the reader may wish to refer there for a rigorous treatment.

In the special case, under the microscopic models $\PP_N$, there will initially only be particles 
of unit mass, and each will diffuse at rate two. Moreover, as time evolves and particles collide in pairs, the concerned particles will disappear, without the appearance of any new particle. Thus we discuss an interaction of annihilation rather than coagulation.  The particles will initially be placed not in $\R^d$ but rather in the unit $d$-dimensional torus, each placed independently and uniformly with respect to Lebesgue measure. This choice forces the whole dynamics of $\PP_N$ to be invariant under any given spatial translation.

Let $\tord$ denote the $d$-dimensional unit torus, namely the quotient of 
$\R^d$ by $\Z^d$, or the unit cube in 
$\R^d$ with periodic boundary conditions. In this section, $\PP_N$ refers to the annihilating particle dynamics described in the preceding paragraph; we write $\alpha = \alpha(1,1)$ for the microscopic interaction strength of the single particle mass pair in question. In the formal language specifying the Markov generator that we saw for our main object of study in Section~\ref{sec.micromodels}, we are instead setting
the free-motion operator equal to 
$\freem F (q) = \sum_{i \in I_q} \Delta_{x_i} F (q)$,  and  the collision operator is equal to
$\collop F (q) = - \tfrac{\alpha}{2} \sum_{i,j \in I_q}  \eps^{-2} V\big( \tfrac{x_i - x_j}{\eps} \big) F(q)$; note that the absence of a collision gain term is due to our working with annihilation rather than coagulation.

We also make a further minor simplification, choosing the constant $Z$ in (\ref{e.intrange}) to equal one.

The task of making the kinetic limit derivation in this case is to explain how it is that statistics summarising the densities of particles in the microscopic model converge to the appropriate macroscopic evolution, which in this case is given by a function $f:\tord \times [0,\infty) \to [0,\infty)$ satisfying the PDE
\begin{equation}\label{feqn}
  \frac{\partial}{\partial t} f(x,t) =  \Delta f(x,t) - \beta f(x,t)^2
\end{equation}
with initial condition $f(x,0) = 1$ for all $x \in \tord$.  One simplification in our analysis is readily apparent: 
the initial condition has no dependence on the spatial parameter, and this property will be maintained in time. So we may define $h:[0,\infty) \to [0,\infty)$ by setting $h(t) = f(x,t)$ for any choice of $x \in \tord$ and thereby recast (\ref{feqn}) as an ordinary differential equation 
\begin{equation}\label{heqn}
  \frac{\dd}{\dd t} h(t) =  - \beta h(t)^2,
\end{equation}
with $h(0) = 1$. 

In what sense is the evolution of the microscopic models $\PP_N$ approximately summarised by the ODE~(\ref{heqn})?

The spatial homogeneity present in the special case in question offers a simple form for the answer to this question. We introduce a 
{\it microscopic candidate density} $h_N:[0,\infty) \to [0,\infty)$, a quantity which summarises the density of particles in the microscopic model $\PP_N$ and which we hope to show approximates $h:[0,\infty) \to [0,\infty)$ when $N$ is high. We set $h_N(t) = N^{-1} \mathbb{E} s_N(t)$, where $s_N(t)$ is the number of {\it surviving} particles at time $t$, namely, the mean number of particles which have not been annihilated before time $t$.  (In the general, spatially inhomogeneous, setting, we will also define a microscopic candidate density, but its definition will be a little more involved, it will be random rather than deterministic, and it will depend not only on the time variable but also on the macroscopic location.)

In the special case, our kinetic limit derivation amounts to explaining how it is that $h_N:[0,\infty) \to [0,\infty)$ converges as $N \to \infty$ to the unique solution $h:[0,\infty) \to [0,\infty)$ of~(\ref{heqn}). The principal aim of this section is to justify heuristically the relation (\ref{recp}) between the quantity $\beta$ appearing in~(\ref{heqn}) and the microscopic parameters. In the present case, this relation takes the following form.

\begin{proposition}\label{prop.heur}
Assume that $V:\R^d \to [0,\infty)$ is continuous and compactly supported.
Then the functions $h_N:[0,\infty) \to [0,\infty)$ converge pointwise as $N \to \infty$ to the unique solution of~(\ref{heqn}), where $\beta$ is specified by
\begin{equation}\label{e.betaheur}
 \beta =  \alpha  \int_{\R^d} \big( 1 - u(x) \big) V(x) \dd x \, ,
\end{equation}
with  $u:\R^d \to [0,1]$ being the unique solution (provided by Lemma~\ref{l.s}) subject to $u(x) \to 0$ as $x \to \infty$ of the modified Poisson equation
$$
 - 2 \Delta u (x) = \alpha V(x) (1-u)(x) \, .
$$  
\end{proposition}

As we prepare to justify the proposition, it is convenient to recast the definition of $h_N$ in terms of the evolution of a particle, which we will call the tracer particle, picked uniformly at random at time zero. Since the distribution of particles at the initial time is invariant under particle reindexing, the next definition is suitable. 
\begin{definition}
The tracer particle is the particle whose time zero index is $1$. 
\end{definition}
The relationship between $h_N$ and tracer particle survival probability is straightforward.
\begin{lemma}\label{lemhn}
Let $N \in \N$ and $t \geq 0$. Then the microscopic candidate density $h_N(t)$ is equal to the $\PP_N$-probability that the tracer particle has survived until time $t$. 
\end{lemma}
\noindent{\bf Proof.}  Note that 
$$
 h_N(t) =  N^{-1} \sum_{i=1}^N \PP_N \Big( \textrm{the $i\textsuperscript{th}$ indexed particle survives until time $t$} \Big) \, .
$$
However, the summand is independent of $i \in [1,N]$ due to the symmetry in both the initial placement of particles and in their dynamics. \qed

In seeking to argue that $h_N$ converges to the unique solution $h$ of~(\ref{heqn}) in some appropriate sense, 
it is natural to try to find an expression for $\tfrac{\dd}{\dd t} h_N(t)$. We may hope to show that in fact this derivative equals $- \beta h_N(t)^2$ up to some error term which in some way tends to zero in the high $N$ limit. Since the initial conditions $h_N(0)$ (for $N \in \N$), and $h(0)$, all coincide (with $1$), we might then argue that $h_N \to h$ in some sense as $N \to \infty$.

With this aim in mind, we find an expression for the derivative of $h_N$ in terms of the behaviour of the tracer particle in the microscopic model $\PP_N$:
\begin{lemma}\label{lemtpexpress}
Let $N \in \N$ and $t \geq 0$. Let $\surv_t$ denote the event that the tracer particle under $\PP_N$ survives until time $t$. Then
\begin{equation}\label{hnexp}
 \frac{\dd}{\dd t} h_N(t) =  
  - \lim_{\delta \searrow 0} \delta^{-1} \PP_N \big( \surv_t \big) \PP_N \big( \surv^c_{t+ \delta} \big\vert \surv_t \big) \, ,  
\end{equation} 
should the limit on the right-hand side exist.
\end{lemma}
\noindent{\bf Proof.} By Lemma~\ref{lemhn}, $h_N(t) = \PP_N(\surv_t)$, and thus $h_N(t + \delta) - h_N(t) = \PP_N\big( \surv_t \cap \surv^c_{t + \delta} \big)$. \qed

The expression (\ref{hnexp}) gives us a probabilistic means of thinking about the derivative of $h_N$. We should consider
\begin{itemize}
 \item the $\PP_N$-probability that the tracer particle survives until time $t$; and
 \item given that it does so, the conditional probability that it is instantaneously annihilated in a collision.
\end{itemize}

There is a particular value of $t$ for which the probability of this event is easier to evaluate: $t = 0$. In this case, the survival probability $\PP_N(\surv_0)$ is trivially equal to one. What then is the $\PP_N$-probability that the tracer particle is annihilated before a very short time $\delta$ has passed? As our discussion will now tend to be heuristic rather than rigorous, we write subsequent statements as claims rather than lemmas.

\begin{claim}\label{c.collprob}
Let $t = t_N$ satisfy $t = o(\eps^2)$ as $N \to \infty$; equivalently, by~(\ref{e.intrange}), $t= o(N^{2(2 - d)})$.
Then $\PP_N \big( \surv_t^c \big) = t \alpha \int_{\R^d} V(x) {\rm d} x \big( 1 + o(1) \big)$ as $N \to \infty$. 
\end{claim}
\noindent{\bf Sketch of proof.}
We begin by estimating the probability that the tracer particle collides with a given other particle in a very short interval $[0,t]$. In the case of the particle with index two, this probability is by definition given by 
\begin{equation}\label{e.collprob}
1 - \exp \Big\{ - \alpha \int_0^t V_\eps \big( X_2(s) - X_1(s) \big) \mathbf{1}_{ \{1,2 \} \subseteq I_{q(s)}} \, {\rm d} s \Big\} \, , 
\end{equation}
since the exponential term here is the probability that the Poisson process for collision of particles indexed by $1$ and $2$ has yet to ring by time $t$. Note the presence of the indicator function for the event  $\{1,2 \} \subseteq I_{q(s)}$ that two particles have yet to be annihilated by time $s$. (The notation $I_{q(s)}$, in which we now explicitly refer to the time parameter $s$, was introduced back in Section~\ref{sec.micromodels}: it is the set of particles that are surviving at time $s$.) In fact, this indicator function may be dropped from the expression~(\ref{e.collprob}) at the expense of a lower order term as $t \searrow 0$ because, as we explained in Section~\ref{sec.meanfreepath}, there is asymptotically zero probability that either of the two particles are annihilated during $[0,s]$ as $s \searrow 0$. Recalling that $V_\eps(\cdot) = \eps^{-2} V(\cdot/\eps)$ and that $V:\R^d \to [0,\infty)$ is supposed to be continuous, and noting that the difference $X_2 - X_1$ is a rate four Brownian motion, we see that 
$$
\int_0^t V_\eps \big(X_2(s) - X_1(s) \big) \, {\rm d} s = 
\int_0^t V_\eps \big(X_2(0) - X_1(0) \big) \big( 1 + o(1) \big) \, {\rm d} s
$$
provided that $t = t_N$ is chosen so that $t \eps^{-2} \searrow 0$ as $N \to \infty$. That is to say, as the total particle number $N$ tends to infinity, the collision probability~(\ref{e.collprob}) on $[0,t]$ is accurately approximated by 
$1 - e^{-t \alpha V_\eps \big(X_2(0) - X_1(0)\big)}$ provided that $t$ tends to zero more quickly than $\eps^2$, because, in this limiting regime, the locations $X_1$ and $X_2$ are asymptotically static on scale~$\eps$.  Since $t V_\eps(x)$ converges to zero uniformly in $x \in \R^d$ in this regime, our asympotic expression for (\ref{e.collprob}) is $t \alpha V_\eps \big(X_2(0) - X_1(0)\big)$.
 Since $X_2(0) - X_1(0)$ is simply uniformly distributed in $\tord$, this quantity, after averaging over $X_1(0)$ and $X_2(0)$, equals $t\alpha  \int_{\R^d} V_\eps (s) {\rm d} s$, which is $t \alpha \eps^{d -2} \int_{\R^d} V(x) {\rm d} x$. 

By symmetry of the particle indices, this estimate applies equally to the probability of collision between particles with any two given indices in $[1,N]$. Since the probability that the tracer particle experiences two collisions during $[0,t)$ (with $t = o(\eps^2)$) behaves as ${N \choose 2} \big( t \alpha \eps^{d -2} \int_{\R^d} V(x) {\rm d} x \big)^2 = o(\eps^4)$, which is much smaller than the $\eps^2$-order probability of a single such collision, the probability that the tracer particle experiences some collision during $[0,t)$ is well approximated by the mean number of collisions that it experiences, which is 
$$
 (N - 1) t \alpha \eps^{d -2} \int_{\R^d} V(x) {\rm d} x \, .
$$ 
Recalling that $N = \eps^{2 - d}$, we obtain the claim. \qed

We are ready to return to~(\ref{hnexp}) and record a limiting expression in high $N$ for the time zero derivative of the microscopic candidate density:

\begin{claim}\label{lemhnexpzero}
Let $N \in \N$. Then
\begin{equation}\label{hnexpzero}
  \lim_N \frac{\dd h_N}{\dd t} (0) =  - \alpha  \int_{\R^d} V (x) \dd x \, .
\end{equation} 
\end{claim}
\noindent{\bf Proof.}  Since the expression $\lim_N \frac{\dd}{\dd t} h_N(0)$
equals $\lim_{N \to \infty} \lim_{\delta \searrow 0} \delta^{-1} \PP_N(\surv_\delta^c)$, the claim follows from Claim~\ref{c.collprob}. \qed

Recall that our aim is to show that, in some appropriate sense, $h_N$ converges to the solution $h$  
of~(\ref{heqn}) as $N \to \infty$. Claim~\ref{lemhnexpzero} gives us a guess for the derivative at time zero of $h$: if it is the limit of the derivatives of its anticipated approximants, then $h'(0) =  - \alpha  \int_{\R^d} V (x) \dd x$. However, by its definition~(\ref{heqn}), $h'(0)$ is also $- \beta h(0)^2$, which is simply $- \beta$. In other words, the preceding argument has given us a guess for the recipe by which the macroscopic coagulation propensity $\beta$ is to be computed from the parameters in the underlying microscopic models. Namely, the argument points to the conclusion that
\begin{equation}\label{eqnbetanaive}
 \beta = \alpha \int_{\R^d} V(x) \dd x \, . 
\end{equation}
\begin{subsection}{Surviving and ghost particles}
However, this guess is wrong. The formula (\ref{eqnbetanaive}) is not the correct relation between the microscopic and macroscopic coagulation propensities. To see why this is so,
it is useful to introduce a coupling of our annihilating Brownian dynamics $\PP_N$ with a system of independent non-interacting Brownian particles $\PP_N'$.
\begin{definition}
Let $\coup$ denote a coupling of the law $\PP_N$ with a further law $\PP_N'$. Under $\PP_N'$, $N$ particles are scattered in $\tord$ at time zero with the same law as in $\PP_N$, and under $\coup$ the two initial conditions are always equal. In both marginals under $\coup$, each particle pursues a given rate two Brownian motion, independently of the others. In the $\PP_N$ marginal, particles disappear on collision according to the rule for that dynamics; in the $\PP_N'$ marginal, the collision has no effect on either particle, and each continues its Brownian trajectory undisturbed.   
\end{definition}
We say that under $\coup$ each particle is initially {\it surviving}. When a collision event occurs between two surviving particles, each becomes a ghost. In this way, the collection of all particles has the law~$\PP_N'$ while the collection of surviving particles has the law $\PP_N$. 

For now, we need only one consequence of the coupling, namely that, for any $t \geq 0$, the time-$t$ marginal of $\PP_N$ is stochastically dominated by its time zero marginal. To see this,  note that this time-$t$ marginal is dominated by the time-$t$ marginal of all particles which has the law of the time-$0$ marginal of $\PP_N$.

\end{subsection}
\begin{subsection}{The naive guess is wrong}\label{s.torussec}
We now give an intuitive explnation of why (\ref{eqnbetanaive}) is the wrong relation between $\beta$, $\alpha$ and $V$. Suppose for convenience that $\int_{\R^d} V(x) \dd x = 1$. Now choose $\alpha$ to be fixed but very high. We are left with the formula $\beta = \alpha$, so that the solution $h:[0,\infty) \to [0,\infty)$ of~(\ref{heqn}) with $h(0)=1$ equals $h(t) = \frac{1}{1 + \alpha t}$.

Our high choice of $\alpha$ means that $h$ drops towards zero quickly after time zero: specifically, $h(\alpha^{-1/2}) \leq \alpha^{-1/2}$. Here, however, we encounter a difficulty when we think of $h$ as a limit of its approximants $h_N$. Expecting that $h_N(t) \to h(t)$ pointwise, recall from Lemma~\ref{lemhn} that $h_N(t) = \PP_N(\surv_t)$ is the survival probability of the tracer particle until time $t$ under $\PP_N$. 
If the annihilation event $\surv^c_t$ is to occur, then it is necessary that at some time before~$t$, some other particle enters the interaction range of the tracer particle. By the coupling~$\coup$ with an independent system of Brownian particles, to find an upper bound on $\PP_N (\surv^c_t )$, it is enough to bound the probability that, among $N$ uniformly scattered particles in $\tord$ each performing Brownian motion (of rate two), a given particle comes at some time in $[0,t]$ to distance $\eps$ of some other.
Since the probability of such an approach at any given time is $N \eps^d = \eps^2$, and such an approach occurs for a mean duration of order $\eps^2$, this probability is bounded above by a constant multiple of $t$ for small $t$. (An explicit bound on this mean time is that the time spent overlapping during $[0,\infty)$ by two rate-two radius-$\epsilon$ Brownian spheres in $\R^d$, with $d \geq 3$, which are tangent at time zero is in expectation at most $\big( 1 + \tfrac{2^{1-d}}{(d-2)\Gamma(d/2 + 1)} \big) \eps^2/2$.)  Thus, $\PP_N \big( \surv_t \big) \geq 1 - C t$, where $C > 0$ may be chosen  uniformly in  both $N \in \N$ and $\alpha \in (0,\infty)$. 
For all $\alpha > 0$, the pointwise convergence of $h_N$ to $h$ forces $h(t) \geq 1 - Ct$. However, this is inconsistent with $h(\alpha^{-1/2}) \leq \alpha^{-1/2}$ for $\alpha > (C + 1)^2$. We conclude then that the guessed formula~(\ref{eqnbetanaive}) is in fact wrong.

What is wrong with the derivation of~(\ref{eqnbetanaive}) is that, in fact, the initial Poissonian distribution of particles is in a certain sense unstable, making an inference based on an analysis at time zero misleading. 
Although at later times the particle distribution is Poissonian in the large, there is a microscopic repulsion effect which modifies this: in the order $\eps$ vicinity of the tracer particle at some positive time, there is a diminished probability for presence of another particle. Indeed, this other particle may have already collided with the tracer particle, in which case, the tracer particle and the other particle would not in fact be located close to each other because each would have vanished.  We now turn to quantifying the effect of this mechanism of curtailment of interaction due to collision and so derive the correction needed to~(\ref{eqnbetanaive}). To make sense of the notion of a particle in the vicinity of the tracer particle which may already have vanished at a certain time, we will make use of our coupling of surviving and ghost particles.

\end{subsection}
\begin{subsection}{Quantifying the correction}
The discussion in the preceding section reveals that our computation of $\tfrac{\dd h_N(t)}{\dd t}$ for $t = 0$
may have exceptional features which change as $t$ increases. In fact, as we will see, this change will be apparent already when $t$ reaches the order of $\eps^2$. To understand the change quantitatively, we want to return to
Lemma~\ref{lemtpexpress} and use the right-hand side of (\ref{hnexp})
to compute $\tfrac{\dd h_N(t)}{\dd t}$ for $t > 0$. 
This involves computing the terms $\PP_N(\surv_t)$ and $\PP_N\big( \surv_{t + \delta}^c \big\vert \surv_t \big)$ appearing on the right-hand side of (\ref{hnexp}), when $t > 0$ is fixed, and $\delta > 0$ is infinitesimally small. Reexpressing the events in terms of the coupling $\coup$, $\surv_t$ is the event that the tracer particle survives until time $t$, and  $\surv_t \cap \surv_{t + \delta}^c$ is the event that the tracer particle changes its status from surviving to ghost during the short time interval $[t,t + \delta]$.

\begin{claim}\label{cl.big}
Suppose that $\delta \eps^{-2} \searrow 0$ as $N \to \infty$. Then
$$
 \PP_N\big( \surv_{t + \delta}^c \cap \surv_t \big) = \delta \alpha h(t)^2 \int_{\R^d} \big( 1 - u(x) \big) V(x) \dd x \, \Big( 1 + o(1) \Big) \, , 
 $$
as $N \to \infty$, where $u:\R^d \to [0,1]$ is specified in Proposition~\ref{prop.heur}.
\end{claim}

Seeking to justify this claim, we let $C_{i,j}(t,\delta)$ denote the event that
\begin{itemize}
 \item $\vert X_i(t) - X_j(t) \vert \leq \eps$, 
 \item and these two particles collide during $[t,t+\delta]$.
\end{itemize}
As we now record, the event $\surv_t \cap \surv_{t + \delta}^c$ is characterized up to a probability of smaller order by the intersection of the following events:
\begin{itemize}
\item the survival of the tracer particle until time $t$; 
\item the presence at time $t$ of some other surviving particle in $X_1(t) + \ball_\eps$, (a set that contains the interaction range of the tracer particle);
\item and the collision of that other particle with the tracer particle  during $[t,t + \delta]$. 
\end{itemize}
(Here, $\ball_\eps$ denotes the Euclidean ball of radius $\eps$ about the origin, so that, since the support of $V$ is contained in the Euclidean unit ball, the set  $X_1(t) + \ball_\eps$ indeed contains the interaction range of the tracer particle.) 
\begin{claim}\label{cl.union}
For each $t > 0$, there exists $C_t > 0$ such that
$$
 \coup \bigg(  \big( \surv_t \cap \surv_{t + \delta}^c  \big)  \, \Delta \, \Big( \bigcup_{j=2}^N \big\{ \{ 1,j \} \in I_{q(t)} \big\} \cap C_{1,j}(t,\delta) \Big)  \bigg)
  \leq  C_t \coup \Big(   \surv_t \cap \surv_{t + \delta}^c \Big)  
   \big(    \delta +  \alpha \eps \delta^{1/2} \big) \, .  
$$
\end{claim}
\noindent{\bf Sketch of proof.} 
The claim will emerge from two assertions. First,
\begin{equation}\label{e.intone}
 \coup \bigg(  \big( \surv_t \cap \surv_{t + \delta}^c  \big) \, \Delta \, \Big( \bigcup_{j=2}^N \big\{ \{ 1,j \} \in I_{q(t)} \big\} \cap C_{1,j}(t,\delta) \Big) \bigg)
  \leq  
    C  \delta^2 + C \alpha \eps \delta^{3/2}  \, ,
\end{equation}
and, second, for some $t$-dependent constant $c > 0$,
\begin{equation}\label{e.inttwo}
\coup \Big(  \surv_t \cap \surv_{t + \delta}^c \Big) \geq  c \delta \, .
\end{equation}
The second bound holds because, as we described in Section~\ref{sec.meanfreepath}, the tracer particle will survive to any time $t$ with some positive, $t$-dependent probability, and it is then liable to collide with some other particle at a rate of order one. 

Regarding the symmetric difference in~(\ref{e.intone}), note that, if one event occurs without the other, the cause must be one of the following:
\begin{itemize}
 \item although  $\{ 1,j \} \in I_{q(t)}$ and  $C_{1,j}$ occur for some $j \in [2,N]$, there is a third particle which collides with particle $j$ after time $t$ but before the collision of $j$ with $1$ that happens before time $t + \delta$;  
 \item or, $\surv_t \cap \surv_{t + \delta}^c$ occurs due to some particle $j \in [2,N]$, which is not in the interaction range $X_1 + \ball_\eps(t)$ at time $t$, entering this range and colliding with $1$ during $[t,t+\delta]$. 
\end{itemize}
Regarding the first possibility, for given $j \in [2,N]$, the probability of $C_{1,j}$ is at most $C \eps^d \delta \eps^{-2}$, while, given this event, the conditional probability of some third particle behaving as described is at most $C \delta$. Summing over $j$, the probability is at most $C N \delta^2 \eps^{d-2} = C \delta^2$. 

For given $j \in [2,N]$, the probability of the second occurrence is at most a constant multiple of~$\eps \delta^{3/2}$. We only sketch how this bound is obtained. Should particle $j$ at time $t$ lie within a distance of order $\delta^{1/2}$ of the boundary of the interaction range $X_1 + \ball_\eps$ of the tracer particle, there is positive probability that particle $j$ enters this range during the ensuing $\delta$ units of time, and should this happen, there is conditional probability at most 
$1 - \exp \{ - \alpha \vert V \vert_{\infty} \delta \} \leq \alpha  \vert V \vert_{\infty} \delta \leq C \alpha \delta$ of collision between the two particles. The probability of this turn of events is thus $\big\vert \ball_{\eps + \delta^{1/2}} \setminus \ball_{\eps}\big\vert C \alpha \delta \leq C \eps^{d-1} \delta^{3/2}$. On the other hand, it is easily checked that there is negligible probability of such a collision should particle $j$ at time $t$ lie at much greater distance than $\delta^{1/2}$ from the boundary of $X_1 + \ball_\eps$. Summing over the $N -  1 \sim \eps^{2 - d}$ choices of $j \in [2,N]$, we see that the probability of the second listed event is at most $C \eps \delta^{3/2}$. \qed

The next two claims estimate the probability of $\big\{ \{ 1,j \} \in I_{q(t)} \big\} \cap C_{1,j}(t,\delta)$ for $j \in [2,N]$ and will lead to Claim~\ref{cl.big}.

\begin{claim}\label{cl.one}
Suppose that $\delta \eps^{-2} \searrow 0$ as $N \to \infty$. Then, for each $x \in \ball_\eps$, 
  $$
   \coup \Big( C_{1,2}(t,\delta)  \, \Big\vert \,  X_1(t) - X_2(t)  = x  \Big) = \eps^{-2} \delta \alpha  V(x/\eps)   \, \big( 1 + o(1) \big) \, .
  $$
\end{claim}
\begin{claim}\label{cl.two}
For each $x \in \ball_\eps$, we have that
  $$
   \coup \Big( \{ 1,2 \} \subseteq I_{q(t)} \, \Big\vert \,  X_1(t) - X_2(t)  = x  \Big) =  \big( 1 - u(x/\eps) \big) h(t)^2  \, .
  $$ 
\end{claim}

\noindent{\bf Proof of Claim~\ref{cl.big}.} 
Note that, conditionally on $X_2(t) - X_1(t) \in \ball_\eps$, $X_2(t) - X_1(t)$ is uniform on $\ball_\eps$. In light of this, and 
Claims~\ref{cl.one} and~\ref{cl.two}, we see that, in the limit in question, 
$$
 \coup \big( C_{1,2}(t,\delta) , \{ 1,2 \} \subseteq I_{q(t)}   \big) =  \eps^{-2} \delta \alpha  h(t)^2 \int_{\ball_\eps} V(x/\eps) \big( 1 - u(x/\eps) \big)  {\rm d} x \, \Big( 1 + o(1) \Big) \, ,   
$$
whose integral term may be also written as $\eps^d \int_{\R^d} V(x) \big( 1 - u(x) \big) {\rm d} x$. By particle symmetry and $N - 1 \sim \eps^{2 - d}$, we obtain
$$
\sum_{j=2}^N  \coup \Big( C_{1,j}(t,\delta) , \{ 1,j \} \subseteq I_{q(t)}   \Big) =  \delta \alpha  h(t)^2 \int_{\R^d} V(x) \big( 1 - u(x) \big)  {\rm d} x  \, \Big( 1 + o(1) \Big) \, .
$$
It is easy to convince oneself that typically the occurrence of $\cup_{j=2}^N C_{1,j}(t,\delta) \cap \big\{  \{ 1,j \} \subseteq I_{q(t)} \big\}$ entails the occurrence of exactly one of the constituent events. For this reason, the preceding equality holds equally for the probability of $\cup_{j=2}^N C_{1,j}(t,\delta) \cap \big\{  \{ 1,j \} \subseteq I_{q(t)} \big\}$. Thus, Claim~\ref{cl.big} follows from Claim~\ref{cl.union}. \qed

\medskip

\noindent{\bf Proof of Claim~\ref{cl.one}.} Conditionally on $X_1(t) - X_2(t)$ being a given $x \in \ball_\eps$, 
the probability of collision between $1$ and $2$ during $[t,t+\delta]$ is $\delta \alpha V_\eps(x) \big( 1 + o(1) \big)$, since $V:\R^d \to [0,\infty)$ is assumed continuous, and $X_1 - X_2$ is asymptotically static as $\e \searrow 0$ on scale~$\eps$ during this duration of length~$\delta = o(\e^2)$.  \qed 

\medskip 

\noindent{\bf Proof of Claim~\ref{cl.two}.}
To reiterate the problem, given that at time $t$ particles $1$ and $2$ have displacement $x \in \ball_\eps$, what is the probability that both survive to this time? There are two reasons why one or other may be a ghost particle at time $t$:
\begin{itemize}
 \item it may be that, during $[0,t]$, at a moment when each of $X_1$ and $X_2$ are surviving, a collision between this pair occurs;
 \item it may be that one or other of $X_1$ and $X_2$, at a moment when this particle is surviving, collides with some other surviving particle.
\end{itemize}
Calling these two events $E_1^t$ and $E_2^t$, we want to gauge the probability of $\big( E_1^t \big)^c \cap \big( E_2^t \big)^c$.

In considering these possibilities, it is convenient to reverse time, running time backwards from $t$ to $0$. We will now use these time coordinates, where the forward time evolution from $0$ to $t$ corresponds to the usual evolution backwards from $t$ to $0$. Note that, because we do not condition on $X_1$ or $X_2$ surviving until time $t$, the conditional distribution of the trajectories $X_1,X_2:[0,t] \to \tord$ in the new time coordinates is a pair of independent rate two Brownian motions, where $X_1(0)$ is uniformly distributed on $\tord$ and $X_2(0) = X_1(0) + x$.  Phrased in these terms, $E_1^t$ is the event of collision between $X_1$ and $X_2$ during $[0,t]$.

In light of Lemma~\ref{l.gfin}, we learn that
\begin{claim}\label{cl.eone}
If $t = t_N$ satisfies $t/\eps^2 \to \infty$ as $N \to \infty$, then  
 $$
   \coup \Big(  \big( E_1^t \big)^c \Big) = \big( 1 - u(x/\eps) \big) \big( 1 + o(1) \big) \, ,
 $$  
 where  $u:\R^d \to [0,1]$ is specified in Proposition~\ref{prop.heur}.
 \end{claim}
 We also need to estimate the conditional probability that $E_2^t$ occurs given that $E_1^t$ does not.
 \begin{claim}\label{cl.etwo} 
If $t = t_N$ satisfies $t/\eps^2 \to \infty$ as $N \to \infty$, then  
$$
   \coup \Big( \big( E_2^t \big)^c \, \Big\vert \, \big( E_1^t \big)^c \Big) = h(t)^2 \big( 1 + o(1) \big) \, .
$$
 \end{claim}
\noindent{\bf Sketch of proof.} Should $E_1^t$ not occur, the two particle trajectories $X_1,X_2:[0,t] \to \tord$, begun at the points $0$ and $x$ at distance of order $\eps$, will not experience collision, and will separate to a distance much greater than $\eps$ in a time whose order is large compared to $\eps^2$; by our assumption on the time $t$, this separation occurs on a time scale (called $s$) much shorter than $t$, and, after it has done so, it is reasonable to think that the collision behaviour of the two trajectories will become effectively independent of one another. Any given one of these trajectories experiences collision with some other particle with probability $1 - h(t-s)$, by the definition of $h$; assuming this independence, both particles survive collision during $[0,t]$ with probability $h(t-s)^2$. However, since $s \ll t$, $h(t - s) = h(t) \big( 1 + o(1) \big)$ as $N \to \infty$. \qed

Claim~\ref{cl.two} follows from Claims~\ref{cl.eone} and~\ref{cl.etwo}. \qed

We are now able to complete our sketch proof of the formula for $\beta$ for the model in question.

\medskip

\noindent{\bf Sketch proof of Proposition \ref{prop.heur}.}
Applying Lemma~\ref{lemtpexpress} and Claim~\ref{cl.big}, we find that
 $$
  \frac{\dd h_N(t)}{\dd t} = -  \alpha h(t)^2 \int_{\R^d} \big( 1 - u(x) \big) V(x) \dd x \, \Big( 1 + o(1) \Big) \, ,
 $$
provided that $t/\eps^2 \to \infty$ as $N \to \infty$. 
This estimate does not control the behaviour of $\tfrac{\dd h_N(t)}{\dd t}$ on the small time scale where $t$ is of order $\epsilon^2$. Here, however, arguments in the style of those leading to Claim~\ref{lemhnexpzero} justify that this derivative is non-positive and bounded below by
 $-\alpha \int_{\R^d} V(x) \dd x$. Recalling that $h$ is the unique solution of~(\ref{heqn}), we see that, since $h_N(0) = h(0) = 1$ for all $N \in \N$,   
$h_N$ converges  to $h$ pointwise as $N \to \infty$. \qed
\end{subsection}
\subsection{Bose-Einstein condensates, and a parallel macroscopic interaction}
Our computation of the macroscopic interaction rate finds a parallel in the quantum mechanical problem of the dynamics of a collection of $N$ bosons in three dimensions that interact via a short-range pair potential, that was investigated during Erd\"os, Schlein and Yau's derivation~\cite{ErdosSchleinYau} of the macroscopic evolution of the system.  We briefly discuss this now.
The dynamics of the system is governed by the Schr{\"o}dinger equation 
$$
 i \partial_t  \psi_{N,t} = H_N \psi_{N,t} \ ,
$$
for the wave function $\psi_{N,t} \in L^2_s(\R^{3N})$, the subspace of $L^2(\R^{3N})$ consisting of all functions symmetric under permutations of the $N$ particles. Short-range repulsive interaction is modelled by the choice of Hamiltonian $H_N = H_{\beta,N}$, 
$$
 H_{\beta,N} = - \sum_{j=1}^N \Delta_j \, + \, \tfrac{1}{N} \sum_{1 \leq i < j \leq N} N^{3\beta} V \big( N^\beta(x_i - x_j) \big) \, ,
$$
where $V:\R^3 \to [0,\infty)$ is a compactly supported interaction potential, and $\beta > 0$ is a parameter. (Note that this notation is in conflict with our use of $\beta(\cdot,\cdot)$.) The choice of $\beta=1$ provides the closest parallel with the main discussion in this survey and is the principal object of study in~\cite{ErdosSchleinYau}. The macroscopic evolution of the system may be summarized by a decoupling property enjoyed by the $k$-particle reduced density matrices; these matrices in~\cite{ErdosSchleinYau} are shown by an analysis of the BBGKY hierarchy  to factorize asymptotically in the high $N$ limit, with the factor governed by the non-linear Gross-Pitaevskii equation
$$
 i \partial_t \psi_t = \Delta \psi_t + \sigma \vert \psi_t \vert^2 \psi_t \, ,
$$
where the coupling constant $\sigma$ is given by 
\begin{equation*}
 \sigma = 
\begin{cases}
  \, b_0  &  \textrm{ if $0 < \beta < 1$} \, ,\\
  \, 8 \pi a_0 & \textrm{ if $\beta = 1$}\,  .
\end{cases}
\end{equation*}

Here, $b_0 = \int_{\R^3} V(x) \dd x$, while $a_0$ satisfies
$$
 a_0 =  \tfrac{1}{8 \pi} \int_{\R^3} V(x) \big( 1 - \omega_0(x)\big) \dd x \, ,
$$
with $\omega_0$ being the unique solution to 
$$
 \Big[ - \Delta +  \tfrac{1}{2} V(x) \Big]  \big( 1 - \omega_0(x) \big) = 0 
$$
that satisfies $\lim_{x \to \infty}\omega_0(x) = 0$. (Note that $\omega_0$ is nothing other than the solution of~(\ref{pdew}) if we take $\tfrac{\alpha(n,m)}{d(n) + d(m)}$ equal to $1/2$.)

That is, the macroscopic interaction coefficient undergoes a transition as $\beta \in (0,1)$ changes to $\beta = 1$ in precise correspondence to the modification from the naive guess $\beta = \int_{\R^d} V(x) \dd x$ to $\beta$ given by~(\ref{recp}) which we have devoted this section to discussing. 

Indeed, we may specify a collection of random models $\PP_{N,\beta}$, with $\beta > 0$, in such a way that our models $\PP_N$ coincide with $\PP_{N,1}$, while $\PP_{N,\beta}$, $\beta \in (0,1)$, form counterparts to the interacting bosonic systems at such values of $\beta$. Maintaining the relation~(\ref{e.intrange}) between $\eps$ and $N$, we modify the pairwise collision rule from Section~\ref{sec.micromodels} from one under which the particles indexed by $i,j \in [1,N]$ collide at rate $\alpha(m_i,m_j) \eps^{-2} V\big( \tfrac{x_i - x_j}{\eps} \big)$ to one whose rate is $\alpha(m_i,m_j) \eps^{-2 + d(1-\beta)} V\big( \tfrac{x_i - x_j}{\eps^\beta} \big)$. The rule for $\PP_{N,\beta}$ is determined in order that a typical particle maintain a unit-order  interaction with all the others per unit time, so that the new models remain in the regime of constant mean free path. The derivation of Theorem~\ref{thmo} may be reprised for choices of $\beta \in (0,1)$, with the formula for the macroscopic interaction rates 
$\beta(n,m)$ now given by 
$\beta(n,m) = \alpha(n,m) \int_{\R^d} V(x) \dd x$. The new formula holds in essence because microscopic pairwise repulsion is absent asymptotically in high $N$ in these models.

See \cite{TaoBlog} for a blog post by Terry Tao, written after a talk by Natasa Pavlovic, which provides a more informative summary of this quantum problem, including at its end, and in the ensuing comments, mention of the dichotomy between interaction coefficient in the cases $\beta \in (0,1)$ and $\beta =1$.    

\section{The route to Theorem 1.1 
}\label{s.routetothm}
In this section, we explain the overall plan for proving the main theorem, and reduce it to a fundamental proposition, the {\it Stosszahlansatz}, which concerns the total particle coagulation propensity in the microscopic models. 
\begin{subsection}{Approximating the PDE using microscopic candidate densities}
\subsubsection{A microscopic counterpart to the PDE in weak form}
Recall the weak formulation~(\ref{syspdeweak}) of the Smoluchowksi PDE. Our aim is to show that the particle densities in the microscopic model $\PP_N$ converge in a suitable sense to this weak solution. To do so, we find a microscopic counterpart to the equation~(\ref{syspdeweak}), namely an equation expressed in terms of the law $\PP_N$. Note that~(\ref{syspdeweak}) expresses the change in the quantity $\int_{\R^d} J_n(x,t) f_n(t) \dd x$ that occurs between times $0$ and $T$ as an integral over the intervening duration $[0,T]$ of the differential changes caused by variation in the test function $J_n$, and by the diffusion and coagulation of the particles being modelled. The quantity $\int_{\R^d}J_n(x,t) f_n(t) \dd x$  is an expression for the total number of particles of mass $n$ at time $t$, where each particle is weighted by $J_n$. As such, it has a clear microscopic analogue: under the law $\PP_N$, the random variable $\sum_{i \in I_{q(t)}} J_n \big( x_i, m_i \big) \mathbf{1}_{m_i(t) = n}$, which is the sum over 
mass-$n$ particles at time $t$ where 
each particle carries a weight given by $J_n$.
The form of the Markov generator for the dynamics of $\PP_N$ now provides us with an analogue of the weak formulation~(\ref{syspdeweak}) of the Smoluchowski PDE:
\begin{eqnarray}
 & & \eps^{d-2} \sum_{i \in I_{q(T)}} J_n \big( x_i, T \big) \mathbf{1}_{m_i(T) = n} -  \eps^{d-2}  \sum_{i \in I_{q(0)}} J_n \big( x_i, 0 \big) \mathbf{1}_{m_i(0) = n} \label{emicro} \\
  & = &  \eps^{d-2}  \int_0^T \bigg( \sum_{i \in I_{q(t)}} \frac{\partial J_n}{\partial t} \big( x_i, t \big) \mathbf{1}_{m_i(t) = n}
    \, + \,   \sum_{i \in I_{q(t)}} d(n) \Delta J_n \big( x_i, t \big) \mathbf{1}_{m_i(t) = n} \nonumber \\
  & & \quad \quad \quad \, + \,   \sum_{i,j \in I_{q(t)}}  \alpha(m_i,m_j) V_\eps \big( x_i - x_j \big)  \collision_{i,j,t,n} J_n  \bigg) \, \dd t \, \, + \, M_T \, . \nonumber
\end{eqnarray}
In the integrands on the right-hand side, we see the infinitesimal mean changes caused by the time-variation of $J$, by the diffusion of the individual particles, and by their collision in pairs; the final term is a martingale (which we will argue to be suitably small). 
In the collision term,
the real-valued quantity $\collision_{i,j,t,n} J_n$ is the instantaneous change in the value of  
$$
\sum_{i \in I_{q(t)}} J_n \big( x_i, t \big) \mathbf{1}_{m_i(t) = n}
$$
that is caused by the collision of particles $(x_i,m_i)$ and $(x_j,m_j)$ at time $t$. As such, it has the expression
\begin{eqnarray*}
\collision_{x_i,x_j,t,n} J_n & = & 
\tfrac{m_i}{m_i + m_j} J_n(x_i,t)  \mathbf{1}_{m_i(t) + m_j(t) = n} \, +  \, \tfrac{m_j}{m_i + m_j} J_n(x_j,t)  \mathbf{1}_{m_i(t) + m_j(t) = n} \\ 
 & & \quad - \, \, J_n(x_i,t) \mathbf{1}_{m_i(t) = n} \,  -  \,     J_n(x_j,t)  \mathbf{1}_{m_j(t) = n}  \, :
\end{eqnarray*}
two gain terms arise from the appearance of a new particle at one or other of the locations of the disappearing pair, and two loss terms correspond to the disappearance of each element of this pair.

\subsubsection{Introducing microscopic candidate densities}
Our plan for the kinetic limit derivation of the Smulochowski PDE is
 to argue that~(\ref{syspdeweak}) emerges in a suitable sense when we take the high $N$ limit, with the martingale term $M_T$ vanishing in this limit. To implement this plan, we introduce {\it microscopic candidate densities} $f^{\eps,\delta}_{n,t}:\R^d \to [0,\infty)$ of mass $n$ particles at time $t$ under the law $\PP_N$. Here, $\delta > 0$ is a fixed positive quantity, while the interaction radius $\eps$ is determined from $N$ as usual by $N \eps^{d-2} = Z$. The candidate density is given by
\begin{equation}\label{ecandidate}
 f_{n,t}^{\eps,\delta}(u)  = \eps^{d-2} \sum_{i \in I_{q(t)}} \delta^{-d} \eta \big( \tfrac{x_i - u}{\delta} \big)  \mathbf{1}_{m_i(t) = n} \, , \, \, u \in \R^d \, ,
\end{equation}
where $\eta: \R^d \to [0,\infty)$ is a smooth compactly supported function for which $\int_{\R^d} \eta(x) \dd x = 1$.
That is, $f_{n,t}^{\eps,\delta}(u)$ is a statistic reporting a smoothed count of the number of particles in a small macroscopic region about the point $u \in \R^d$ at time $t$ in the model $\PP_N$; note that the time zero microscopic candidate density $f_{n,0}^{\eps,\delta}(u)$ has a high $N$ pointwise limit which as a function of $u$ is given by the convolution of the initial condition $h_n$ of~(\ref{syspde}) and $\eta^\delta(\cdot) = \delta^{-d} \eta\big(\cdot/\delta\big)$. Taking a $\delta \searrow 0$ limit after this limit, we see that, at time zero at least, the appropriate initial condition $h_n(u)$ is obtained at all $u \in \R^d$.

Our aim is to argue that something similar happens at all later times $t > 0$. To do so, we will replace the various terms appearing in the expectation value of equation~(\ref{emicro}) with approximating terms expressed in terms of the microscopic candidate densities, and then take the high $N$ and then the low $\delta$ limit. If we are to reach~(\ref{syspdeweak}) as a result, we will need to understand that the new terms approximate the old ones appropriately.

\subsubsection{Replacing old terms by new: simple cases}

The first term  $\e^{d-2} \sum_{i \in I_{q(T)}} J_n \big( x_i, T \big) \mathbf{1}_{m_i(T) = n}$ has a simple counterpart expressed in the fashion we seek: $\int_{\R^d} J_n(x,T) f_{n,T}^{\eps,\delta}(x) \dd x$. Nor is any nontrivial estimate needed to find a suitable bound on the difference of the terms in this case. Indeed, the two expressions differ by
$$
 \eps^{d-2} \sum_{i \in I_{q(T)}} \Big( J_n(x_i,T) - \int_{\R^d} J_n(y,T) \delta^{-d}  \eta \big( \tfrac{x_i - y}{\delta} \big) \dd y \Big)   \mathbf{1}_{m_i(T) = n} ,
 $$
 which in absolute value is at most $Z \delta \dist \nabla J_n \dist_\infty \leq C \delta$, since total particle number at time $T$ is less than the initial total $N = Z \eps^{2-d}$. 
 
 The first and second terms on the right-hand side of (\ref{emicro}) similarly have counterparts  
 $$
 \int_{\R^d \times [0,T)} \tfrac{\partial J_n}{\partial t}(x,t) \cdot f_{n,t}^{\eps,\delta}(x) \dd x \dd t \quad \textrm{     and    } \quad 
  \int_{\R^d \times [0,T)} d(n) \Delta J_n(x,t) \cdot f_{n,t}^{\eps,\delta}(x) \dd x \dd t \, .
  $$
  Each pair of term and counterpart likewise has a difference which in absolute value is deterministically bounded above by some constant multiple of $\delta$. 
  
\subsubsection{Replacing the coagulation term by using the Stosszahlansatz}  
  
  Given the form of the coagulation term present in~(\ref{syspdeweak}), there is a clear candidate for the form of the term which will form a counterpart to the interaction term appearing in the third line of~(\ref{emicro}): namely, $\int_{\R^d \times [0,T)}   J_n(x,t) \big( Q_1^n \big( f_n^{\eps,\delta}(x) \big) -  Q_2^n \big( f_n^{\eps,\delta}(x) \big) \big) \dd x \dd t$. In stark contrast to the other cases, proving that the replacement of the collision term with this counterpart involves a suitably small error is a major undertaking. We now state the key estimate in this regard, a proposition which we will sometimes call the Stosszahlansatz. Recall that the coefficients $\beta:\N^2 \to (0,\infty)$ are specified in~(\ref{recp}).
  \begin{prop}\label{propsz}
   For each $n,m \in \N$, we have that 
  \begin{eqnarray*}
    & & \eps^{d-2} \E_N \int_0^T \sum_{i,j \in I_{q(t)}} \alpha(m_i,m_j) V_\eps(x_i - x_j) J_n \big( x_i, t \big) \mathbf{1}_{m_i(t) = n,m_j(t) = m} \\
    & = & \beta(n,m) \int_0^T \int_{\R^d} J_n(x,t) f_n^{\eps,\delta}(x,t) f_m^{\eps,\delta}(x,t) \, \dd x \dd t \, \, + \, {\rm Err}_{n,m}(\eps,\delta) \, , 
  \end{eqnarray*}
   where the error ${\rm Err}_{n,m}$ satisfies 
   $$
    \lim_{\delta \searrow 0} \limsup_{\eps \searrow 0} \sum_{m \in \N} \E_N \big\vert {\rm Err}_{n,m}(\eps,\delta) \big\vert = 0 \, . 
   $$
  \end{prop}
  
  Setting $J_n = 1$ for ease of description, note that the integral on the left-hand side is the cumulative rate at which particle pairs of masses $n$ and $m$ are liable to coagulate during all of $[0,T]$; by the relation~(\ref{e.intrange}) and the anticipated survival of a positive fraction of particles at any given positive time, we see that the normalization $\eps^{d-2}$ on the left-hand side is chosen so that the overall expression is of unit order in the high $N$ limit.
  Proposition~\ref{propsz} asserts that this expression is closely approximated by the integral over space-time of the product of the microscopic candidate densities multiplied by the constant coefficient $\beta(n,m)$. As such, this $\beta(n,m)$ is a macroscopic coagulation propensity of pairs of particles of these masses.

Proposition~\ref{propsz} is an expression of the type of precollisional particle independence that we discussed for elastic billiards in Subsection~\ref{s.loschmidt}; here, this independence is manifested by the presence of the product $f_n^{\eps,\delta}(x,t) f_m^{\eps,\delta}(x,t)$. 
  
We now explain how Proposition~\ref{propsz} may be invoked to show that the coagulation term in (\ref{emicro}) is suitably approximated by its counterpart. Recall that the instantaneous change $\collision_{x_i,x_j,t,n} J_n$ is comprised of four terms: two gain terms and two loss terms. Consider the third of these terms, which is the first loss term. This term expresses the instantaneous loss of terms $J_n(x_i,t)$ due to the collision at time $t$ of $x_i$ with some other particle $x_j$. This other particle may have any mass $m_i \in \N$. Writing this term as a sum over that mass, we obtain that the term equals
$$
 - \, \eps^{d-2} \, \E_N \int_0^T \sum_{m=1}^\infty  \sum_{i,j \in I_{q(t)}}  \alpha(m_i,m_j) V_\eps \big( x_i - x_j \big)  J_n(x_i,t) \mathbf{1}_{m_i(t) = n, m_j(t) = m}   \, \dd t \, .
$$
Applying Proposition~\ref{propsz}, we find that the term equals
$$
 - \sum_{m=1}^\infty \beta(n,m) \int_0^T \int_{\R^d} J_n(x,t) f_{n,t}^{\eps,\delta}(x) f_{m,t}^{\eps,\delta}(x) \, \dd x \dd t \, \, + \, {\rm Err}_n(\eps,\delta) \, ,
$$
where this error term, after the sum over $m \in \N$, is known to satisfy 
\begin{equation}\label{eerror}
\lim_{\delta \searrow 0} \limsup_{\eps \searrow 0} \E_N \big\vert {\rm Err}_n(\eps,\delta) \big\vert = 0 \, .
\end{equation}
Exactly the same considerations apply to the second of the loss terms because the two terms are equal due to the symmetry of the interaction kernel $V$. 

The comparable estimate for each of the gain terms is slightly easier to handle, because a particle of mass $n$ may be produced by only finitely many mass pairs -- $(1,n-1)$, $(2,n-2)$, \ldots, $(n-1,1)$ -- rather than the infinite number of choices -- $(n,1)$,$(n,2)$,\ldots -- which may cause such a particle to disappear. Regarding the first term, an application of Proposition~\ref{propsz} yields that
$$
 \eps^{d-2} \E_N \int_0^T \sum_{m=1}^{n-1}  \sum_{i,j \in I_{q(t)}}  \alpha(m_i,m_j) V_\eps \big( x_i - x_j \big) \tfrac{m}{n} J_n(x_i,t) \mathbf{1}_{m_i(t) = m, m_j(t) = n- m}   \, \dd t  
$$  
equals
$$
 \sum_{m=1}^n \beta(m,n-m) \int_{\R^d \times [0,T)} J_n(x,t) \tfrac{m}{n} f_{n,t}^{\eps,\delta}(x) f_{m,t}^{\eps,\delta}(x) \, \dd x \dd t \, \, + \, {\rm Err}_n(\eps,\delta) \, ,
$$
where likewise the error satisfies (\ref{eerror}). 
The second gain term differs only in that  $\tfrac{n- m}{n}$ replaces  $\tfrac{m}{n}$; thus, the total gain term satisfies the same statement with this term omitted. 

\subsubsection{The martingale term is replaced by zero}
The martingale term $M_T$ in (\ref{emicro}) is treated by arguing that it is typically small in absolute value:
\begin{prop}\label{p.mart}
There exists $C > 0$ such that, for each $N \in \N$, $\sup_{T \in (0,\infty)} \E_N M(T)^2 \leq C \eps^{d-2}$.
\end{prop}
The martingale term is in a sense much smaller than the collision term treated by the Stosszahlansatz: it vanishes before the low $\delta$ limit is even taken. The ideas in the proof of Proposition~\ref{p.mart} are already in large part seen in the proof of the more substantial Proposition~\ref{propsz} and we will not explain their specifics; it is in Section~$5$ that the martingale term is treated in~\cite{HR3d}. 
\subsubsection{The counterpart to the PDE using microscopic candidate densities}
By using Propositions \ref{propsz} and \ref{p.mart} and the other, easier,  estimates, we are able to replace each term in~(\ref{emicro}) with its counterpart, expressed in terms of the microscopic candidate densities, and obtain the following bound on the error in the resulting near identity:

\begin{eqnarray}
 & & \int_{\R^d} J_n(x,T) f_{n,T}^{\eps,\delta}(x) \, \dd x - \int_{\R^d} J_n(x,0) f_{n,0}^{\eps,\delta}(x) \, \dd x \label{emicrocand} \\
  & = &  
   \int_{\R^d \times [0,T)} \tfrac{\partial J_n}{\partial t}(x,t) f_{n,t}^{\eps,\delta}(x) \, \dd x \dd t
    \, + \,
  \int_{\R^d \times [0,T)} d(n) \Delta J_n(x,t) f_{n,t}^{\eps,\delta}(x) \, \dd x \dd t  \nonumber \\
  & & \quad \quad \quad \, + \, 
   \sum_{m=1}^n \beta(m,n-m) \int_{\R^d \times [0,T)} J_n(x,t)  f_{n,t}^{\eps,\delta}(x) f_{m,t}^{\eps,\delta}(x) \, \dd x \dd t \nonumber \\
    & & \qquad \qquad
 \, - \,   2 \sum_{m=1}^\infty \beta(n,m) \int_{\R^d \times [0,T)} J_n(x,t) f_{n,t}^{\eps,\delta}(x) f_{m,t}^{\eps,\delta}(x) \, \dd x \dd t \, \, + \, {\rm Err}_n(\eps,\delta) \, . \nonumber
\end{eqnarray}
where the error satisfies~(\ref{eerror}) because each of the errors used in the five estimates which we applied does.  
\end{subsection}
\begin{subsection}{Taking the limit to obtain the Smoluchowski PDE}\label{s.takingthelimit}
 Our approximate identity~(\ref{emicrocand}) closely resembles  the equation~(\ref{syspdeweak}) satisfied by a weak solution of the Smoluchowski PDE: we simply replace the solution of the latter with the microscopic candidate densities, and add in the error term, to obtain the former. However, to pass to~(\ref{syspdeweak}) from~(\ref{emicrocand}) in the limit of low $\eps$ followed by low $\delta$, we must carry out a short further analysis that will make use of some additional information about the microscopic models~$\PP_N$. 
 
 \subsubsection{The approximate identity rewritten using empirical measures}
 Recall from Theorem~\ref{thmo} that in fact we express approximation by $\PP_N$ for high $N$ of the Smoluchowski PDE by using the empirical measures $\empmeas_N$ valued in space-mass-time $\R^d \times \N \times [0,\infty)$. In Section~\ref{s.empirical}, we let $\empmeas_{N,n}$ denote the mass $n$ marginal of $\empmeas$, the empirical measure in space-time for particles of mass $n$,  for each $n \in \N$. That is, 
 $$
  \empmeas_{N,n} 
= \eps^{d-2} \sum_{i\in I_{q(t)}}\d_{(x_i , t )} \mathbf{1}_{m_i = n} \, \dd t \, .
$$

On $\PP_N$, the microscopic candidate densities are expressed in terms of the empirical measures by
$f_{n,t}^{\eps,\delta} (x) =  \big( \empmeas_{N,n} * \eta^\delta  \big) (x,t)$ for all $(x,t) \in \R^d \times [0,\infty)$, where the convolution is in the space variable.

Recall further that in Theorem~\ref{thmo}, $\mathcal{M}$
is the space of measures  $\mu$ on $\R^d \times \N \times [0,\infty)$
such that $\mu \big( \R^d \times \N \times [0,T] \big) \in [0,TZ]$ for each $T \geq 0$; and that $\mu_N$ is a random measure taking values in the space $\mathcal{M}$ whose law we denote by $\mathcal{P}_N$. By equipping $\mathcal{M}$
with the topology of vague convergence, we give meaning to the notion of weak convergence of a sequence of measures such as $\big\{ \mathcal{P}_N : N \in \N \big\}$. Now, in the theorem, $\mathcal{P}$ is supposed to be a weak limit point of this sequence (the existence of which is assured by the tightness of this sequence of measures). Thus, there is a subsequence $\big\{ N_i: i \in \N \big\}$ of natural numbers for which $\mathcal{P}_{N_i}$ converges weakly to~$\mathcal{P}$.

Now since under $\PP_N$, we have that 
$f_{n,t}^{\eps,\delta} (x) =  \big( \empmeas_{N,n} * \eta^\delta  \big) (x,t)$ for all $(x,t) \in \R^d \times [0,\infty)$, we see that (\ref{emicrocand})
asserts that 

\begin{eqnarray}
 & & \int_{\R^d} J_n(x,T) \big( \empmeas_n * \eta^\delta  \big) (x,T) \, \dd x - \int_{\R^d} J_n(x,0)  \big( \empmeas_n * \eta^\delta  \big) (x,0) \, \dd x  \label{emicrocandlimit} \\
  & = &  
   \int_{\R^d \times [0,T)} \tfrac{\partial J_n}{\partial t}(x,t)  \big( \empmeas_n * \eta^\delta  \big) (x,t) \, \dd x \dd t
    \, + \,
  \int_{\R^d \times [0,T)} d(n) \Delta J_n(x,t) \big( \empmeas_n * \eta^\delta  \big) (x,t) \, \dd x \dd t  \nonumber \\
  & & \quad \quad \quad \, + \, 
   \sum_{m=1}^n \beta(m,n-m) \int_{\R^d \times [0,T)} J_n(x,t)   \big( \empmeas_n * \eta^\delta  \big) (x,t)  \big( \empmeas_m * \eta^\delta  \big) (x,t) \, \dd x \dd t \nonumber \\
    & & \qquad \qquad
 \, - \,   2 \sum_{m=1}^\infty \beta(n,m) \int_{\R^d \times [0,T)} J_n(x,t)  \big( \empmeas_n * \eta^\delta  \big) (x,t)  \big( \empmeas_m * \eta^\delta  \big) (x,t) \, \dd x \dd t \, \, + \, {\rm Err}_n(\eps,\delta) \, . \nonumber
\end{eqnarray}
 Here, $\big\{ \mu_n: n \in \N \big\}$
be a sequence of random measures, with $\mu_n$ supported on $\R^d \times \{ n \} \times [0,\infty)$ for each $n \in \N$, such that the $\mathcal{M}$-valued random measure $\sum_{n=1}^\infty \mu_n$ is~$\mathcal{P}_N$-distributed. The error is now a real-valued function  defined on the space $\mathcal{M}$ which in view of~(\ref{eerror}) is seen to satisfy 
$\lim_{\delta \searrow 0} \limsup_{\eps \searrow 0} \E_{\mathcal{P}_N} \big\vert {\rm Err}_n(\eps,\delta) \big\vert = 0$, where here we write $\E_{\mathcal{P}_N}$
for expectation with respect to the law $\mathcal{P}_N$.
We may now take a high $N$ limit of this identity along the subsequence $\big\{ N_i: i \in \N \big\}$. We learn that the identity continues to hold, but where now the $\mathcal{M}$-valued random measure $\sum_{n=1}^{\infty}\empmeas_n$ has the distribution of~$\mathcal{P}$. 
What estimate does the error term satisfy when this limit is taken? Note  that in~(\ref{emicrocandlimit}) the error  ${\rm Err}_n(\eps,\delta)$ may be viewed as a function of $\sum_{n=1}^\infty \mu_n \in \mathcal{M}$. Indeed, this formula for the error is a continuous function of   $\sum_{n=1}^\infty \mu_n$ given that $\mathcal{M}$ carries the vague topology. As such, in the $N_i \to \infty$ limit (in which $\e \searrow 0$),
the error term loses its $\e$-dependence and now satisfies $\lim_{\delta \searrow 0} \mathbb{E}_{\mathcal{P}} \big\vert {\rm Err}_n(\delta) \big\vert = 0$. Here, $\mathbb{E}_{\mathcal{P}}$ denotes expectation with respect to the law~$\mathcal{P}$.

\subsubsection{Preparing for the final step towards the PDE: uniform integrability}
In order to conclude the proof of Theorem~\ref{thmo}, two steps are needed. First, 
\begin{prop}\label{prop.firststep}
Let $\big\{ \mu_n: n \in \N \big\}$
be a sequence of random measures, with $\mu_n$ supported on $\R^d \times \{ n \} \times [0,\infty)$ for each $n \in \N$, such that $\sum_{n=1}^\infty \mu_n$ is~$\mathcal{P}$-distributed. Then, almost surely, we may express each $\mu_n$ in the form $f_n \,  \dd x \times \delta_n \times \dd t$, where  $f_n:\R^d \times [0,\infty) \to [0,\infty)$.
\end{prop}
Second, we must argue  that the collection $\{ f_n: n \in \N \}$ solves~(\ref{syspdeweak}).
Given the first step, it is the taking of the low~$\delta$ limit in~(\ref{emicrocandlimit}) which will yield the second. 
However, to successfully carry out this limit, an extra piece of information will be needed, namely, for each fixed $n \in \N$, the uniform integrability of the family $f_n*\eta^\delta$ as $\delta$ ranges over $(0,1)$.
The next proposition is sufficient in this regard.
\begin{prop}\label{prop.secondstep}
There exists a sequence $\big\{ k_n: n \in \N \big\}$ of positive constants such that the collection of random functions provided by Proposition~\ref{prop.firststep} satisfies
 $\dist f_n \dist_{L^\infty(\R^d \times [0,\infty))} \leq k_n$ for each $n \in \N$ almost surely.
\end{prop}
We now confirm that these two elements applied to~(\ref{emicrocandlimit}) are enough to yield Theorem~\ref{thmo}.

\subsubsection{Taking the final step: Proof of Theorem~1.1
.} By Proposition~\ref{prop.firststep},~(\ref{emicrocandlimit}) after the $\e \searrow 0$ limit is taken holds with $f_n$ in place of $\mu_n$ for each $n \in \N$. The uniform boundedness provided by Proposition~\ref{prop.secondstep}
and the Lebesgue differentiation theorem imply that $f_n * \eta^\delta$ converges pointwise to $f_n$ almost everywhere on $\R^d \times [0,\infty)$ for each $n \in \N$. Recall that our test functions~$J_n$ are compactly supported in space-time. Using this alongside the same uniform boundedness as above, we may apply the  dominated convergence theorem to find that each integral appearing in the identity converges to its counterpart where $f_n * \eta^\delta$ is replaced by $f_n$. Recall that $\lim_{\delta \searrow 0} \mathbb{E}_{\mathcal{P}} \big\vert {\rm Err}_n(\delta) \big\vert = 0$. Thus the error term converges to zero in probability in the low $\delta$ limit. It thus converges to zero almost surely along a subsequence of $\delta \searrow 0$. In this way, we see that~(\ref{syspdeweak}) holds $\mathcal{P}$-almost surely. \qed

Our remaining task then is to prove Propositions~$4.i$ for integer~$i$ satisfying $1 \leq i \leq 4$. Figure~\ref{f.flowchart} depicts how the derivations will be presented.

  \begin{figure}
    \begin{center}
      \includegraphics[width=1.0\textwidth]{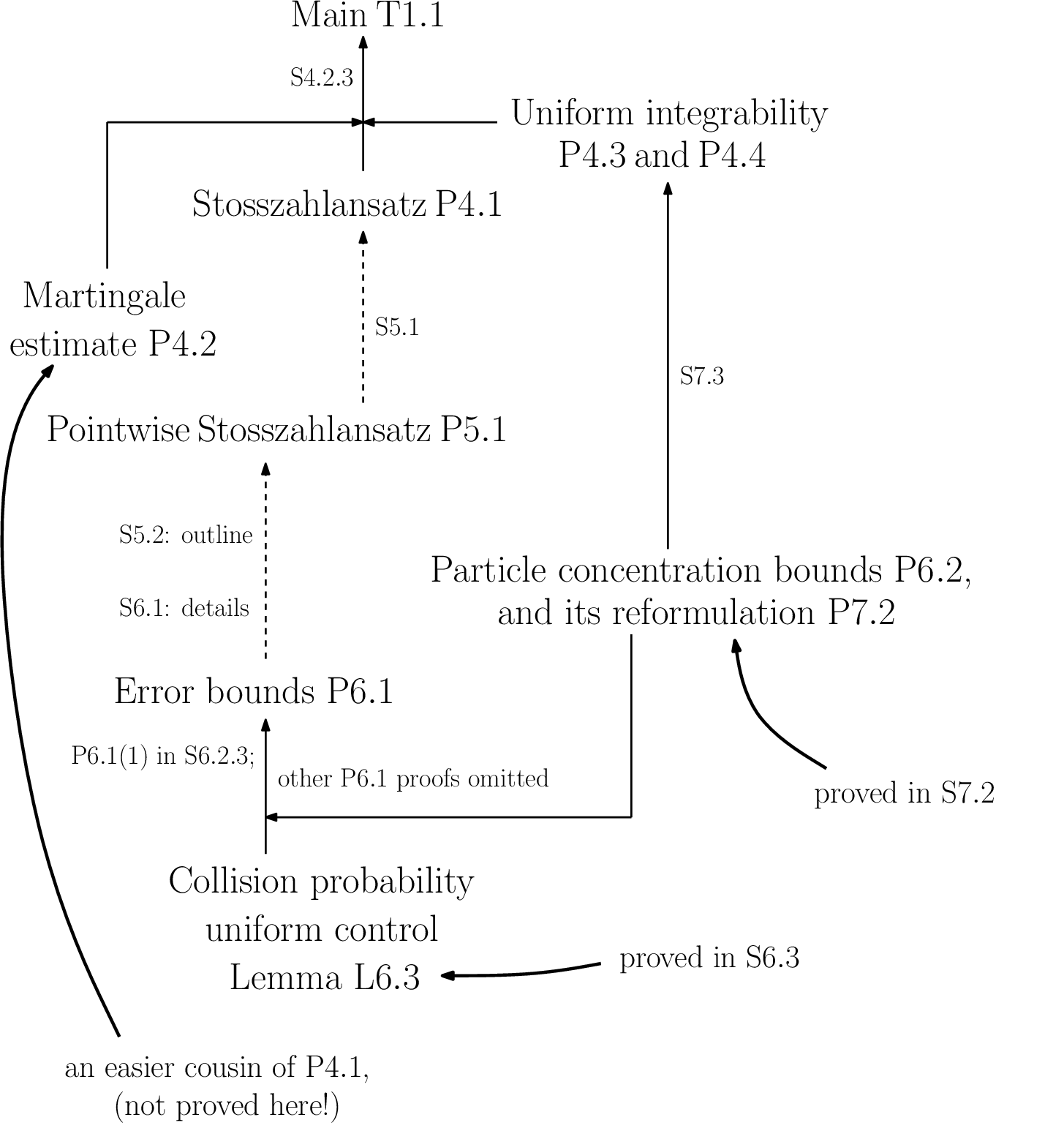}
    \end{center}
    \caption{The road ahead: the structure of the proof of the main elements of the main theorem. $T =$ Theorem, $P =$ Proposition, $L =$ Lemma, and $S=$ Section (or Subsection). An arrow indicates that one result is used to prove the other; a dashed arrow indicates that some details are omitted in the proof. The $S$-labellings of arrows indicate the section where the derivation takes place.}\label{f.flowchart}
  \end{figure}

\end{subsection}

\section{An outline of the proof of the Stosszahlansatz}\label{s.outlinestoss}
Here we explain in outline how we will prove Proposition~\ref{propsz}.
\subsection{Coagulation propensity, and particle pairs at small macroscopic distance}\label{s.coagprop}
For $z \in \R^d$ and $n,m \in \N$, define under the law $\PP_N$ the stochastic process $Q_z = Q_{z,n,m}:[0,\infty) \to \R$ whose value at time $t \in [0,\infty)$ is given by 
\begin{equation}\label{defq}
 \tfrac{1}{2} \eps^{d-2} \sum_{i,j \in I_{q(t)}} \alpha(m_i,m_j) V_\eps(x_i - x_j + z) J_n \big( x_i, t \big) \mathbf{1}_{m_i(t) = n,m_j(t) = m} \, .
\end{equation}
In seeking to prove the Stosszahlansatz, then, it is our aim to show that $\E_N \int_0^T Q_0(t) \dd t$
is close to a $\beta(n,m)$ multiple of the time-integrated product of microscopic candidate densities for particles of mass $n$ and $m$. Since the microscopic coagulation density is binary in nature, it is unsurprising that $\E_N \int_0^T Q_0(t) \dd t$, the cumulative rate of coagulation between pairs of particles of such masses (at least if $J_n = 1$),
should be approximated by the time integral of such a product of empirically defined densities. We have already explained heuristically in Section~\ref{s.torus} why we might expect the macroscopic coagulation propensity $\beta(n,m)$ to have the form~(\ref{recp}). The challenge now is to find a rigorous means of approximating  $\E_N \int_0^T Q_0(t) \dd t$; as we outline this approach, we will see an alternative explanation for the formula~(\ref{recp}) emerge. 

Consider for a moment the expression $Q_z$, where $z \in \R^d$ is a small macroscopic quantity; which is to say, $z$ is fixed at a given small value as we take a high $N$ (or low $\eps$) limit. We see that the quantity  $\E_N \int_0^T Q_z(t) \dd t$ is a time-averaged count of all instances of pairs of particles, of mass~$n$ and mass $m$, for which  the mass $m$ particle lies in the tiny $\eps$-ball whose centre is displaced from the mass $n$ particle by the small quantity $z$. Such instances at any given moment of time are weighted by the factor $\alpha(m_i,m_j) V_\eps(x_i - x_j + z)$; for later convenience, it is useful to also define $\hat{Q}_z$, where in the formula for $Q_z$, we replace $V_\eps(x_i - x_j + z)$ by $\hat{V}_\eps (x_i - x_j + z)$. Here, $\hat{V}_\eps(\cdot) = \eps^{-2}\hat{V}(\cdot/\eps)$, where $\hat{V}:\R^d \to [0,\infty)$ is a smooth and compactly supported function. (In fact, we will use two such variants, also writing $\overline{Q}_z$ when $\overline{V}_\eps$ replaces $V_\eps$.) The quantity $\E_N \int_0^T \hat{Q}_z(t) \dd t$ qualitatively meets the same description as does $\E_N \int_0^T Q_z(t) \dd t$, a time-averaged count of instances of $z$-displacements of particle pairs. The next assertion shows that  $\E_N \int_0^T Q_0(t) \dd t$
is well approximated by an appropriately weighted count of particle pairs at small macroscopic distance $z$:
\begin{prop}\label{propqqbar}
For $n,m \in \N$, recall that $u_{n,m}:\R^d \to [0,1]$ is specified in~(\ref{pdew}), and let $\overline{V} = V \big( 1 - u_{n,m} \big)$. Then 
\begin{equation}\label{e.overl}
   \int_0^T Q_0(t) \, \dd t =   \int_0^T \overline{Q}_z(t) \, \dd t \, \, + \, {\rm Err}_{n,m}(\eps,z) \, .
\end{equation}

Regarding the error term: defining ${\rm Err}_\delta$ to be the supremum over $z \in \R^d$ for which $\vert z \vert =\delta$ of $$\limsup_{\eps \searrow 0} \sum_{m \in \N} \E_N \big\vert {\rm Err}_{n,m}(\eps,z) \big\vert \, ,
$$
we have that 
${\rm Err}_\delta \to 0$ as $\delta \to 0$.
\end{prop}

This result may be called the {\em pointwise} Stosszahlansatz, because the small vector~$z$ may be treated as fixed rather than used as a variable for averaging. Indeed, we easily deduce  the Stosszahlansatz Proposition~\ref{propsz} from this pointwise version by averaging over~$z$: 

\medskip

\noindent{\bf Sketch of proof of Proposition~\ref{propsz}.} 
We verify the statement only in the case $J_n = 1$; the general case invokes a simple additional estimate.

To derive  Proposition~\ref{propsz}, we begin by averaging the information in Proposition~\ref{propqqbar} over small macroscopic $\delta$. 
In what follows, the relation $f \almeq g$ asserts that $f(n,m,\delta)$ and $g(n,m,\delta)$ are random variables on $\PP_N$ 
such that, for each $n \in \N$, 
$$
\lim_{\delta \searrow 0} \limsup_{\eps \searrow 0} \sum_{m \in \N} \E_N \vert f_{n,m,\delta} - g_{n,m,\delta} \vert = 0 \, .
$$ 

For $\delta > 0$, write $\eta^\delta: \R^d \to [0,\infty)$ for $\eta_\delta = \delta^{-d}\eta(\cdot/\delta)$.
 Proposition~\ref{propqqbar}, the definition of $\overline{Q}$ and substitutions $x_i - z_1 = \omega_1$ and $x_j - z_2 = \omega_2$ imply then that
\begin{eqnarray*}
 & &  \int_0^T Q_0(t) \, \dd t \\
  & \almeq & \int_0^T \int_{\R^d \times \R^d} \overline{Q}_{z_2 - z_1}(t) \eta^\delta(z_1) \eta^\delta(z_2) \, \dd z_1 \dd z_2 \dd t \\
  & = & \eps^{d-2}  \alpha(n,m)  \int_0^T \int_{\R^d \times \R^d}  \sum_{i,j \in I_{q(t)}} \overline{V}_\eps\big( (x_i - z_1) - (x_j - z_2)  \big) 
  \mathbf{1}_{m_i(t) = n,m_j(t) = m}  \eta^\delta(z_1) \eta^\delta(z_2) \, \dd z_1 \dd z_2 \dd t \\
    & = & \eps^{d-2}  \alpha(n,m)  \int_0^T \int_{\R^d \times \R^d}  \sum_{i,j \in I_{q(t)}}
    \overline{V}_\eps\big( \omega_1  - \omega_2  \big) 
    \mathbf{1}_{m_i(t) = n,m_j(t) = m}  \eta^\delta(x_i - \omega_1) \eta^\delta(x_j - \omega_2) \, \dd \omega_1 \dd \omega_2 \dd t  \, . 
 \end{eqnarray*}
 The virtue of this last expression is that the two particle sums may be passed inside to yield the microscopic candidate densities. Indeed, the expression equals
\begin{equation}\label{e.doubint}    
\eps^{d-2}  \alpha(n,m)  \int_0^T \int_{\R^d \times \R^d}
    \overline{V}_\eps\big( \omega_1  - \omega_2  \big) 
    f_n^{\eps,\delta}(\omega_1,t) f_m^{\eps,\delta}(\omega_2,t) \,  \dd \omega_1 \dd \omega_2 \dd t \, . 
\end{equation}
Note that the double integral in $(\omega_1,\omega_2)$ is almost on the diagonal, because $\overline{V}_\eps$ is supported in the $\eps$-ball. When $\omega_1,\omega_2 \in \R^d$ satisfy $\dist \omega_1 - \omega_2 \dist \leq \e$, 
we have that $\big\vert \eta \big( \tfrac{x_i - \omega_2}{\delta} \big)  - \eta \big(  \tfrac{x_i - \omega_1}{\delta} \big) \big\vert \leq \e \delta^{-1} \dist \nabla \eta \dist_\infty$, so that 
$\big\vert f_n^{\eps,\delta}(\omega_2,t) - f_n^{\eps,\delta}(\omega_1,t) \big\vert \leq \eps^{d-1} \delta^{-d -1}  \dist \nabla \eta \dist_\infty$. 

Thus,  at the expense of an error that is small in the sense of the $\almeq$ relation,  $f_m^{\eps,\delta}(\omega_2,t)$ may be replaced by  $f_m^{\eps,\delta}(\omega_1,t)$ in~(\ref{e.doubint}); this done, 
the $\omega_2$ integral may be detached, so that we see that~(\ref{e.doubint}) satisfies
$$
      \almeq  \eps^{d-2}  \alpha(n,m)  \int_{\R^d }
    \overline{V}_\eps \big( x  \big) \dd x  \int_0^T  
     \int_{\R^d} f_n^{\eps,\delta}(\omega,t) f_m^{\eps,\delta}(\omega,t) \, \dd \omega  \dd t \, ,
$$
which since $\overline{V}_\eps(\cdot) = \eps^{-2} \overline{V}(\cdot/\eps)$ equals
$$
       \eps^{2(d-2)}  \alpha(n,m)  \int_{\R^d }
    \overline{V} \big( x  \big) \dd x  \int_0^T  
     \int_{\R^d} f_n^{\eps,\delta}(\omega,t) f_m^{\eps,\delta}(\omega,t) \, \dd \omega  \dd t \, ,
$$

Combining the above estimates, we confirm that Proposition~\ref{propsz} holds with 
$$
    \qquad \qquad  \qquad \qquad  \qquad \qquad  \beta(n,m) = \alpha(n,m) \int_{\R^d} \overline{V}(x) \, \dd x \, .  \qquad \qquad  \qquad \qquad \qquad \qquad \Box
$$
\subsection{An outline of the proof of Proposition 5.1 
}
 In an attempt to find a convenient representation of the quantity $\int_0^T Q_z(t) \, \dd t$, both when  $z \in \R^d$ is zero and when it is non-zero and small, we define a $z$-dependent random variable $S_z$ under $\PP_N$ for which the action of the free motion operator $\freem$ on $S_z$ produces, among others, the term $Q_z$. 
  For each pair $(n,m) \in
\mathbb{N}^2$, we define $\phi_{n,m}^{\epsilon}:
\mathbb{R}^d \to (0,\infty)$ so that 
$$
- \Delta \phi^\epsilon_{n,m} (x) =  \ald  \epsilon^{-d} V(x/\epsilon) \, ,
$$
subject to $\lim_{x \to \infty} \phi_{n,m}^{\epsilon}(x) = 0$.
We then define a non-negative stochastic process $S_z: [0,\infty) \to [0,\infty)$ on $\PP_N$: for each $t \geq 0$, we set 
\begin{equation}\label{exz}
S_z(t) (q) = \epsilon^{2(d-2)} \sum_{i,j \in I_{q(t)}}
\phi_{n,m}^{\epsilon} (x_i- x_j + z) J_n (x_i,t)
\mathbf{1}_{ m_i = n , m_j = m } \, .
\end{equation}
The action $\freem(S_z - S_0)$ of the free motion operator on $S_z - S_0$ is itself a random variable on $\PP_N$ which maps non-negative time $t \in [0,\infty)$ to $\R$. The term $Q_z(t) - Q_0(t)$ appears in the expression 
$- \freem(S_z - S_0)(t)$, in the case where each of the derivatives in the Laplacian operator falls on $\phi_{n,m}^{\epsilon}$ rather than on the test function $J_n$. (Note the minus sign attached to $\freem(S_z - S_0)(t)$; it arises from our choice that the functional $S_z$ be positive rather than negative.)

For $T > 0$, consider the $\PP_N$-almost sure identity 
\begin{eqnarray}\label{hac}
    \big( S_z - S_0 \big) \big( T \big) & = &
    \big( S_z - S_0 \big) \big( 0 \big) +
\int_{0}^{T}{ \left(\tfrac{\partial}{\partial t}+\freem 
\right)(S_z - S_0) (t) \, \dd t} \\
    & & \qquad + \,
\int_{0}^{T}{ \collop (S_z - S_0) (t) \, \dd t} \, + \, M_T \, , \nonumber
\end{eqnarray}
and note that the process $\big\{ M_T : T \geq 0 \big\}$ is a $\PP_N$-martingale.
As we have noted, each of the terms  $- \int_0^T Q_z(t) \dd t$ and $\int_0^T Q_0(t) \dd t$ appears in the free motion term on the right-hand side.
The quantity  $\int_0^T Q_0(t) \dd t$ remains of unit order in the low $\eps$ limit, as we discussed after the statement of the Stosszahlansatz Proposition~\ref{propsz}. For similar reasons,  $\int_0^T Q_z(t) \dd t$ may be expected to have this property for any given $z \in \R^d$. 
Suppose for a moment that it were the case that all the other terms appearing in~(\ref{hac}) were of smaller order, as a low $\eps$ and then low $z$ limit is taken. More precisely, suppose that, after the removal of the two terms mentioned above, the remainder satisfies the estimate on the error ${\rm Err}_{n,m}$ given in Proposition~\ref{propqqbar}. Then in fact Proposition~\ref{propqqbar} would hold, but with the term $\overline{V}$ on the right-hand side of~(\ref{e.overl}) replaced by $V$. Reviewing the proof of Proposition~\ref{propsz} from Proposition~\ref{propqqbar}, we would find that the Stosszahlansatz indeed holds, but with the formula $\beta = \alpha \int_{\R^d} V(x) {\rm d} x$ rather than~(\ref{recp}). In other words, the reasoning that there are no further unit-order terms appearing in~(\ref{hac}) -- which the authors of~\cite{HR3d} believed when first studying this approach -- leads to the fallacious guess $\beta = \alpha  \int_{\R^d} V(x) {\rm d} x$.


The formula's incorrectness means that we should expect some further term in~(\ref{hac}) to remain of unit order as $\eps \to 0$ and then $z \to 0$. 
In Section~\ref{s.torussec}, it was explained that this guess is wrong because of an effect of the curtailment of the interaction clock associated with a particle pair at the moment of the concerned particles' collision.

The further unit-order term in (\ref{hac}) does indeed exist, and its form reflects this mechanism of curtailment of the interaction clock on particle collision. The term arises from the action of
the collision operator $\collop$ on $S_0$. In the expression $\int_0^T \collop S_0(t) \dd t$,
when the form of the collision operator~$\collop$ and of the functional $S_0$ is substituted, a sum is obtained. For each summand, four particles are concerned, in two pairs: two particles in the first pair arise from $\collop$, and it is this pair whose infinitesimal interaction rate is being integrated over the time period $[0,T]$, while the difference of the locations of the particles in the second pair form the argument of $\phi^\eps_{n,m}$ arising from $S_0$. Although each of the two pairs is formed of distinct particles, there may be one or two coincidences between members of the first and of the second pair. When both of these coincidences occur, and the second pair equals the first, the contribution made to  $- \int_0^T \collop S_0(t) \dd t$ by such terms is given by 
\begin{equation}\label{geezee}
     \, \epsilon^{2(d-2)}  \micp(n,m) \int_0^T \sum_{k,l \in
I_{q(t)}}{  V_{\epsilon}\big( x_k - x_l \big)}
\phi_{n,m}^{\epsilon}(x_k -
x_l ) J_n(x_k,t)  \, \dd t \, .
\end{equation}

This term witnesses the abrupt curtailment of the propensity to
coagulate of a pair of particles at that moment when the particles do
coagulate. That it is this term which remains of unit order reflects the role of the microscopic repulsion about a given particle in determining the relation~(\ref{recp}) which is discussed in Section~\ref{s.torus}.

(We mention a possible confusion relating to our convention regarding double sums. It might seem that a factor of one-half should multiply (\ref{geezee}) because such a factor is present in the definition of the collision operator~(\ref{e.collisionoperator}). Recall however that double sums are really over ordered pairs of distinct indices. 
The two pairs at stake in arriving at~(\ref{geezee}) may be labelled $\{ i,j \}$ and $\{ k,l \}$. 
Ordering the first pair $(i,j)$, there are two  ways, $(i,j) = (k,l)$ and $(i,j) = (l,k)$,  that the coincidence of both terms of the first pair with those of the second may occur. This factor of two cancels the one-half from~(\ref{e.collisionoperator}), yielding the expression~(\ref{geezee}).)

It turns out that the sum of the remaining terms in (\ref{hac}) is indeed negligible in that it satisfies 
the estimate that ${\rm Err}_{n,m}(\eps,z)$  does in Proposition~\ref{propsz}. Only the three unit-order terms identified above remain in the limit of low $\eps$ and then $z$. That is, we have found that 
$\int_0^T Q_z(t) \dd t$ differs from
\begin{eqnarray}\label{geeque}
     &  &  \epsilon^{2(d-2)}  \micp(n,m) \int_0^T  \sum_{i,j \in
I_{q(t)}}{ V_{\epsilon}(x_i - x_j)}
\\
    & & \qquad \qquad \qquad \times \
     \Big[
1 +  \phi_{n,m}^{\epsilon} \big( x_i - x_j \big) \Big] J_n(x_i,t)
 \mathbf{1}_{ m_i = n,m_j = m} \, \dd t \, . \nonumber
\end{eqnarray}
by an error of the form ${\rm Err}_{n,m}(\eps,z)$ in Proposition~\ref{propsz}.
Note that 
the `$1$' that appears in the square bracket corresponds to $Q_0$,
and the other term to the unit-order term (\ref{geezee}). 
In the language of Proposition~\ref{propqqbar}, we have learnt that  
\begin{equation}\label{efirstconc}
   \int_0^T \hat{Q}_0(t) \, \dd t =   \int_0^T Q_z(t) \, \dd t \, \, + \, {\rm Err}_{n,m}(\eps,z) \, ,
\end{equation} 
where $\hat{V} = V \big( 1 + \phi_{n,m} \big)$ and $\hat{Q}_0$ is defined by the formula~(\ref{defq}) that specifies $Q_0$ with $\hat{V}$ replacing~$V$. As a check of working, note that, since $\phi_{n,m} \geq 0$, we are asserting that the positive $\int_0^T Q_z(t) \dd t$ exceeds the positive $\int_0^T Q_0(t) \, \dd t$ by a further positive term of the same order. This is consistent with the explanation offered in Section \ref{s.torussec}: we expect $\int_0^T Q_z(t) \dd t$ to exceed $\int_0^T Q_0(t) \dd t$, because the size of the latter term (measuring the cumulative interaction clock of $\eps$-displaced particles) is limited by the disappearance of particles on collision, while the former (measuring a comparable quantity for the much more distant $z$-displayed particles) experiences no such limitation. 

Of course, (\ref{efirstconc}) is not quite the conclusion we sought: to prove Proposition~\ref{propqqbar}, we want to approximate  $\int_0^T Q_0(t)$ by $\int_0^T \overline{Q}_z(t)$, so that the modification $Q \to \overline{Q}$ falls in the $z$-displaced term; but so far we have obtained such a result where the modification is made to the $z = 0$ term. 

In light of this analysis, we may however revisit the approach. Consider a variant $X_z$ of the process~$S_z$: for each $z \in \R^d$, under $\PP_N$, $X_z:[0,\infty) \to \R$ is given by 
\begin{equation}\label{exzo}
X_z(t) (q) = \epsilon^{2(d-2)} \sum_{i,j \in I_{q(t)}}
 u_{n,m}^{\epsilon} (x_i- x_j + z) J_n (x_i,t)
\mathbf{1} \big\{ m_i = n , m_j = m \big\} \, ,
\end{equation}
where here, for each pair $(n,m) \in
\mathbb{N}^2$, we define $u_{n,m}^{\epsilon}:
\mathbb{R}^d \to (0,\infty)$ so that, for $z \in \R^d$, 
\begin{equation}\label{eueps}
 - \, \Delta u^\epsilon_{n,m} (z) =  \ald  \epsilon^{-d} U(z/\epsilon) \, ,
\end{equation}
subject to $\lim_{x \to \infty} u_{n,m}^{\epsilon}(x) = 0$. The function $U: \R^d \to (0,\infty)$ is at yet unspecified; of course the choice $U = V$ would specify the earlier functional $S_z$. Our aim now is to make a different choice of~$U$, for which the solution of the problem~(\ref{eueps}) exists uniquely, and for which the earlier analysis may be carried out in such a way that its conclusion is not~(\ref{efirstconc}) but rather the desired 
\begin{equation}\label{esecondconc}
   \int_0^T Q_0(t) \, \dd t =   \int_0^T \overline{Q}_z(t) \, \dd t \, \, + \, {\rm Err}_{n,m}(\eps,z) \, ,
\end{equation} 
for some $\overline{Q}: \R^d \to [0,\infty)$. Regarding scaling, note that, whenever $U: \R^d \to (0,\infty)$ is such that~(\ref{eueps}) has a unique solution for some $\eps > 0$, then this holds in fact for all $\eps > 0$; indeed, writing $u_{n,m}:\R^d \to (0,\infty)$ for $u_{n,m} = u_{n,m}^1$, we have that, for each $\eps > 0$, and for all $z \in \R^d$,
\begin{equation}\label{euscaling}
u^\eps_{n,m}(z) = \eps^{2 -d} u_{n,m}\big( z/\eps \big) \, .
\end{equation}

In order to find a candidate for $U$ that may make this plan work, we may hope that, for some suitable class of $U$, the earlier discussion continues to apply to the extent that the unit-order terms that survive in the passage of low $\eps$ and then low $z$ are the natural counterparts to the three terms identified there.

We want the analogue of the term $\int_0^T \hat{Q}_0(t)$ to be $\int_0^T Q_0(t)$ in the new calculation. 
Recall that $\hat{Q}_0(t)$ equals~(\ref{geeque}). When we reprise the earlier discussion with $X_z - X_0$ in place of $S_z - S_0$, the counterpart of the expression~(\ref{geeque}) is
\begin{eqnarray}
     &  &  \epsilon^{2(d-2)} \sum_{i,j \in
I_{q(t)}}{   \Big[ - \big( d(n) + d(m) \big) \Delta
 u^\eps_{n,m} \big( x_i - x_j \big) + \micp \big( n , m \big)
    V_{\epsilon} \big( x_i - x_j \big) u_{n,m}^{\epsilon} \big( x_i - x_j) \Big]
} \nonumber \\
    & & \qquad \qquad \qquad \qquad \qquad  \qquad \qquad \qquad \qquad  \qquad \qquad \qquad  \times \, \, J_n(x_i,t) \,  
 \mathbf{1}_{m_i
    = n , m_j = m} \, .  \label{geerevd} 
\end{eqnarray}

By~(\ref{eueps}),  $V_\eps(\cdot) = \eps^{-2} V(\cdot/\eps)$ and~(\ref{euscaling}), the quantity in the square brackets above equals 
$$
 \eps^{-d} \alpha(n,m) \Big( U\big( \tfrac{x_i - x_j}{\eps} \big) +  V \big( \tfrac{x_i - x_j}{\eps} \big)  u_{n,m} \big( \tfrac{x_i - x_j}{\eps} \big) \Big) \, . 
$$

Our aim is that~(\ref{geerevd}) will equal $Q_0(t)$; we see that this demand is equivalent to 
the identity $U + V u_{n,m} = V$. That is, the function $U$ must be specified by $U = V \big( 1 - u_{n,m} \big)$.
Assuming for now that such a choice may be made, consider the term which is analogous to $\int_0^T Q_z(t) \dd t$ in~(\ref{efirstconc})
when the earlier analysis is replayed with $X_z$ in place of $S_z$. This new term equals
\begin{equation}\label{ppo}
     - \,  \epsilon^{2(d-2)} \int_0^T \sum_{i,j \in
I_{q(t)}} \big( d(n) + d(m) \big) \Delta u^\eps_{n,m}
\big( x_i - x_j + z \big)
     J_n(x_i,t)   \mathbf{1}_{
     m_i(t)= n ,m_j(t)= m }  \, \dd t \, .
\end{equation}

Recalling the definition of $\overline{Q}$ from the statement of Proposition~\ref{propqqbar},
and noting that  
$$
- \big( d(n) + d(m) \big) \Delta u^\eps_{n,m}(z) = \alpha(n,m) V\big(z/\eps \big) \big( 1 - u^\eps_{n,m}(z) \big)
$$ 
for $z \in \R^d$, we see that~(\ref{ppo}) is precisely  $\int_0^T \overline{Q}_z(t) \dd t$.

That is, setting $U$ as described above, our reprisal of the method yields
\begin{equation}\label{efirstconcnew}
   \int_0^T Q_0(t) \, \dd t =   \int_0^T \overline{Q}_z(t) \, \dd t \, \, + \, {\rm Err}_{n,m}(\eps,z) \, ,
\end{equation} 
in place of~(\ref{efirstconc}), which is precisely the form of the statement asserted by Proposition~\ref{propqqbar}.

To turn these ideas into a proof of Proposition~\ref{propqqbar}, note first that making use of our desired choice of $U$ entails that we argue that  the PDE 
$$
 - \, \Delta u_{n,m}  =  \ald  V \big( 1 - u_{n,m} \big) 
$$
has a unique solution $u_{n,m}:\R^d \to [0,1)$ satisfying $u_{n,m}(z) \to 0$ as $z \to \infty$. This we have already taken care of: see Lemma~\ref{l.g}.  

Our more substantial remaining task is the following.
Defining the functional $X_z$ with this choice of $u_{n,m}$, we must argue that the dominant terms in the  identity~(\ref{hac}) (with $X_z$ in place of $S_z$) are indeed $- \int_0^T \big( \overline{Q}_z(t) - Q_0(t) \big) \, \dd t$; more precisely, we must show that both the left-hand side of (\ref{hac}), and the difference of its right-hand side with  $- \int_0^T \big( \overline{Q}_z(t) - Q_0(t) \big) \, \dd t$, satisfy the demand made of the error ${\rm Err}_{n,m}(\eps,z)$ in the statement of Proposition~\ref{propqqbar}.
\section{Proof of Proposition 5.1}\label{s.esterr}
We now present the proof of the pointwise Stosszahlansatz, or rather, reduce it to certain key estimates (and we do so making some simplifications which in no way diminish the essentials of the argument). These estimates are gathered at the end of the Section~\ref{s.action}, in Proposition~\ref{p.errorbounds}. 
To make this reduction, the job at hand is to carry out the task mentioned in the preceding paragraph. 
In order to analyse the various error terms, we begin by providing formulas for them.
\subsection{The action of the free motion and collision operators on the functional}\label{s.action}
Recall that, for each $z \in \R^d$, under $\PP_N$, we are defining $X_z:[0,\infty) \to \R$ by means of 
\begin{equation}\label{exzonew}
X_z(t) (q) = \epsilon^{2(d-2)} \sum_{i,j \in I_{q(t)}}
 u_{n,m}^{\epsilon} (x_i- x_j + z) J_n (x_i,t)
\mathbf{1}_{ m_i = n , m_j = m} \, ,
\end{equation}
where, for each pair $(n,m) \in
\mathbb{N}^2$, we define $u_{n,m}^{\epsilon}:
\mathbb{R}^d \to (0,\infty)$ so that, for $z \in \R^d$, 
\begin{equation}\label{e.ue}
 - \, \Delta u^\eps_{n,m}(z)  =  \eps^{-d} \ald  V(z/\eps) \big( 1 - u_{n,m}(z/\eps) \big) 
\end{equation}
subject to $\lim_{x \to \infty} u_{n,m}^{\epsilon}(x) = 0$; moreover, $u_{n,m}^\epsilon$ and $u_{n,m}$ enjoy the scaling relationship~(\ref{euscaling}).

Our aim is to analyse the high-$N$ behaviour of the terms in the $\PP_N$-almost sure identity 
\begin{eqnarray}\label{xidentity}
    \big( X_z - X_0 \big) \big( T \big) & = &
    \big( X_z - X_0 \big) \big( 0 \big) +
\int_{0}^{T}{ \left(\tfrac{\partial}{\partial t}+\freem 
\right)(X_z - X_0) (t) \, \dd t} \\
    & & \qquad + \,
\int_{0}^{T}{ \collop (X_z - X_0) (t) \, \dd t} \, + \, M_T \, , \nonumber
\end{eqnarray}
where the process $\big\{ M_T : T \geq 0 \big\}$ is a $\PP_N$-martingale.

To simplify our presentation, we will consider only the case that the test function $J:\R^d \times \N \times [0,\infty) \to [0,\infty)$ takes the form $J(x,m',t) = J_n \mathbf{1}_{m'=n}$, where $J_n \in (0,\infty)$ is a constant; in the general case, $J_n:\R^d \times [0,\infty) \to [0,\infty)$ is a map on space-time. (Technically, this simplication is not a special case. After all, the test functions $J_n$ in~(\ref{syspdeweak}) are compactly supported. However, the transition back to valid choices of $J_n$ from constant functions is a minor technical point, because it entails only the analysis of some extra terms which are better behaved than terms we will anyway have to treat.) We also write $\overline{J}: \R^d \times \N \times [0,\infty) \to [0,\infty)$ for the function $\overline{J}(x,m',t) = \mathbf{1}_{m' = m}$. In this way, we may write
\begin{equation}\label{exzonewnew}
X_z(t) (q) = \epsilon^{2(d-2)} \sum_{i,j \in I_{q(t)}}
 u_{n,m}^{\epsilon} (x_i- x_j + z) J (x_i,m_i,t) \overline{J}(x_j,m_j,t) \, .
\end{equation}

We are about to label terms arising from the action of the free motion operator on the functional~$X_z - X_0$. For this purpose, an extra piece of notation is useful. Recalling that we write $V_\e(\cdot) = \e^{-2} V(\cdot/\e)$, we now set $V^\e(\cdot)$ equal to~$\e^{-d}V(\cdot/\e)$. 
Recalling~(\ref{e.ue}) as well  the identity~(\ref{euscaling}), we see then that we may write
\begin{equation}\label{e.freem}
\left(\tfrac{\partial}{\partial t}+\freem 
\right)(X_z - X_0) (t)  = H_1 + H_2 + H_3 \, , 
\end{equation}
where 
\begin{eqnarray*}
 H_1 & = & - \, \eps^{2(d-2)} \sum_{i,j \in I_q}   \alpha(m_i,m_j) \Big[ V^\eps(x_i - x_j + z) - V^\eps(x_i - x_j) \Big] 
  J(x_i,m_i,t) \overline{J}(x_j,m_j,t) \, , \\
 H_2 & = &  - \, \eps^{2(d-2)} \sum_{i,j \in I_q}   \alpha(m_i,m_j) V_\eps(x_i - x_j) u_{m_i,m_j}^\eps(x_i - x_j)  
  J(x_i,m_i,t) \overline{J}(x_j,m_j,t) \, ,
  \end{eqnarray*}
  and
$$
 H_3  =    \, \eps^{2(d-2)} \sum_{i,j \in I_q}   \alpha(m_i,m_j) V_\eps(x_i - x_j + z) u_{m_i,m_j}^\eps(x_i - x_j + z)  
  J(x_i,m_i,t) \overline{J}(x_j,m_j,t) \, .
$$

The terms arising from the collision operator may be labelled
$$
\collop (X_z - X_0) (t) = G_z(1) + G_z(2) - G_0(1) - G_0(2) \, ,
$$
where $G_z (1)$ equals
\begin{eqnarray}
& & \tfrac{1}{2} \sum_{k,l \in
I_q}{\micp(m_k,m_l) V_{\epsilon}( x_k - x_l)} \epsilon^{2(d-2)} \sum_{i
\in I_q}{} \\
    & & \quad \bigg\{
    \frac{m_k}{m_k + m_l} \Big[  \ue_{m_k,m_i} ( x_k - x_i
    + z) J( x_k , m_k + m_l,t)
    \overline{J}(x_i,m_i,t) \nonumber \\
    & & \qquad \qquad \qquad  + \, \ue_{m_i,m_k} ( x_i - x_k
    + z) J( x_i , m_i ,t) \overline{J}(x_k,m_k + m_l,t)  \Big]
     \nonumber \\
    & & \quad \, \, + \, \frac{m_l}{m_k + m_l}  \Big[ \ue_{m_l,m_i} ( x_l - x_i
    + z) J( x_l , m_k + m_l,t)
    \overline{J}(x_i,m_i,t) \nonumber \\
    & &  \qquad \qquad \qquad  + \, \ue_{m_i,m_l} ( x_i - x_l
    + z) J( x_i , m_i ,t) \overline{J}(x_l,m_k + m_l,t) \Big]
     \nonumber \\
    & & \qquad  - \ \
    \Big[  \ue_{m_k,m_i} (x_k - x_i + z) J( x_k , m_k,t)
    \overline{J}(x_i,m_i,t) \nonumber \\
    & &  \qquad \qquad \qquad \qquad + \,  \ue_{m_i,m_k} (x_i - x_k + z)
    J( x_i , m_i ,t) \overline{J}(x_k,m_k,t) \Big]
     \nonumber \\
    & & \qquad - \ \
     \Big[  \ue_{m_l,m_i} (x_l - x_i + z) J( x_l , m_l,t)
    \overline{J}(x_i,m_i,t) \nonumber \\
    & & \qquad \qquad \qquad \qquad + \,  \ue_{m_i,m_l} (x_i - x_l + z) J( x_i ,
    m_i ,t) \overline{J}(x_l,m_l,t) \Big] \bigg\}
     \nonumber ,
\end{eqnarray}
and where
\begin{equation}\label{eqnzeta}
G_z (2)  =  - \epsilon^{2(d-2)} \sum_{k,l \in
I_q}{ \micp(m_k,m_l)  V_{\epsilon}(x_k - x_l)}
\ue_{m_k,m_l} (x_k -
x_l + z ) J(x_k,m_k,t)  \overline{J}(x_l,m_l,t).
\end{equation}

In the triple sum over distinct particle indices $(k,l,i)$ appearing in $G_z(1)$, the particles indexed by $k$ and $l$ are interacting at rate $\micp(m_k,m_l) V_{\epsilon}( x_k - x_l)$; when this pair collides, there is an instantaneous change in the value of those terms in $X_z(q)$ that include the location $x_i$ of a given third particle indexed by $i$ and not involved in the collision. There are two gain terms, associated to the appearance of a new particle at one or other of $x_k$ and $x_l$, and two loss terms, associated to the disappearance of the particles indexed by $k$ and $l$. 

The term $G_z(2)$ is a double sum over distinct particle indices $(k,l)$ that records the instantaneous change caused by collision of such a particle pair in the value of those terms in $X_z(q)$ expressed in terms only of the elements of that pair. That is, the collision occurs at infinitesimal rate  $\micp(m_k,m_l)  V_{\epsilon}(x_k - x_l)$; when collision happens, the particles with indices $k$ and $l$ disappear, so that the term $\eps^{2(d-2)} \ue_{m_k,m_l} (x_k -
x_l + z ) J_n \mathbf{1}_{m_j = m}$ no longer appears in $X_z(q)$. 

Note that
$$
H_2 = G_0(2) \, .
$$

We find then that
\begin{eqnarray}\label{abc}
& & \bigg\vert
\int_{0}^{T}{ H_1 \big(t\big) \, \dd t} +
\int_{0}^{T}{ H_3 \big(t\big) \, \dd t} \, \bigg\vert \\
& \leq &  \big\vert X_z - X_0 \big\vert \big( q(T) \big) + \big\vert
X_z - X_0 \big\vert \big( q(0) \big) \nonumber \\
& & \, + \,
 \int_{0}^{T}{  \big\vert G_z(1) - G_0(1) \big\vert (t) \, \dd t} +
\int_{0}^{T}{  \big\vert G_z(2) \big\vert (t) \, \dd t} \, + \, \big\vert
M(T) \big\vert \, . \nonumber
\end{eqnarray}
We will now state bounds on these error terms which are sufficient for the purpose of proving Proposition~\ref{propqqbar} (and thus Proposition~\ref{propsz} and Theorem~\ref{thmo}). 
\begin{prop}\label{p.errorbounds}
Suppose that the survey assumptions are in force.
There exists a constant $C > 0$ such that
\begin{enumerate}
 \item 
for all $N \in \N$,
$\sum_{m \in \N}  \int_{0}^{T}{ \E_N \big\vert G_z(1) - G_0(1) \big\vert (t) \, \dd t} \leq C T^{3d/2} \vert z \vert \big( \log 1/ \vert z \vert \big)^{3d/2}$; 
 \item
for any $t \geq 0$,
$\sum_{m \in \N}  \E_N \big\vert X_z(t) - X_z(0) \big\vert \leq C \vert z \vert$;
 \item 
for all $N \in \N$,
$\sum_{m \in \N}  \int_{0}^{T}{ \E_N \big\vert G_z(2)  \big\vert (t) \, \dd t} \leq C \big( \tfrac{\eps}{z} \big)^{d-2}$;
 \item and,  for each $t \geq 0$,
$\sum_{m \in \N}  \E_N \big[ M(t)^2 \big]  \leq C  \eps^{d-2}$.
\end{enumerate}
\end{prop}
\subsection{Proving the error bounds}

Two important tools are needed to prove the above bounds. Here, we present these two tools (but do not yet prove the assertions we state about them), and use them to give a proof of Proposition~\ref{p.errorbounds}(1); the three other estimates in this proposition follow in a roughly similar way, and we do not give the proofs of these estimates here. 

\subsubsection{Particle concentration bounds}
The first tool is an assertion that, at any given time, the joint density of any given number of particles is uniformly bounded above. Recall that $\big\{ h_n: n \in \N \big\}$ is the initial density profile of particles under~$\PP_N$, and that $\ell_n$ denotes $\dist h_n \dist_{L^\infty(\R^d)}$.
\begin{prop}\label{p.pl}
Suppose that $d:\N \to (0,\infty)$ is non-increasing.
For $k \in \N$, let $g_k: \R^{dk} \times [0,\infty) \to [0,\infty)$ be such that $g_k(y_1,\ldots,y_k,t)$ is the density at $(y_1,\ldots,y_k) \in (\R^d)^k$ for the ordered presence of the lowest-$k$ indexed particles in $\PP_N$ at time $t$. 
Then 
$$
 \vert\vert g_k \vert\vert_{L^\infty(\R^{d k} \times [0,\infty))} \leq 
  \big( Z^{-1}  L s \big)^k \, ,
$$
where 
 $L = \sum_{n=1}^\infty \ell_n n d(n)^{d/2}$
and $s = \sup_{m \geq 1} m^{-1} d(m)^{-d/2}$. 
\end{prop}
Recall that the quantities $L$ and $s$ are both supposed to be finite under the survey assumptions.
\subsubsection{Uniform control on pairwise collision probabilities}\label{s.killing}

Recall from Section~\ref{s.killedbm} the notion of Brownian motion on $\R^d$ killed at rate $W$, where  $W:\R^d \to [0,\infty)$ is a smooth and compactly supported function; recall from there that $u_W:\R^d \to [0,1]$ is such that $u_W(x)$ is the  
 probability that rate-two Brownian motion killed at rate $W$ and begun at~$x$ is killed at some time, and also that we set $u_W^\eps(\cdot) = \eps^{2 - d} u_W\big( \cdot/\eps \big)$. 

\begin{lemma}\label{l.unifcont}
For each $d \geq 3$, there exist a constant $C_d > 0$ such that, for all continuous  $W: \R^d \to [0,\infty)$ with support in the Euclidean unit ball,  
\begin{itemize}
\item for $x \in \R^d$,
$u_W(x) \leq \dist x \dist^{2-d}$ and
$\dist \nabla u_W (x) \dist \leq C_d \, \dist x
\dist^{1-d}$ \, ;
\item for $x,z \in \R^d$ and $\eps > 0$ such that $\dist x \dist \geq  \max \big\{ 2
\dist z \dist + \epsilon , 2  \epsilon \big\}$,
\begin{equation}\label{lept}
    \Big\vert u^\eps_W \big( x + z \big) - u^\eps_W \big( x
    \big) \Big\vert \leq 2^{3d - 6} \frac{\dist z \dist}{\dist x \dist^{d-1}}
\end{equation}
and
\begin{equation}\label{lepta}
    \Big\vert \nabla u^\eps_W \big( x + z \big) - \nabla
u^\eps_W \big( x
    \big) \Big\vert \leq  4^d \big( 2^{d-1} d + 1 \big) \frac{\dist z \dist}{\dist x \dist^d} \, .
\end{equation}
\end{itemize}
\end{lemma}

\subsubsection{Applying the tools}
For the proof of Proposition~\ref{p.errorbounds}(1), we need one further simple estimate, 
on long-range particle displacement. Note that the survey assumptions assure that the initial data is supported in a shared compact set as demanded by the next lemma.
\begin{lemma}\label{l.bigdisp}
Suppose that each element of the initial data $h_n:\R^d \to [0,\infty)$, $n \in \N$, is supported in a given compact region $B$, and that $\overline{d} : = \sup_{n \in \N} d(n) < \infty$. Then, for some constant $C > 0$ and for all $r > 0$,
$$
 \PP_N \Big( \vert x_1(T) \vert \geq r   \Big) \leq C \exp \left\{ - \frac{r^2}{2 \overline{d} T} \right\} \, .
$$ 
\end{lemma}
\noindent{\bf Proof.} If $R > 0$ is an upper bound on the radius of the region $B$, then note that 
$x_1(T)$ under $\PP_N$ is stochastically dominated by the maximum modulus during $[0,T]$ of a rate-$\overline{d}$ Brownian motion begun at distance $R$ from the origin. From this and the reflection principle the result follows. \qed
\noindent{\bf Proof of Proposition~\ref{p.errorbounds}(1).}
Note that 
$$
 \int_{0}^{T}{ \E_N \big\vert G_z(1) - G_0(1) \big\vert (t) \, \dd t} \leq \sum_{i=1}^8 D_i \, ,
$$
where 
\begin{eqnarray}
D_1 & = & \tfrac{1}{2} \, \mathbb{E}_N \int_0^T{} dt \sum_{k,l \in
I_q}{\micp(m_k,m_l)V_{\epsilon}(x_k - x_l)} 
 \big\vert J(x_k,m_k) \big\vert  \\
& & \quad \epsilon^{2(d-2)} \sum_{i \in I_q}{
    \big\vert\overline{J}(x_i,m_i)\big\vert  \cdot \Big\vert
\ue_{m_k,m_i}(x_k - x_i + z) - \ue_{m_k,m_i}(x_k - x_i )
\Big\vert} \, , \nonumber
\end{eqnarray}
and the later $D_i$ terms differ from $D_1$ only in inessential ways. Note that $D_1$ depends implicitly on $n$ and $m$, in that the test functions $J$ and $\overline{J}$ have been chosen to charge only particles of mass $n$ and $m$. The above sum over particle index triples $(k,l,i)$ has at most $N^3 = 
Z^3 \eps^{3(2-d)}$ summands; recalling that $V_\eps = \eps^{-2} V(\cdot/\eps)$ is supported in the $\eps$-ball about the origin, we see that $\sum_{m \in \N} D_1$ is at most
\begin{equation}\label{e.done}
 \eps^{3(2-d)} \cdot \eps^{-2} \cdot \eps^{2(d-2)} Z^3 C_1 \, \int \int_0^T \mathbf{1}_{\vert x_1 - x_2 \vert \leq \eps} \cdot \sup_{m \in \N} \Big\vert \ue_{n,m} \big( x_3 - x_2 + z \big) - \ue_{n,m} \big( x_3 - x_2  \big) \Big\vert \cdot \mu_t\big(\dd x_1,\dd x_2,\dd x_3 \big) \, ,
\end{equation}
where we set $C_1 = \dist V \dist_\infty \cdot \, \sup_{n,m} \alpha(n,m) \cdot \, \sup_n J_n \cdot \, \sup_m \overline{J}_m$, and 
where the law $\mu_t$ is the joint distribution at time $t$ of the first three indexed particles under $\PP_N$ (so that the outer integral in~(\ref{e.done}) is an $\mathbb{E}_N$-expectation over this triple of locations).
By Lemma~\ref{l.bigdisp}, and the uniform bound $\ue \leq 1$, the contribution to the right-hand side made by the integral over 
$(x_1,x_2,x_3) \in \R^{3d} \setminus [-R, R]^{3d}$ is at most $C \exp \big\{ - \tfrac{R^2}{2 \overline{d} T} \big\}$, which equals $\vert z \vert$ if we choose $R^2 = 2 \overline{d} T \log \big( C/\vert z \vert \big)$. We bound the integrand on  $[-R, R]^3$
by noting that the support of $V_\eps$ has volume $\eps^d$ and using Lemma~\ref{l.unifcont} to bound the functions $\big\vert \ue_{n,m}(x+z) - \ue_{n,m}(x)  \big\vert$ simultaneously. Thus,
\begin{eqnarray*}
\sum_{m \in \N} D_1 & \leq &  Z^3 C_1 \big( Z^{-1} L s \big)^3 \int_{[- R,R]^{3d}}  
 \sup_{m \in \N} \big\vert \ue_{n,m} \big( x_3 - x_2 + z \big) - \ue_{n,m} \big( x_3 - x_2  \big) \big\vert \, \dd x \, + \, \vert z \vert \\
& \leq &  \big(  L s \big)^3  \dist V \dist_\infty \cdot \, \sup_{n,m} \alpha(n,m) \cdot \, \sup_n J_n \cdot \, \sup_m \overline{J}_m \cdot
 \big(  2 \overline{d} T \log \big( C/\vert z \vert \big) \big)^{3d/2} \, \cdot C \vert z \vert \, \, + \, \vert z \vert \, ,
\end{eqnarray*}   
where Proposition~\ref{p.pl} was applied in the first inequality in order to replace the measure $\mu_t\big(\dd x_1,\dd x_2,\dd x_3 \big)$ by Lebesgue measure ${\rm d} x$ on $(\R^d)^3$. (Recall that the survey assumptions imply the finiteness of $L$ and $s$.)  Note also that   
 it is (\ref{lept}) that determines the value of the constant $C > 0$.
This completes the proof of~Proposition~\ref{p.errorbounds}(1). \qed

\subsection{Uniform control on pairwise collision probabilities: proofs}

\noindent{\bf Proof of Lemma~\ref{l.unifcont}.} By the uniqueness element of Lemma~\ref{l.s}, we have that, for $x \in \R^d$,
\begin{equation}\label{e.uw}
  u_W(x) =  c_0  \int_{\R^d}  W(y) \big( 1 - u_W(y) \big)  \dist x - y \dist^{2 - d} \, \dd y \, ,
\end{equation}
where we recall from the proof of Lemma~\ref{l.s} that, for each $d \geq 3$, $c_0^{-1} = d(d-2)\omega_d$, with $\omega_d$ equal to the volume of the Euclidean unit ball in $d$ dimensions.

Write $u_\infty:\R^d \to [0,1]$ so that for $x \in \R^d$, $u_\infty(x)$ is the probability that Brownian motion begun at $x \in \R^d$ visits the Euclidean unit ball; our notation is used because formally this object coincides with $u_W$ for $W = \infty \mathbf{1}_{\dist x \dist \leq 1}$. From the interpretation of $u_W$ and $u_\infty$ as killing probabilities, it is evident that $u_W(x) \leq u_\infty(x)$ for all $x \in \R^d$ and for any continuous $W$ supported in the unit ball. However, for $x \in \R^d$,
\begin{equation}\label{e.uinfty}
 u_\infty(x) =   \min \big\{  1, \dist x  \dist^{2-d} \big\} \, ;
\end{equation}
thus, $u_W(x) \leq \dist x \dist^{2 -d}$ for all $x \in \R^d$ and $W$ as above, as Lemma~\ref{l.unifcont} firstly asserts.

By using this monotonicity to compare the formulas~(\ref{e.uw}) and~(\ref{e.uinfty}) along a sequence of $x \in \R^d$ for which $x \to \infty$, we obtain that, for all such potentials $W$,
\begin{equation}\label{e.wclaim}
   \int_{\R^d} W(y) \big( 1 - u_W(y) \big) \, \dd y \leq c_0^{-1} \, . 
\end{equation}

Note that
\begin{equation*}
    \nabla u_W(x) = c_0 (2-d) \int_{\R^d}  W(y)
     \frac{x-y}{\dist x - y \dist^d} \big(  1 - u_W ( y )  \big)  \, \dd y \, .
\end{equation*}

From (\ref{e.wclaim}), we see that $\dist \nabla u_W(x) \dist \leq  (d-2) \int_{\dist y \dist \leq 1} \dist x - y \dist^{1 - d} \dd y \leq C_d \dist x \dist^{1 - d}$ whenever $\dist x \dist \geq 2$, as we also asserted. On the other hand, that $\dist \nabla u_W(x) \dist \leq C_d$ when $\dist x \dist \leq 2$ is straightforward. We have obtained Lemma~\ref{l.unifcont}'s second assertion.

As we turn to derive~(\ref{lept}), 
we mention that, for the rest of the proof, we will denote the Euclidean norm on $\R^d$ by $\vert \cdot \vert$ rather than by $\dist \cdot \dist$.
Note that, for $x,z \in \R^d$,
\begin{eqnarray*}
 & & \big\vert u^\eps_W(x+z) - u^\eps_W(x) \big\vert = \eps^{2-d} \big\vert u_W\big(\tfrac{x+z}{\eps}\big) -  u_W\big(\tfrac{x}{\eps}\big) \big\vert \\
 & \leq & c_0    \int_{\R^d}   W \big( y \big) \big( 1 - u_W(y) \big)  \Big\vert \big\vert x + z - \eps y \big\vert^{2 - d} -  \big\vert x  - \eps y \big\vert^{2 - d} \Big\vert \, \dd y \\
 & \leq  &  \sup_{\vert y \vert \leq \eps} \, \Big\vert \big\vert x + z - y \big\vert^{2 - d} -  \big\vert x  - y \big\vert^{2 - d} \Big\vert \\
  & = &   \sup_{\vert y \vert \leq \eps} \,  \frac{\left| \vert x - y
    \vert^{d-2} -
    \vert x + z - y \vert^{d-2} \right|}{\vert x + z  - y \vert^{d-2}
    \vert x -y \vert^{d-2}} \, ,  
\end{eqnarray*}
the second inequality by~(\ref{e.wclaim}).

Note that
\begin{equation}\label{chon}
    \big\vert x + z - y \big\vert^{d-2} - \big\vert x - y \big\vert^{d-2}
    \leq  \Big( \big\vert x - y \big\vert  + \vert z \vert \Big)^{d-2} -
    \vert x - y \vert^{d-2}
\end{equation}
and that
\begin{equation}\label{chtwo}
    \big\vert x  - y \big\vert^{d-2} - \big\vert x + z - y \big\vert^{d-2}
    \leq  \vert x - y \vert^{d-2} -  \Big( \big\vert x - y \big\vert   - \vert z
    \vert \Big)^{d-2} \, .
\end{equation}
The right-hand sides of (\ref{chon}) and (\ref{chtwo}) each take the form $\alpha^{d-3} \vert z \vert$, for some $\alpha \in \big[ \vert x-y
\vert - \vert z \vert, \vert x- y \vert + \vert z \vert \big]$. Note that
if $| y | \leq \eps$, then
$$
\vert x - y \vert - \vert z \vert \geq \vert x \vert -  \epsilon -
    \vert z \vert\ge 0,
    $$
since $|x| \geq | z | + \eps$. As a result, we have that
$\vert x - y \vert + \vert z \vert \leq 2 \vert x-y \vert$,
so that
\begin{displaymath}
    \Big\vert \big\vert x - y \big\vert^{d-2} - \big\vert x + z - y
    \big\vert^{d-2} \Big\vert \leq 2^{d - 3} \vert z \vert \, \vert x-y \vert^{d-3}.
\end{displaymath}
Hence,
\begin{displaymath}
    \Big\vert u_W^\eps \big( x + z \big) - u_W^\eps \big( x
    \big) \Big\vert  \leq 
     2^{d-3} \vert z \vert \sup_{\vert y \vert \leq \eps}  \vert x -y  \vert^{-1} \vert x + z - y
    \vert^{2-d} \, .
\end{displaymath}
   From $\vert x \vert \geq \max \big\{ 2 \vert z \vert +
 \epsilon , 2  \epsilon \big\}$ and $| y | \leq \eps$ follows  $\big\vert x +
z - y \big\vert \geq \vert x -y \vert/2$ and $\vert x - y \vert \geq
\vert x \vert/2$; thus, the above supremum is at most $2^{2d - 3} | x |^{1 - d}$. We obtain (\ref{lept}).

In seeking to prove (\ref{lepta}), note that
\begin{displaymath}
    \frac{x + z - y}{\big\vert x + z - y \big\vert^d} - \frac{x -
    y}{\big\vert x - y \big\vert^d} = \frac{ \big(x + z - y \big)
    \big\vert x - y \big\vert^d - \big( x - y \big) \big\vert x + z - y
    \big\vert^d}{ \big\vert x + z - y \big\vert^d \big\vert x - y \big\vert^d }.
\end{displaymath}
Note that, for any $a \in \mathbb{R}^d$,
\begin{equation}\label{cule}
    \Big\vert \big( a + z \big) \vert a \vert^d - \vert a + z \vert^d a
    \Big\vert \leq  \vert a \vert \Big\vert \vert a + z \vert^d -
    \vert a \vert^d \Big\vert + \vert z \vert \vert a \vert^d
     \leq  \big( 2^{d-1} d + 1 \big) \vert z \vert \vert a \vert^d \, ,
\end{equation}
as long as $|z|\le|a|$. Given that
\begin{displaymath}
    \nabla u_W^\eps(x) = - (d-2)c_0   \int_{\R^d}  W(y) \big( 1 - u_W(y) \big)
      \frac{x - y}{\big\vert x - \eps y \big\vert^{d}} \, \dd y \, ,
\end{displaymath}
we may apply (\ref{cule}) with the choice $a = x - \eps y$ and then use~(\ref{e.wclaim}) to obtain
\begin{eqnarray*}
    & & \Big\vert  \nabla u_W^\eps ( x + z ) - \nabla u_W^\eps (x)
     \Big\vert \\ 
     & \leq & \big( 2^{d-1} d + 1 \big)  \vert z \vert  c_0 \int_{\R^d}
    W(y) \big( 1 - u_W(y) \big)  \big\vert x + z - \eps y \big\vert^{-d} \, \dd  y \, \\
     & \leq &  \big( 2^{d-1} d + 1 \big)  \vert z \vert \sup_{\vert y \vert \leq \eps}  \big\vert x + z -  y \big\vert^{-d} \, .
\end{eqnarray*}

From $\vert x \vert \geq \max \big\{ 2 \vert z \vert +
 \epsilon , 2  \epsilon \big\}$ and $| y | \leq \eps$, we see that  $\big\vert x +
z - y \big\vert \geq \vert x -y \vert/2$ and $\vert x - y \vert \geq
\vert x \vert/2$. We conclude that
\begin{displaymath}
     \Big\vert  \nabla u_W^\eps ( x + z ) - \nabla u_W^\eps (x)
     \Big\vert \leq   4^d \big( 2^{d-1} d + 1 \big)  \frac{\vert z \vert}{\vert x \vert^d} ,
\end{displaymath}
as required. \qed
\section{Particle concentration bounds and uniform integrability}\label{s.partcon}

The principal aim of this section is to prove the particle concentration upper bound, Proposition~\ref{p.pl}, and the closely related uniform integrability assertions, Propositions~\ref{prop.firststep} and~\ref{prop.secondstep}; though the section also includes, at its end, an analytic derivation of mass conservation under certain assumptions.

Proposition~\ref{p.pl} asserts that, if $d:\N \to (0,\infty)$ decreases, but not too rapidly, then supremum norm bounds enjoyed initially by the particle profile propagate to all later times up to factors determined by the diffusion rates $d$. Such particle concentration results play an essential role in our derivation of the Smoluchowski PDE; unlike in our preceding work, here we present a proof of these results using probabilistic techniques. 
The derivation will occupy several pages and invokes moderately restrictive hypotheses on $d(\cdot)$. 
We make use of this approach because of the attractive probabilistic perspective that it offers one of the more technical aspects of our kinetic limit derivation of the Smoluchowski PDE; and because the uniform integrability Proposition~\ref{prop.secondstep} is an immediate corollary.  

It is much quicker to describe the supremum norm propagation effect in terms of solutions to the PDE. In order to illustrate the effect succinctly to begin with, and perhaps also for the benefit of analytically minded readers who may wish to skip some details in the upcoming proof of Proposition~\ref{p.pl}, we first present the statement and proof of \cite[Lemma 4.1]{momentbounds}. Such a reader may also wish to consult~\cite{YRH}, where an analogous particle distribution result, Theorem 3.1, is proved by means not unlike, but more analytic than, our approach to establishing Proposition~\ref{p.pl}.


\subsection{An analytic bound on particle concentration}

Our analytic lemma concerns a weak solution $\big\{ f_n: n \in \N \big\}$, one that solves the system~(\ref{syspdeweak}).
\begin{lemma}\label{lem3.1} 
Assume $d(\cdot)$ is non-increasing. Then, for all $x \in \R^d$ and $t \geq 0$,
 \begin{equation}\label{eqn3.1} 
 \sum_{n=1}^{\i} n d(n)^{d/2} f_n(x,t) \leq d(1)^{d/2} u(x,t) \, ,
\end{equation}
where $u$ is the unique solution to $u_t = d(1) \D u$ subject to the
initial condition $u(x,0) = \sum_{n=1}^{\i} n f_n(x,0)$.
\end{lemma}
\noindent {\bf Proof.} For $D > 0$, let $\big\{ S_t^D : t \geq 0 \big\}$ denote the diffusion rate $D$ heat semigroup. That is, for any continuous function $f:\R^d \to \R$, $S_t^D f : \R^d \to \R$ is given by
 $$
 S_t^D f (x) =  \int_{\R^d} f(x-y) \cdot \big( 2\pi Dt \big)^{-d/2} \exp \big\{ - \tfrac{y^2}{2Dt} \big\} \, \dd y \, .
 $$
 
 The heat semigroup satisfies the property that, if $D_1 \ge
D_2$ and $g \ge 0$, then
\begin{equation}\label{e.heatbound} 
D_1^{d/2} S_t^{D_1} g \ge  D_2^{d/2}
S_t^{D_2}g \, ;
\end{equation}
this is a consequence of an elementary bound on the normal density, which in fact we will shortly state as~(\ref{e.normalbound}).

 Using the shorthand $Q_n(x,s) = Q_n(f)(x,s)$, $s \in [0,\infty)$, note that, for the collision operator in~(\ref{syspde}), Duhamel's principle implies the basic relation that, for each $n \in \N$, and for all $(x,t) \in \R^d \times [0,\infty)$,
\begin{equation}\label{e.duhamel}
f_n(x,t) = S_t^{d(n)} h_n (x) + \int_0^t
S_{t-s}^{d(n)} Q_n(x,s) \, \dd s \, ,
\end{equation}
where recall that $h_n = f_n(\cdot,0)$ denotes the initial condition.
 
We will argue that, for each $\ell \in \N$, and for all $(x,t) \in \R^d \times [0,\infty)$,
\begin{equation}
\label{eqn3.2} \sum_1^{\ell} n d(n)^{d/2}f_n(x,t) \le d(1)^{d/2}
S_t^{d(1)} \left( \sum_1^{\ell} n h_n \right) (x) + d(\ell)^{d/2}
\int_0^t S_{t-s}^{d(\ell)} \left( \sum_1^{\ell} n Q_n(x,s)\right)\, \dd s \, ;
\end{equation}
first, let us show that this claim proves the lemma.

Consider the expression  $\sum_1^{\ell} n Q_n(x,t)$ for any $(x,t) \in \R^d \times [0,\infty)$. Coagulations at $(x,t)$ between pairs of particles whose combined mass is at most $\ell$, or each of whose masses are at least $\ell + 1$, do not contribute to the expression, while the remaining coagulations contribute negatively. Thus,  $\sum_1^{\ell} n Q_n(x,t) \leq 0$. We find then from~(\ref{eqn3.2}) that
\begin{equation}
\label{eqn3.3} \sum_1^{\ell} n d(n)^{d/2}f_n \le d(1)^{d/2}
S_t^{d(1)} \left( \sum_1^{\ell} n f_n^0\right) \, .
\end{equation}
In this way, we see that, to prove the lemma, it suffices to derive~(\ref{eqn3.2}).

We will establish this bound by induction on $\ell \in \N$.  When $\ell =1$, the bound holds as an equality due to~(\ref{e.duhamel}).

 Supposing that (\ref{eqn3.2}) is valid at index $\ell$, we now derive it at index $\ell + 1$.  
By (\ref{e.heatbound}) and (\ref{eqn3.2}), we learn that,  for all $(x,t) \in \R^d \times [0,\infty)$,
\begin{equation}\label{eqn3.5} 
\sum_1^{\ell} n d(n)^{d/2}f_n(x,t) \le d(1)^{d/2}
S_t^{d(1)} \left( \sum_1^{\ell} n h_n \right)(x) + d(\ell+1)^{d/2}
\int_0^t S_{t-s}^{d(\ell+1)} \left( \sum_1^{\ell} n Q_n(x,s)\right) \, \dd s
\end{equation}
because $d(\ell) \ge d(\ell+1)$ and $\sum_1^{\ell} n Q_n(x,t) \le 0$.

Applying (\ref{e.heatbound}) to (\ref{e.duhamel}) with $n = \ell + 1$ yields
\begin{equation}\label{eqn3.6} 
 f_{\ell+1} (x,t) \le \left( \frac
{d(1)}{d(\ell+1)}\right)^{d/2} S_t^{d(1)} h_{\ell+1}(x) + \int_0^t
S_{t-s}^{d(\ell+1)} Q_{\ell+1}(x,s)ds.
\end{equation}
We multiply both sides of $(\ref{eqn3.6})$ by $(\ell+1)d(\ell+1)^{d/2}$ and
add the result to $(\ref{eqn3.5})$.  The outcome is
\[
\sum_1^{\ell+1} n d(n)^{d/2} f_n(x,t) \le d(1)^{d/2} S_t^{d(1)} \left(
\sum_1^{\ell+1} n h_n \right)(x) + d(\ell+1)^{d/2}
 \int_0^t S_{t-s}^{d(\ell+1)} \left( \sum_1^{\ell+1} n
Q_n(x,s)\right)ds.
\]
This completes the proof.\qed

\subsection{Proof of Proposition 6.2}

We begin the proof by reformulating the proposition as Proposition~\ref{proppartplace} and proving this. 
%
%
Recall that $\mu$ denotes Lebesgue measure on $\R^d$.

\begin{prop}\label{proppartplace}
Suppose that $d:\N \to (0,\infty)$ is non-increasing.
For any $N,k \in \N$  with $N \geq k$ and $T \geq 0$, 
\begin{equation*}
 \PP_N \bigg(  \bigcap_{i=1}^k \big\{   x_i(T) \in    A_i   \big\}      \bigg)  \, \leq \,
   K^k \prod_{i=1}^k \mu (A_i) \, ,
\end{equation*}
where 
 $\big\{ A_i: 1 \leq i \leq k \big\}$ is any collection of open sets in $\R^d$. The positive constant $K$ equals $Z^{-1} \sum_{n \geq 1}   n   
d(n)^{d/2} \ell_n \sup_{m \geq n} m^{-1}  d(m)^{-d/2}$, where recall that $\ell_n$ denotes $\vert\vert h_n  \vert\vert_{L^\infty(\R^d)}$ 
for $n \in \N$.
\end{prop}

\noindent{\bf Proof of Proposition~\ref{p.pl}.}
Using Proposition~\ref{proppartplace}, note that 
$g_k\big( (y_1,\cdots,y_k) , t \big) \leq  K^k$
for any $ (y_1,\cdots,y_k)  \in (\R^d)^k$ and $t \geq 0$.
Recalling from Proposition~\ref{p.pl}'s statement that we set $L = \sum_{n=1}^\infty \ell_n n d(n)^{d/2}$ and $s =  \sup_{m \geq 1} m^{-1} d(m)^{-d/2}$, we see that 
$K \leq Z^{-1} L s$. \qed

The proof of Proposition \ref{proppartplace}  relies on the specific details of pairwise collision seen in the dynamics under $\PP_N$ which we defined in Section~\ref{sec.micromodels}; the reader may wish to recall these now by consulting the paragraphs in Section~\ref{sec.micromodels} under the heading ``the precise mechanism of collision''. 

\subsubsection{The method of proof: the tracer particle}

We will first prove Proposition \ref{proppartplace} with $k=1$. 
To do so, it is convenient to interpret the lowest indexed particle, with index one, as the {\em tracer particle}. (In fact, we used this term already, in treating the translation invariant model on a torus seen in Section~\ref{s.torus}.) 
By symmetry of the initial particle placements under~$\PP_N$, the tracer particle is indistinguishable from a particle selected uniformly at random at time zero. 

\subsubsection{A stronger inductive hypothesis}

We will prove Proposition  \ref{proppartplace} with $k=1$ by formulating a stronger inductive hypothesis where the parameter for the induction is the initial total particle number $N$. As a shorthand, we write
\begin{equation}\label{e.normaldensity}
 \nu_{x,s}(\dd y) =  \big( 2 \pi s \big)^{-d/2} \exp \big\{ - \tfrac{\vert\vert x - y \vert\vert^2}{2s} \big\} \, \dd y
\end{equation}
for the law of a normal random variable of mean $x \in \R^d$ and variance $s \geq 0$.
\begin{lemma}\label{l.indhyp}
 Let $N \geq 1$. For given $x \in \R^d$, $n_0 \in \N$ and $\chi \in \big( \R^d \times \N \big)^{N-1}$, 
let  $\PP_N^{x,n_0,\chi}$ denote the law $\PP_N$ conditionally on 
$\xtp(0) = x$, $\mtp(0) = n_0$ and on the other $N-1$ particles at time zero having locations and masses given by $\chi$.
Then, for all such $(x,n_0,\chi)$, and for any $T \geq 0$ and $A \subseteq \R^d$ open,
\begin{equation}\label{e.indhyponenew}
  \PP_N^{x,n_0,\chi} \big( x_1(T) \in A \big) \leq   \sup_{m \geq n_0}  \tfrac{n_0}{m} \big( \tfrac{d(n_0)}{d(m)} \big)^{d/2}  \cdot \nu_{x,2 d(n_0) T}(A) \, .
\end{equation}  
\end{lemma}

Before beginning this lemma's proof, we give an overview of the argument, which is an induction on $N \geq 1$. The case $N=1$ may seem to be a triviality. Collision being impossible for a single particle, $x_1(T)$ is normally distributed with mean $x$ and variance $2 d(n_0) T$, while $m_1(T)$ equals $n_0$ almost surely; from which~(\ref{e.indhyponenew}) follows. One might object however that conditioning on non-collision on the part of the tracer particle will bias the law of its trajectory; formally, we will treat the case $N=1$ as an instance of the generic step of the induction.

In the case of several particles, a key role is played by the following uniform bound on normal densities: for all $x,y \in \R^d$ and $s \geq 0$, and for all $m,m' \in \N$ such that $m' > m$,
\begin{equation}\label{e.normalbound}
   \frac{\dd \nu_{x,2d(m')s}}{\dd \nu_{x,2d(m)s}}(y) \leq  \bigg( \frac{d(m)}{d(m')} \bigg)^{d/2} \, . 
\end{equation}
This bound is a consequence of $d:\N \to (0,\infty)$ being non-increasing.

We make some informal comments about the case $N=2$ in order to illustrate the idea of the proof of Lemma~\ref{l.indhyp}. The tracer particle may make only one collision in this case. Consider the tracer particle dynamics until the first collision time for this particle at some time $t$, which we may assume to be less than $T$ (for the other case is in effect the $N=1$ case). Suppose that the tracer particle has mass $n_0$ just before time $t$ and collides with a mass $n_1$ particle at that time. Compare the subsequent dynamics to an altered one in which the two particles in the model do not interact. In the ordinary dynamics, the tracer particle survives the collision with probability $\tfrac{n_0}{n_0 + n_1}$, and then pursues a Brownian trajectory of diffusion rate $2 d(n_0 + n_1)$; thus, the conditional probability given the tracer particle trajectory until first collision at time $t$ that $\xtp(T) \in A$ equals $\tfrac{n_0}{n_0 + n_1} \nu_{\xtp(t),2d(n_0 + n_1)(T - t)}(A)$. On the other hand, in the altered 
dynamics, 
the 
tracer particle remains of mass $n_0$ at time $t$, and thus has conditional probability $\nu_{\xtp(t),2d(n_0 )(T - t)}(A)$ of achieving $\xtp(T) \in A$. The uniform bound~(\ref{e.normalbound}) on normal densities
implies that the conditional probability of $x_1(T) \in A$ for the ordinary dynamics exceeds that for the altered dynamics by a factor of at most $\tfrac{n_0}{n_0 + n_1} \big(  \tfrac{d(n_0)}{d(n_0 + n_1)} \big)^{d/2}$. Noting that the tracer particle in the altered dynamics is simply a Brownian particle of diffusion rate $2d(n_0)$ for all time, we may average over the tracer particle trajectory until first collision, and, in doing so, we see that $\PP_N\big( x_1(T) \in A \big)$ is at most a $\tfrac{n_0}{n_0 + n_1} \big( \tfrac{d(n_0)}{d(n_0 + n_1)} \big)^{d/2}$-multiple of the probability $\mu_{x,2d(n_0)T}(A)$ that the altered dynamics tracer particle reaches $A$ at time $T$. Taking a supremum in $n_1$ heuristically explains~(\ref{e.indhyponenew}) when 
$N=2$. When $N > 2$, the tracer particle may collide several times, with particles of successive masses $n_1,n_2,\cdots,n_k$, say. In essence, the same line of argument works, with a comparison factor of $\tfrac{\sum_{j=0}^{i-1}n_j}{\sum_{j=0}^i n_j} \Big( \tfrac{d(\sum_{j=0}^{i-1} n_j)}{d(\sum_{j=0}^i n_j)} \Big)^{d/2}$ being associated to the $i\textsuperscript{th}$ collision. The product of these telescoping factors, $\tfrac{m_0}{\sum_{j=0}^k n_j} \Big( \tfrac{d(n_0)}{d(\sum_{j=0}^k n_j)} \Big)^{d/2}$, is then an upper bound on the ratio of the probabilities of $x_1(T) \in A$ in the interacting model~$\PP_N$ and in the model formed from $\PP_N$ by the suppression of all collisions on the part of the tracer particle. Heuristically this also explains the form~(\ref{e.indhyponenew}) when $N \geq 2$. 

\subsubsection{Deriving Proposition~\ref{proppartplace} when $k=1$}

\noindent{\bf Proof of Lemma~\ref{l.indhyp}.} We turn to the rigorous argument establishing the general inductive step. Let $N \geq 1$ be given. Assume then that the statement of Lemma~\ref{l.indhyp} is known for values of the inductive parameter strictly less than $N$. We will analyse $\PP_N^{x,n_0,\chi}\big( x_1(T) \in A \big)$ as an average of the conditional probability of $x_1(T) \in A$ given the tracer particle trajectory until immediately before the first collision. To this end, for $t \geq 0$, we write $\FF_t$ for the $\sigma$-algebra generated by the $\PP_N$-random variables $x_1:[0,t] \to \R^d \cup \{ c \}$ and $m_1:[0,t] \to \N \cup \{ c \}$, so that the information available in $\FF_t$ is the data given by monitoring the tracer particle during $[0,t]$. 
We also write $\FF_t^-$ for the $\sigma$-algebra generated by $\big\{ \FF_s: 0 \leq s < t \big\}$, representing the information concerning the tracer particle's history which is available immediately before time $t$. 
Let $\sigma_1 \geq 0$ denote the time of the tracer particle's first collision. We also set $\sigma_1^T = \sigma_1 \wedge T$.

Note then that we may express
$$
\PP_N^{x,n_0,\chi}\big( x_1(T) \in A \big) = \E \, \PP_N^{x,n_0,\chi}\Big( x_1(T) \in A  \, \Big\vert \, \FF^-_{\sigma_1^T} \Big) \, .  
$$
Note that in the right-hand side, we are writing the probability that $x_1(T) \in A$ as an average over tracer particle histories up to, but not including, the first collision time. 
We consider separately the cases where the first tracer particle occurs {\em before}, or {\em after}, time $T$ (and certainly the first case will be the more demanding). That is,  we work with the identity
\begin{equation}\label{e.beforeafter}
\PP_N^{x,n_0,\chi} \big( x_1(T) \in A  \big)  =   \E_N^{x,n_0,\chi} \bigg[ \, \E \Big[  \big(  \mathbf{1}_{\sigma_1 < T} +  \mathbf{1}_{\sigma_1 \geq T}   \big) \mathbf{1}_{x_1(T) \in A} \, \big\vert \, \FF^-_{\sigma_1^T} \Big]\bigg] \, .
\end{equation}

To begin treating the more difficult, `before', case, we note that,
given any instance of data in $\FF^-_{\sigma_1^T}$ for which $\sigma_1 < T$, it is known that the tracer particle is about to experience a collision at time~$\sigma_1$, even though the mass of the second particle participating in the collision and the collision's outcome -- the survival or perishing of the tracer particle -- remain random events. Consider the instantaneous future of the tracer particle trajectory under the law $\PP_N^{x,n_0,\chi}\big( \cdot  \big\vert \FF^-_{\sigma_1^T} \big)$. We may denote by $n_1$ the mass of the particle with which the tracer particle collides at time $\sigma_1$; note that under the conditional law, $n_1$ is a random variable. Let $\surv$ denote the event that the tracer particle survives this collision, so that the conditional probability of $\surv$ under $\PP_N^{x,n_0,\chi}\big( \cdot  \big\vert \FF^-_{\sigma_1^T}, n_1 \big)$ equals 
$\tfrac{\mass_0}{\mass_0 + n_1}$. If $\surv^c$ occurs, then $x_1$ immediately arrives in the cemetery state~$c$, so that, if $\sigma_1 < T$, there is no possibility that $x_1(T) \in A$ in this event. On the other hand, should $\surv$ occur,
$m_1(\sigma_1) = \mass_0 + n_1$, so that we are able to note that
\begin{eqnarray}
\E \Big[   \mathbf{1}_{\sigma_1 < T}  \mathbf{1}_{x_1(T) \in A} \, \big\vert \, \FF^-_{\sigma_1^T} \Big] & = & 
   \mathbf{1}_{\sigma_1 < T} \cdot
\E \Big[  \mathbf{1}_{x_1(T) \in A} \mathbf{1}_{\surv} \, \big\vert \, \FF^-_{\sigma_1^T} \Big] \label{e.condexp} \\
& = & 
  \mathbf{1}_{\sigma_1 < T} \int  \tfrac{n_0}{n_0 + n_1} 
\PP_n^{x_1(\sigma_1),\mass_0 + n_1,\phi}\big( x_1(T - \sigma_1) \in A \big) \, \dd \mu_{(n,n_1,\phi)} \, . \nonumber 
\end{eqnarray}
Here, the triple $(n,n_1,\phi)$ records the following random data: 
\begin{itemize}
 \item 
$n$, the number of surviving particles immediately after the collision at time $\sigma_1$;
 \item $n_1$, the mass of the particle with which $x_1$ collides at time $\sigma_1$;
 \item and $\phi \in ( \R^d \times \N )^{n-1}$, the vector of locations and masses of the particles other than the tracer particle at this time;
\end{itemize}
while $\mu_{(n,n_1,\phi)}$ denotes the  $\FF^-_{\sigma_1^T}$-measurable random measure that specifies the conditional distribution of $(n,n_1,\phi)$ under the law 
$\PP_N^{x,n_0,\chi}\big( \cdot \big\vert \, \FF^-_{\sigma_1^T} \big)$.

Since a collision occurs at time $\sigma_1$, $n$ is necessarily at most $N-1$, so that the inductive hypothesis may be applied to bound above the term 
$\PP_n^{x_1(\sigma_1),\mass_0 + n_1,\phi}\big( \xtp(T - \sigma_1) \in A \big)$ appearing in the integrand above. We find that, for each $n \geq 1$, 
\begin{eqnarray*}
 & &  \sup_{1 \leq k \leq N-1, \phi \in \big( \R^d \times \N \big)^{k-1}} \PP_k^{\xtp(\sigma_1),n_0 + n,\phi}\big( \xtp(T - \sigma_1) \in A \big) \\
 & \leq &  \sup_{n' \geq n_0 + n} \tfrac{n_0 + n}{n'} \big( \tfrac{d(n_0 + n)}{d(n')} \big)^{d/2}
 \cdot \nu_{ \xtp(\sigma_1),2 d(n_0+n)(T - \sigma_1)}(A) \, .
 \end{eqnarray*}
The right-hand side is at most 
$$
 \sup_{n' \geq n_0 + n} \tfrac{n_0 + n}{n'} \big( \tfrac{d(n_0)}{d(n')} \big)^{d/2} \cdot \nu_{ \xtp(\sigma_1),2 d(n_0)(T - \sigma_1)}(A) 
$$
due to the uniform bound~(\ref{e.normalbound}) in the guise
$\tfrac{\dd \nu_{x,2d(n_0 + n)s}}{\dd \nu_{x,2d(n_0)s}}(y) \leq \big( \tfrac{d(n_0)}{d(n_0 + n)} \big)^{d/2}$ with $s = T - \sigma_1$ (and for any $x,y \in \R^d$).

Equipped with this information, we may return to the second line of~(\ref{e.condexp}), and note the cancellation arising from the product of the terms $\tfrac{n_0}{n_0 + n_1}$ and $\tfrac{n_0 + n_1}{n'}$. Thus, we obtain
an upper bound on the `before' case term in~(\ref{e.beforeafter}):
$$
\E \Big[   \mathbf{1}_{\sigma_1 < T}  \mathbf{1}_{\xtp(T) \in A} \, \big\vert \, \FF^-_{\sigma_1^T} \Big] \leq  
   \mathbf{1}_{\sigma_1 < T} \cdot  \sup_{n' \geq n_0 + 1} \tfrac{n_0}{n'} \big( \tfrac{d(n_0)}{d(n')} \big)^{d/2} \cdot \nu_{ \xtp(\sigma_1),2 d(n_0)(T - \sigma_1)}(A) \, .
$$

For the `after' case, we merely note that  
$$
\E \Big[   \mathbf{1}_{\sigma_1 \geq T}  \mathbf{1}_{\xtp(T) \in A} \, \big\vert \, \FF^-_{\sigma_1^T} \Big] =    \mathbf{1}_{\sigma_1 \geq T}   \mathbf{1}_{\xtp(T) \in A}  \, .
$$

We may now rejoin (\ref{e.beforeafter}) to learn  that
\begin{eqnarray}
 & & \PP_N^{x,n_0,\chi} \big( \xtp(T) \in A  \big) \label{e.postcondexp} \\ 
 & \leq & \E_N^{x,n_0,\chi} \bigg[   \mathbf{1}_{\sigma_1 < T} \cdot  \sup_{n' \geq n_0 + 1} \tfrac{n_0}{n'} \big( \tfrac{d(n_0)}{d(n')} \big)^{d/2} \cdot \nu_{ \xtp(\sigma_1),2 d(n_0)(T - \sigma_1)}(A)  +  \mathbf{1}_{\sigma_1 \geq T}   \mathbf{1}_{\xtp(T) \in A}  \bigg] \, . \nonumber
\end{eqnarray}
To bound above the right-hand side, note first that
$$
 \E_N^{x,n_0,\chi} \Big[   \mathbf{1}_{\sigma_1 < T} \cdot  \nu_{ \xtp(\sigma_1),2 d(n_0)(T - \sigma_1)}(A)  \Big] \leq  \nu_{ x ,2 d(n_0) T}(A)\, . 
$$
To see this, 
recall that, when $\sigma_1 < T$, $\nu_{ \xtp(\sigma_1^T),2 d(n_0)(T - \sigma_1)}(A)$ is the probability that an independent Brownian motion of diffusion rate $2d(n_0)$ beginning from $\xtp(\sigma_1)$ at time $\sigma_1$  visits $A$ at time $T$; in the above left-hand side, the mean is taken over trajectories of a Brownian motion $\xtp:[0,\sigma_1] \to \R^d$ of diffusion rate $2d(n_0)$ with initial condition $\xtp(0) = x$. That is, the value of this left-hand side is given by the probability that a Brownian motion of diffusion rate $2d(n_0)$ which at time zero is at~$x$ visits $A$ at time $T$ (and that an auxiliary stopping time, $\sigma_1$, occurs before $T$).

Similarly, we have that 
$$
 \E_N^{x,n_0,\chi} \Big[    \mathbf{1}_{\sigma_1 \geq T}   \mathbf{1}_{\xtp(T) \in A}  \Big] \leq  \nu_{ x ,2 d(n_0) T}(A)\, .  
$$
Indeed, the subprobability measure $\PP_N^{x,n_0,\chi} \big( \cdot \cap \{ \sigma_1 \geq T \}  \big)$ is stochastically dominated by Brownian motion begun at $x$ of diffusion rate $2 d(n_0)$, since the trajectory of $\xtp$ suffers no collision on $[0,T]$ under this defective law.  

The last inferences when allied with~(\ref{e.postcondexp}) yield
$$
\PP_N^{x,n_0,\chi} \big( \xtp(T) \in A  \big) \leq  \sup_{n' \geq n_0} \tfrac{n_0}{n'} \big( \tfrac{d(n_0)}{d(n')} \big)^{d/2} \cdot  \nu_{ x ,2 d(n_0) T}(A) \, ; 
$$
this is (\ref{e.indhyponenew}) for index $N$, so that the inductive proof of Lemma~\ref{l.indhyp} is complete. \qed

\noindent{\bf Proof of Proposition~\ref{proppartplace} with $k=1$.} 
Using the fact that the tracer particle has mass $n_0 \in \N$ at time zero with probability $Z^{-1} \vert\vert h_{n_0} \vert\vert_{L^1(\R^d)}$, Lemma~\ref{l.indhyp}
yields the deduction that
$$
\PP_N \big( x_1(T) \in A  \big)
\leq  Z^{-1} \sum_{n_0 \geq 1} \vert\vert h_{n_0} \vert\vert_{L^1(\R^d)} \cdot
 n_0 d(n_0)^{d/2} \sup_{n \geq n_0} n^{-1} d(n)^{-d/2}
 \cdot \E \big( \nu_{X_{n_0},2d(n_0)T}(A) \big) \, .
$$
We need to explain some notation on this right-hand side. First we mention that the summand indexed by $n_0$
corresponds to the tracer particle beginning with mass $n_0$. As such, the quantity~$X_{n_0}$ is intepreted as the initial location of the tracer particle given that it has this initial mass; that is,~$X_{n_0}$ is a random variable having the law of $\xtp(0)$ given that $\mtp(0) = n_0$, 
so that~$X_{n_0}$ has density $\tfrac{h_{n_0}(\cdot)}{\vert\vert h_{n_0} \vert\vert_{L^1(\R^d)}}$ on $\R^d$.
The mean in the final term in the summand is taken over $X_{n_0}$.
 
Now, the heat equation decreases the $L^\infty$-norm, so that $\E \big( \nu_{X_{n_0},2d(n_0)T}(A) \big) \leq 
\tfrac{1}{\vert\vert h_{n_0} \vert\vert_{L^1(\R^d)}}  \ell_{n_0} \cdot \mu(A)$, where recall that $\ell_{n_0} = \vert \vert h_{n_0} \vert\vert_{L^\infty(\R^d)}$. 
In this way, we obtain Proposition~\ref{proppartplace} with $k=1$. \qed

\subsubsection{Monitoring several tracer particles at once}

The proof of Proposition \ref{proppartplace} when $k>1$ extends the argument for the case that has already been proved.
In asking about the distribution at time $T$ of collections of particles of size $k$, rather than about single particles, we generalize  the concept of the tracer particle. Where before the tracer particle was the particle of lowest index, we now consider an ordered list of $k$ tracer particles (for given $k \in \N$), these being the particles with indices $1,2,\cdots,k$. 
The symmetry of the initial particle placement means that the joint law of the tracer particles is that of a $k$-sized collection of particles chosen uniformly
at time zero independently of other randomness.

As did the special case when $k=1$, Proposition~\ref{proppartplace} will follow from a stronger assertion which will be established by induction on $N$.

\begin{lemma}\label{l.indhypmult}
Given $\overline{y} = (y_1,\ldots,y_k) \in \R^d$, $\overline{n} = (n_1,\ldots,n_k) \in \N^k$ and  $\chi \in \big( \R^d \times \N \big)^{N-k}$, we denote by $\PP_N^{\overline{y},\overline{n},\chi}$ the law $\PP_N$, conditionally on the tracer particle initial data taking the form
$x_i(0) = y_i$, $m_i(0) = n_i$ for $1 \leq i \leq k$, and on the other $N-k$ particles at time zero having locations and masses given by $\chi$. Then, for all such data, and for any $T \geq 0$ and open $A_i \subseteq \R^d$, $1 \leq i \leq k$, 
\begin{equation}\label{e.indhyponemult}
 \PP_N^{\overline{y},\overline{n},\chi} \bigg( \bigcap_{1 \leq i \leq k} \big\{ x_i(T) \in A_i \big\} \bigg) \,  \leq \,
    \prod_{i=1}^k \sup_{m \geq n_i}  \tfrac{n_i}{m}  \big( \tfrac{d(n_i)}{d(m)} \big)^{d/2} \nu_{y_i,2d(n_i)T}(A_i) \, .
\end{equation}  
\end{lemma}
\noindent{\bf Proof.}
 The proof of the lemma is in essence the same as that of Lemma~\ref{l.indhyp}. It works by induction on $N$, beginning with $N = k$. Under the law $\PP_N$, let $\sigma \geq 0$ denote the time of the earliest collision experienced by any one of the $k$ tracer particles. Set $\sigma' = \min \{ \sigma, T \}$.
If $\sigma' < T$ then either two tracer particles collide at time $\sigma'$, rendering the event that each such particle reaches its target set~$A_i$ at time $T$ impossible, or one of the tracer particles collides with one of the $N-k$ non-tracer particles. The latter case is the non-trivial one, and it may be analysed exactly as in the proof of Lemma~\ref{l.indhyp}: if the tracer particle in question perishes on collision, it is assigned to the state $c$ and there is nothing to prove; otherwise, it survives, assumes some added mass, among at most $N-1$ other particles. The inductive hypothesis is applied to the new post-collision scenario. \qed  


\noindent{\bf Proof of Propositon~\ref{proppartplace}.}
First note that we may find an upper bound on the $\PP_N$-probability that the tracer particles respectively arrive in the sets $A_i$, $1 \leq i  \leq k$, at time $T$ by using Lemma~\ref{l.indhypmult} 
and averaging the provided bound over the initial data of the $k$ tracer particles. We find that 
\begin{equation}\label{e.boundk}
 \PP_N \bigg( \bigcap_{1 \leq i \leq k} \big\{ x_i(T) \in A_i \big\} \bigg)
  \leq   \sum_{ \bar{n} \in \N^k} \int  \prod_{i=1}^k \sup_{m \geq n_i}  \tfrac{n_i}{m}  \big( \tfrac{d(n_i)}{d(m)} \big)^{d/2} \nu_{x_i,2d(n_i)T}(A_i) \, \dd \mu_{\bar{n}}(\bar{y}) \, .
\end{equation}
A few words of explanation concerning notation and reasoning are necessary. First, we are here writing  
$\overline{n}$ for the $k$-vector $(n_1,\cdots,n_k)$ and $\overline{y}$ for $(y_1,\cdots,y_k)$. The right-hand summand is an upper bound for the probability of the intersection of the events that
\begin{itemize}
\item for each $i \in \{1,\cdots,k\}$, the $i\textsuperscript{th}$ tracer particle is in $A_i$ at time $T$,
\item and the initial mass vector of the tracer particles equals $\overline{n}$.
\end{itemize} 
The law $\mu_{\overline{n}}$ denotes the sub-probability measure given by the initial distribution of tracer particle locations  when the initial tracer particle $k$-vector equals $\overline{n}$. Note that the density at $(y_1,\cdots,y_k)$ of~$\mu_{\overline{n}}$  equals
$Z^{-k} \prod_{i=1}^k h_{n_i}(y_i)$.

To bound above the right-hand side of~(\ref{e.boundk}), we begin by noting that, for given $\overline{n} \in \N^k$, the expression
$$
 \int  \prod_{i=1}^k \nu_{y_i,2d(n_i)T}(A_i) \, \dd \mu_{\overline{n}}(\overline{y}) 
 $$
may be represented as the product of $k$ terms. The $i\textsuperscript{th}$ of these terms is the probability that a certain particle is present in $A_i$ at time $T$. This particle may or may not be born at time zero; it must be born if it is to appear in $A_i$ at time $T$. The density of the particle's time zero location equals $Z^{-1} h_{n_i}$, so that its probability of birth equals the single tracer particle's probability of being assigned mass $n_i$ initially. The particle then follows a Brownian trajectory at rate $2d(n_i)$. The $i\textsuperscript{th}$ term is at most $Z^{-1} \ell_{n_i} \mu(A_i)$, where recall that  
 $\ell_{n_i}$ is the supremum of $h_{n_i}$ (this due to a short argument that depends principally on the result that the heat equation decreases the supremum norm). Returning to~(\ref{e.boundk}), we find then that
$$
 \PP_N 
 \bigg( \bigcap_{1 \leq i \leq k} \big\{ x_i(T) \in A_i \big\} \bigg)
  \leq  
 \sum_{\bar{n} \in \N^k}  \prod_{i=1}^k \sup_{m \geq n_i}  \tfrac{n_i}{m}  \big( \tfrac{d(n_i)}{d(m)} \big)^{d/2} \cdot \prod_{i=1}^k Z^{-1} \ell_{n_i} \mu(A_i)
$$ 
whose right-hand side equals
$$
  Z^{-k}  \prod_{i=1}^{k} \mu(A_i) \, \cdot \bigg( \sum_{n=1}^\infty \sup_{m \geq n}  \tfrac{n}{m}  \big( \tfrac{d(n)}{d(m)} \big)^{d/2} \ell_n \bigg)^k 
   \, = \,  
 Z^{-k}   
   \Big( \sum_{n=1}^\infty n d(n)^{d/2} \ell_n \sup_{m \geq n}  m^{-1} d(m)^{-d/2}  \Big)^k  \prod_{i=1}^{k} \mu(A_i)  \, .
$$
That is, the $\PP_N$-probability of respective tracer particle occupation of the sets $A_i$, $1 \leq i \leq k$,
is at most $K^k  \prod_{i=1}^{k} \mu(A_i)$, as we sought to show. \qed

\subsection{Uniform integrability}\label{s.ui}
We now use Proposition~\ref{proppartplace} to prove Propositions~\ref{prop.firststep} and~\ref{prop.secondstep}.

\medskip

\noindent{\bf Proof of Propositions~\ref{prop.firststep} and~\ref{prop.secondstep}.}
We prove the results together, seeking a contradiction to the assumption that one or other proposition is false. We are thus supposing that, for some $m \in \N$, the density $f_m$ either fails to exist, or that it exists but has no finite supremum bound. Either way, we find that, for all $C > 0$,
there exist $\eps > 0$, $T > 0$ and $A \subseteq \R^d$ open such that, for some subsequence $\big\{ N_i: i \in \N \big\} \subseteq \N$,
\begin{equation}\label{e.nisubseq}
\PP_{N_i} \Big( \vert \chi_T^m \cap A \vert \geq C N_i \mu(A) \Big) \geq \eps \, , 
\end{equation}
where here we write $\chi_T^m$ for the set of locations of mass~$m$ particles at time~$T$.  In other words, if we call a particle in the $N_i$-indexed system that remains alive at time $T$, and then has mass $m$ and is located in $A$, a {\em target} particle,  then there is probability at least $\e$ that the proportion of the original particles that are target particles is at least $C \mu(A)$. We may express this event in terms of the $k$-vector of tracer particles from Proposition~\ref{proppartplace}. If $\vert \chi_T^m \cap A \vert$ equals a given  $\ell \in \N$, then the proportion of ordered $k$-vectors of initial particles all of whose elements are target particles is ${\ell \choose k} {N_i \choose k}^{-1}$. When $\ell$ equals $r \mu(A) N_i$ for $r$ bounded away from zero as $i \to \infty$, this last expression is asymptotic to $r^k\mu(A)^k$. We infer from~(\ref{e.nisubseq}) that, for any given $k \in \N$ and all sufficiently high~$i \in \N$,
$$
 \PP_{N_i}\bigg(  \bigcap_{i=1}^k \Big\{ x_i(T) \in A \Big\} \bigg) \geq \tfrac{1}{2} \e \cdot C^k \mu(A)^k \, ,
$$ 
where here we also used symmetry of the initial particle locations in the form that this probability is the same for the special $k$-sized set $\{1,\cdots,k\}$ as it is any given $k$-sized subset of $\{ 1,\cdots,N_i\}$.

Of course, Proposition~\ref{proppartplace} with each set $A_i$ set equal to $A$ provides an upper bound on the above probability, of $K^k \mu(A)^k$, 
where $K$'s value is specified in the proposition. Thus, $\e/2 \cdot  C^k$ is at most~$K^k$ for each $k \in \N$. The constant $C$ is seen to be at most $K$, contradicting our assumption.   \qed
\subsection{Deriving mass conservation}\label{s.aux}
Recall from Subsection~\ref{s.masscons} the notion that a solution of the Smoluchowski PDE conserves mass on $[0,T)$ for some $T \in [0,\infty]$.
In this short auxiliary section, we employ the analytic Lemma~\ref{lem3.1} concerning absence of particle concentration alongside a short futher analytic argument in order  
to show that, under certain conditions, there exists a weak solution of~(\ref{syspdeweak}) that conserves mass. 
\begin{proposition}\label{p.aside}
Suppose that $d: \N \to (0,\infty)$ is non-increasing, and that either
\begin{enumerate}
\item there exists $\alpha \in (0,1)$ and $\alpha_0 \in (0,\alpha)$
such that $\beta(n,m) \leq n^{\alpha_0} + m^{\alpha_0}$ for $n,m \in \N$, and  $d(m) \geq m^{-(1-\alpha)}$ for $m \in \N$;
\item or there exists $c > 0$ such that $\beta(n,m) \leq c(n+m)$ for $n,m \in \N$, and  $\inf_{n \in \N} d(n) > 0$.
\end{enumerate}

Then there exists a weak solution of~(\ref{syspdeweak}) that conserves mass on~$[0,\infty)$.
\end{proposition}
We will attempt the proof under for the second set of hypotheses and will omit some details.

\medskip

\noindent{\bf Sketch of proof.}
Since we assume that the second set of hypotheses hold, note that, by Lemma~\ref{lem3.1}, any weak solution $\big\{ f_n: n \in \N \big\}$ is such that $\sum_{n =1}^\infty n f_n(x,t)$ is bounded uniformly in space and time.
\begin{lemma}\label{l.lonebound}
Let $R(t,\cdot) : \R^d \to [0,\infty)$ be given by $R(t,\cdot) = \sum_{n=1}^\infty n^2 f_n(\cdot,t)$. Then, for each $t \geq 0$, $R(t) \in L^1(\R^d)$.
\end{lemma} 
\noindent{\bf Proof.}
Set $X:[0,\infty) \to [0,\infty]$, $X(t) = \int_{\R^d} R(t,x) \, \dd x$.
Then
\begin{equation}\label{e.xderiv}
 \frac{\dd X}{\dd t}(t)  =  2 \int_{\R^d} \sum_{n,m \geq 1} nm \beta(n,m) f_n(x,t) f_m(x,t) \, \dd x \, ;
\end{equation}
  one may interpret this equality by noting that $(n,m) \to n+m$ coagulation occurs at a rate equal to $\beta(n,m) f_n(x,t) f_m(x,t)$ at the space-time location~$(x,t)$, and the $X$ functional registers a change of $(n+m)^2 - n^2 - m^2 = 2nm$ as a result of any such collision.
 The right-hand side of~(\ref{e.xderiv}) is  at most
$$
 4   \int_{\R^d} \Big( \sum_n c n^2 f_n(x,t) \Big) \Big( \sum_m m f_m(x,t) \Big) \, \dd x  
  \leq  4c X(t)  \cdot \, \sup_{x \in \R^d}   \sum_m m f_m(x,t) \, \leq \, C X(t) \, , 
$$
where the first inequality (`is at most') employed $\beta(n,m) \leq c(n+m)$, and the third, uniform boundedness in space and time of $\sum_n n f_n$.

By Gronwall's lemma, $X(t) \leq c_1 \exp \big\{ c_2 t \big\}$ for all $t \geq 0$ and for some $c_1,c_2 > 0$, whence the result. \qed

To establish Proposition~\ref{p.aside}, it is enough to show that
\begin{equation}\label{e.enough}
 \lim_{N \to \infty} \, \frac{\dd}{\dd t} \int_{\R^d} \Big( \sum_{n=1}^N n f_n(x,t) \Big) \dd x = 0 \, .
\end{equation}
To derive this, note that
\begin{eqnarray*}
 \frac{\dd}{\dd t} \int_{\R^d}  \sum_{n=1}^N n f_n(x,t)
 & = & - \, 2 \sum_{n,m\in \N} {\bf 1}_{n \leq N < n + m} n \beta(n,m) f_n(x,t) f_m(x,t) \\
 & \geq & - \, 2c   \sum_{n,m\in \N} {\bf 1}_{n \leq N < n + m} n (n+m) 
f_n(x,t) f_m(x,t)  \, \geq \,  \Omega_1^N(x,t) + \Omega_2^N(x,t) \, , 
\end{eqnarray*}
where 
$$
 \Omega_1^N(x,t) =
 - \, 2c   \sum_{n,m\in \N} {\bf 1}_{n \leq N/2, m > N/2} \cdot \, n (n+m) 
f_n(x,t) f_m(x,t) 
$$
and
$$
 \Omega_2^N(x,t) =
 - \, 2c   \sum_{n,m\in \N} {\bf 1}_{n > N/2} \cdot \, n (n+m) f_n(x,t) f_m(x,t) \, . 
$$
We will derive~(\ref{e.enough}) by showing that $\int_{\R^d} \big\vert \Omega^N_i(x,t) \big\vert \, \dd x \to 0$ as $N \to \infty$, for $i \in \{1,2\}$.

To this end, note that
$$
 \big\vert \Omega^N_2(x,t) \big\vert \, \leq \, \sum_{n \geq N/2} n^2 f_n(x,t) \, \cdot \, \sum_{m \geq 1} f_m(x,t) \, \, + \, \,  \sum_{n \geq N/2} n f_n(x,t) \, \cdot \, \sum_{m \geq 1} m f_m(x,t) \, .
$$
Of the four sums on the right-hand side, the first (and thus the third) has $L_{\R^d}^1$-norm tending to zero by Lemma~\ref{l.lonebound} and the dominated convergence theorem, while, as we have noted, the fourth (and thus the second) is uniformly bounded in space and time.   

Turning to $\Omega_1$, we have that
$$
 \big\vert \Omega^N_1(x,t) \big\vert \, \leq \, \sum_{n \geq 1} n^2 f_n(x,t) \, \cdot \, \sum_{m \geq N/2} f_m(x,t) \, \, + \, \,  \sum_{n \geq 1} n f_n(x,t) \, \cdot \, \sum_{m \geq N/2} m f_m(x,t) \, .
$$
The first and third sums are in $L_{\R^d}^1$; the second and fourth converge to zero uniformly as $N \to \infty$. This establishes~(\ref{e.enough}) and completes our derivation of Proposition~\ref{p.aside}. \qed
 
\section{Review and summary}\label{s.final}
In this section, we list those instances where our derivation of Theorem~\ref{thmo} was incomplete, explaining where the proof is furnished in~\cite{HR3d}; discuss an imprecision in the proof of~\cite{HR3d} which we have sought to clarify by our presentation in this survey; and make a limited comparison between our present method of proof of key estimates by means of particle concentration bounds with the approach adopted in~\cite{HR3d}.
\subsection{The list of shortcuts in the survey's proofs}
In our proof of Theorem~\ref{thmo}, several steps are omitted. Beyond the absence of a proof of the classical Feynman-Kac formula, in the guise that $v(x,t)$ in Lemma~\ref{l.s} satisfies the PDE~(\ref{e.feynmankac}), the missing steps are: 
\begin{itemize}
 \item In the reduction of Proposition~\ref{propsz} to Proposition~\ref{propqqbar} undertaken in Section~\ref{s.coagprop}, the test function $J_n(x,t)$ was chosen to be equal to be identically one. The more general case requires only a few further lines of argument: see \cite[Section 3.5]{HR3d}.
 \item Proposition~\ref{p.mart}, showing smallness of the martingale $M(T)$ in~(\ref{emicro}), is not proved. See \cite[Section 5]{HR3d}.
 \item In the proof of Proposition~\ref{propqqbar} in Section~\ref{s.esterr}, the form of the test functions $J(x,n,t)$ and $\overline{J}(x,m,t)$ was simplified so that they have no space-time dependence. When this simplification is omitted, the action of the free motion operator in~(\ref{e.freem}) generates some extra terms, where one or both of the derivatives in the Laplacian fall on the test functions. The resulting terms tend to be smoother than the existing terms, and the methods of treating them are the same as for their rougher counterparts. See the start of \cite[Section 3]{HR3d} for the full scale version of~(\ref{e.freem}) and \cite[Section 3.3]{HR3d} for bounds on the terms appearing in that version.
 \item We have offered no proof of Proposition~\ref{p.errorbounds}(2), (3) and (4). Similarly to the previous point, this omission is a simplification of the presentation of the proof of the Stosszahlansatz. The key tools -- uniform killing probability bounds and particle concentration results -- apply equally to prove these statements as we saw that they did to prove Proposition~\ref{p.errorbounds}(1). See Sections~$3.3$ and~$3.4$ of \cite{HR3d} for the relevant bounds, valid under the original assumptions.
\end{itemize}
\subsection{A momentary spotlight on an obscurity}
In the opening paragraph of \cite[Section 4]{HR3d}, the step counterpart to that of the present Section~\ref{s.takingthelimit} is discussed: the low $\eps$ limit is taken of the approximate identity that is~(\ref{emicro}) in this survey. However, in the replacement of the collision term by its counterpart expressed using microscopic candidate densities, which happens by means of the Stosszahlansatz, \cite{HR3d} neglects to clarify that this replacement must be made simultaneously over the infinitely many mass pairs $\big\{ (n,m): m \in \N \big\}$, rather than merely being made for one such mass pair. What permits this simultaneous replacement is that the bound satisfied by the error ${\rm Err}_{n,m}$ in Proposition~\ref{propsz} contains a sum over $m \in \N$. As we have seen, the reason why we are able to prove the Stosszahlansatz with such an error bound is the uniform control on the killing probabilities $u_W$ that we saw in Subsection~\ref{s.killing}, specifically that Lemma~\ref{l.unifcont} is 
valid as the kernel $W$ varies over choices having given compact support. Of course, in the present survey, we have not presented a proof of all of the required estimates for Theorem~\ref{thmo}. A complete proof of the result, under the original assumptions made in~\cite{HR3d}, is formed by rendering \cite[Lemma~3.2]{HR3d} uniform over $(n,m) \in \N^2$; the changes needed to do this are contained in the proof of the present Lemma~\ref{l.unifcont}. 
\subsection{Comparison with later kinetic limit derivations of the Smoluchowski PDE}
Several ramifications of the statement and technique of proof of~\cite{HR3d} have been explored by Fraydoun Rezakhanlou, sometimes in collaboration with the author and others. We end by summarising the results so obtained and comparing the approaches to proof in these further articles, both with the original one in~\cite{HR3d} and with that expounded here.
\subsubsection{The particle concentration bound: its robustness and limitations}
The principal technical novelty presented in this survey is the use of the particle concentration bound 
Proposition~\ref{p.pl} to yield the key error bounds Proposition~\ref{p.errorbounds}: the technique used in~\cite{HR3d} was quite different. Our present technique requires stronger hypotheses, but when it may be applied, it yields strong conclusions about diverse aspects of particle dynamics. We now explain this summary by drawing a contrast with the method used in~\cite{HR3d}.

First, to expand, Proposition~\ref{p.pl} offers strong conclusions about the lack of build-up in particle concentration at positives times in the models $\PP_N$. However, it has content only under the fairly restrictive hypothesis that $d(m)$ decays no faster than $m^{-2/d}$, and, regarded as a tool to prove Proposition~\ref{p.errorbounds}, its use must be accompanied by the assumption that $\alpha(\cdot,\cdot)$ is bounded above uniformly. 
Here we make a comment about the one simple aspect of the quite different approach that was adopted in~\cite{HR3d} to prove the key estimates that correspond to the present Proposition~\ref{p.errorbounds}.
We now state a result giving an upper bound on the duration~$[0,T]$ total mean collision rate in the models $\PP_N$. It is \cite[Lemma 3.1]{HR3d}, which we call the ``bound on the collision''.
\begin{lemma}\label{l.boundcollision}
 For any $N \in \N$, and for  $T > 0$, 
 $$
  \eps^{d-2} \E_N \int_0^T \dd s \sum_{i,j \in I_q} \alpha(m_i,m_j) V_\eps(x_i - x_j) \leq Z \, .
 $$
\end{lemma}
\noindent{\bf Proof.} Let $X(T)$ denote the number of surviving particles in $\PP_N$ at time $t$. Consider the variant of~(\ref{emicro}) in which the term $J_n$ is replaced by $\mathbf{1}$ and the result summed over $n \in \N$, so that, for example, the first sum is a total particle count at time $T$. Taking expectations, we find that, for any $T > 0$,
\begin{equation}\label{e.zineq}
 \E_N X(T) = \E_N X(0) + \int_0^t \E_N \freem X(t) \dd t +  \int_0^t \E_N \collop X(t) \dd t \, .
\end{equation}
Particle count is conserved by free motion, so that  $\E_N \freem X(t) = 0$. On the other hand, the integrated mean collision rate 
 $\int_0^t \E_N \collop X(t) \dd t$ is equal to $\E_N \int_0^T \dd s \sum_{i,j \in I_q} \alpha(m_i,m_j) V_\eps(x_i - x_j)$. By~(\ref{e.zineq}), this quantity equals $\E_N X(0) - \E_N X(T)$ which is at most $\E_N X(0) = N$. Since $N = Z \eps^{2 - d}$ by~(\ref{e.intrange}), we obtain the result. \qed
 
 Of course, the proof is almost a triviality. However, it already highlights differences with, and the limitations of,  the particle concentration bound Proposition~\ref{p.pl}. Let us try to emulate this unprepossessing lemma's conclusion by using Proposition~\ref{p.pl}. By merely applying this result without using further tricks, the best we can do is the following.
 \begin{claim}\label{c.boundcollision}
 For any $N \in \N$, $K > 0$ and $T > 0$, we have that
 \begin{eqnarray*}
  & &  \eps^{d-2} \, \E_N \int_0^T \dd s \sum_{i,j \in I_q} \alpha(m_i,m_j) V_\eps(x_i - x_j) {\bf 1}_{\dist x_i(t) \dist \leq K} \\
  & \leq &   T \, \cdot \, (2K)^d \, \cdot \, \dist V \dist_\infty  \, \cdot \,  \sup_{n,m} \alpha(n,m) \, \cdot \, \left( \sup_{m \geq 1} m^{-1} d(m)^{-d/2} \right)^2 \cdot  \,  \left( \sum_{n=1}^\infty n d(n)^{d/2} \dist h_n \dist_\infty \right)^2  \, .
 \end{eqnarray*}
 \end{claim}
\noindent{\bf Proof.} The relation $N = Z \eps^{2-d}$ and 
Proposition~\ref{p.pl} applied with $k=2$ show that, for any time $t \in (0,\infty)$, 
\begin{eqnarray*}
 & & \sum_{i,j \in [1,N]} \PP_N \big( \dist x_j(t) - x_i(t) \dist  {\bf 1}_{\dist x_i(t) \dist \leq K} \leq \eps  \big) 
  \leq  N^2 \, \PP_N \big( \dist x_2(t) - x_1(t) \dist \leq \eps \, , \, \dist x_i(t) \dist \leq K \big) \\
 & \leq & Z^2 \eps^{2(2-d)} \, \cdot \, (2K)^d \eps^d \, \cdot \, Z^{-2}  \left( \sup_{m \geq 1} m^{-1} d(m)^{-d/2} \right)^2   \cdot \, \left( \sum_{n=1}^\infty n d(n)^{d/2} \dist h_n \dist_\infty \right)^2 \, .
\end{eqnarray*}

Note also that 
$$
\sum_{i,j \in [1,N]} \E_N \left[ \alpha\big( m_i(t),m_j(t)\big) V\big(\tfrac{x_i(t) - x_j(t)}{\eps}\big) \right] \leq  \sup_{n,m} \alpha(n,m) \, \cdot \,  \dist V \dist_\infty \, \cdot \,  \sum_{i,j \in [1,N]} \PP_N \big( \dist x_j - x_i \dist \leq \eps  \big)  \, .
$$

If we multiply the above left-hand side by $\eps^{d-2}$, and then further by $\eps^{-2}$ -- which we do because $V_\eps(\cdot) = \eps^{-2} V\big(\cdot/\eps \big)$ -- and integrate over $[0,T]$, then we obtain the left-hand side in the statement of the claim. Thus, the claim follows from the above two inequalities. \qed

Claim~\ref{c.boundcollision} falls short of Lemma~\ref{l.boundcollision} in a number of ways. One of these is a minor point: the use of the localization event $\dist x_i(t) \dist \leq K$, employed to permit the application of Proposition~\ref{p.pl}. Lemma~\ref{l.bigdisp} can easily be used to dispense with this detail, at the expense of an increase in the constant on the right-hand side in the claimed inequality. A more important shortcoming is the appearance of the
factor~$\sup_{n,m} \alpha(n,m)$ on this right-hand side. The reason that such a factor appears is because particle pair presence at distance of order $\eps$ at positive times is penalized due to the microscopic repulsion phenomenon which has been central to this survey, and which in particular we discussed in Section~\ref{s.torus}. Crudely, if $\alpha(n,m)$ is high, then the density for such presence is not of order $\eps^{d-2}$ but of order $\alpha(n,m)^{-1} \eps^{d-2}$. Proposition~\ref{p.pl} is not built to acknowledge this microscopic repulsion effect and the unwanted $\alpha$ goes uncancelled when the proposition is applied.

Lemma~\ref{l.boundcollision} experiences no such limitation. But of course its description of positive time particle distribution is limited to a very specific aspect of the overall dynamics; in comparison, Proposition~\ref{p.pl} is a robust tool, that will say something meaningful about any such aspect, when its hypotheses are satisfied.

Lemmas~3.2 and 3.3 of \cite{HR3d} form the counterpart to the present Proposition~\ref{p.pl} in the sense that they are tools used to prove the present Proposition~\ref{p.errorbounds}.
The two lemmas treat several aspects of particle dynamics other than the bound on the collision given in Lemma~\ref{l.boundcollision}. The proofs of these results generalize that of Lemma~\ref{l.boundcollision} in the sense that the mean value of variants of~(\ref{emicro}) are considered and their terms bounded. In the case of \cite[Lemma 3.3]{HR3d}, some of these terms involve sums over triples of particle indices, with one of the indices having no restriction on the mass parameter. In such cases, (\ref{hypo}) in the original assumptions is invoked during a proof by induction in order to find suitable bounds.

\subsubsection{The planar case: the route to the PDE}

In Section~\ref{s.planar}, we mentioned that, in \cite{HR2d}, the kinetic limit derivation counterpart to that of \cite{HR3d} was undertaken for dimension $d = 2$; we also reviewed the main changes to Theorem~\ref{thmo}'s statement in this case. The technique of proof is the same as in the original work, with the principal technical change concerning the particle distribution result 
\cite[Lemma 3.2]{HR3d}, where the proof may not be directly utilized because the non-negativity of the solution $H$ of Poisson's equation $- \Delta H = J$ (for a given non-negative $J$) is enjoyed in dimension $d \geq 3$ but not $d = 2$. We refer the reader to~\cite{HR2d} for further discussion of this technicality. However, we note that the robustness of the particle concentration bound Proposition~\ref{p.pl} has the virtue of permitting the extension of the proof of Theorem~\ref{thmo} developed in this survey to the case $d=2$, under the survey assumptions, without any comparable technical difficulty arising.

The kinetic limit derivations of \cite{HR3d} and \cite{HR2d} were extended to cases of variable radial dependence for particles in~\cite{RezakhanlouMPRF}; see the present Subsection~\ref{s.gelation}. The technique of proof, including the treatment of particle distribution bounds, is similar to that of the earlier works. 

\subsubsection{The case of continuous mass and a different approach to particle concentration}

As we mentioned after the equations in Section~\ref{s.smolcoag}, the Smoluchowski coagulation-diffusion PDE has a continuous counterpart, where the mass parameter is now a non-negative real. The kinetic limit derivation of the PDE is revisited in~\cite{YRH} in the case that $d \geq 3$ for the PDE with continuous mass parameter. The principal innovation of the article is~\cite[Theorem 3.1]{YRH}, a tool for proving particle concentration which is novel in comparison with that of~\cite{HR3d} and~\cite{HR2d}. This technique has distinct similarities with our tracer particle proof of Proposition~\ref{p.pl}: the comparison of diffusion rate dependent terms in~\cite[(3.6)]{YRH} is a rough counterpart to~(\ref{e.normalbound}).



\subsubsection{The planar case with fragmentation: equilibrium fluctuations}

The Smoluchowski PDE may be modified to include interaction terms corresponding to pairwise particle fragmentation. In the models~$\PP_N$, we may model this fragmentation effect by declaring that a particle of given mass is subject to fragment at the ring times of a Poisson clock that ticks at a mass-dependent rate. On fragmenting, the particle splits in two. The detailed rule for this splitting may be chosen to be a ``dual'' of the rule specified for coagulation under the heading ``the precise mechanism for collision'' in Section~\ref{sec.micromodels}: the fragmenting particle retains its location, and some random proportion of its precollisional mass, while a new particle, bearing the residue of that mass, appears in a randomly selected microscopic vicinity of the fragmenting particle's location. 

Of course, one may attempt to carry out a kinetic limit derivation of the Smoluchowski coagulation-fragmentation-diffusion PDE from microscopic models $\PP_N$ that have been altered in this manner.  Such results have yet to be proved, though important elements for proofs are suggested in~\cite{RanjbarRezak}. In this article, Ranjbar and Rezakhanlou studied a different aspect of particle dynamics in the case that dimension $d$ equals two. The assertion that macroscopic particle densities adhere to a solution of the Smoluchowski PDE is in a sense a weak law of large numbers. What of the analogue of the central limit theorem, a result describing the typical fluctuation of particle density statistics in high indexed $\PP_N$ from the density profile offered by the PDE solution? In~\cite{RanjbarRezak}, the authors define empirical fluctuation fields under $\PP_N$, modelling the discrepancy as a function of space-time of the empirical density of particles from the prediction made by the PDE solution, normalized by a square root of total particle number, in the style of the 
central limit theorem.  In \cite[Conjecture 2.1]{RanjbarRezak}, it is conjectured that, in a high $N$ limit, the fluctuation field converges to a random limit that solves an Ornstein-Uhlenbeck equation under which the density profile diffuses freely and is subject to coagulative and fragmentative forces specified by a linearization of those present in the PDE, as well as to a space-time dependent white noise stimulus determined by the PDE solution. 

Proving this conjecture is likely to be a demanding task, probably much more difficult than that of adapting existing techniques to carry out the kinetic limit derivation of an analogue of Theorem~\ref{thmo}. 
(The reason for the added difficulty is essentially that the rate of decay to zero of the various error terms in the microscopic counterpart to the weak PDE solution must be shown to converge to zero at a sufficiently fast rate.) 
Despite this degree of difficulty, the authors of~\cite{RanjbarRezak} advance a case for the conjecture by rigorously analysing the system at equilibrium. The modified microscopic models~$\PP_N$ have mechanisms for both coagulation of pairs of particles, and fragmentation of particles into pairs; and these mechanisms have been chosen so that the film of the coagulation event when played in reverse shows the fragmentation event. Thus, the equilibrium measures of the laws $\PP_N$ may be explicitly identified: under them, the distribution of particles of any given mass is simply a Poisson process (of some constant intensity determined by the mass), with the clouds of particles of distinct masses being independent. In \cite[Theorem 3.1]{RanjbarRezak}, the conjecture of convergence to the Ornstein-Uhlenbeck process mentioned above is proved for the system at these equilibria. The proof requires 
an understanding of the relation between microscopic and macroscopic interaction propensities which extends~(\ref{recptwod}) to treat fragmentation but also involves an unexpected interaction with the free motion dynamics. 

\bibliographystyle{plain}
\bibliography{biblion}
\end{document}